\documentclass[pdflatex,sn-mathphys-num]{sn-jnl}

\usepackage{graphicx}%
\usepackage{multirow}%
\usepackage{amsmath,amssymb,amsfonts}%
\usepackage{amsthm}%
\usepackage{mathrsfs}%
\usepackage[title]{appendix}%
\usepackage{xcolor}%
\usepackage{textcomp}%
\usepackage{manyfoot}%
\usepackage{booktabs}%
\usepackage{algorithm}%
\usepackage{algorithmicx}%
\usepackage{algpseudocode}%
\usepackage{listings}%
\usepackage{placeins}

\usepackage{subfigure} 
\usepackage{upgreek} 
\usepackage{algpseudocode}

\theoremstyle{thmstyleone}%
\newtheorem{theorem}{Theorem}
\theoremstyle{thmstyletwo}%

\theoremstyle{thmstylethree}%
\newtheorem{definition}{Definition}%

\newtheorem{problem}{Problem}

\newcommand{\laplacian}{\Delta}
\newcommand{\esssup}{\text{ess\,sup}} 
\newcommand{\uvec}{\mathbf{u}}
\newcommand{\vvec}{\mathbf{v}}
\newcommand{\tuBT}{\widetilde{\uvec}_{0;B,T}}
\newcommand{\uBT}{\uvec_{0;B,T}}
\newcommand{\uST}{\uvec_{0;S,T}}
\newcommand{\argmax}{\operatorname{argmax}}
\newcommand{\argmin}{\operatorname{argmin}}

\newcommand{\tuE}{\widetilde{\mathbf{u}}_{\E_0}}

\newcommand{\tuET}{\widetilde{\uvec}_{0;\E_0,T}}

\def\k{{\mathbf k}}
\def\f{{\mathbf f}}
\def\g{{\mathbf g}}
\def\u{{\mathbf u}}
\def\v{{\mathbf v}}
\def\x{{\mathbf x}}
\def\z{{\mathbf z}}

\def\E{{\mathcal{E}}}

\def\K{{\mathcal{K}}}
\def\M{{\mathcal{M}}}
\def\N{{\mathcal{N}}}
\def\Q{{\mathcal{Q}}}
\def\R{{\mathcal{R}}}
\def\T{{\mathcal{T}}}
\def\OO{{\mathcal{O}}}

\def\0{{\mathbf 0}}
\def\CC{{\mathbb{C}}}
\def\P{{\mathcal{P}}}
\def\RR{{\mathbb{R}}}
\def\TT{{\mathbb{T}}}
\def\ZZ{{\mathbb{Z}}}

\def\bnabla{\boldsymbol{\nabla}}
\def\bomega{\boldsymbol{\omega}}

\begin{document}

\title[The Ladyzhenskaya-Prodi-Serrin 
	Conditions and the Search for Extreme Behavior in 3D Navier-Stokes 
	Flows]{The Ladyzhenskaya-Prodi-Serrin Conditions and the Search for 
	Extreme Behavior in 3D Navier-Stokes Flows}


\author[1]{\fnm{Elkin} \sur{Ram\'irez}}

\author*[1]{\fnm{Bartosz} \sur{Protas}}\email{bprotas@mcmaster.ca}


\affil[1]{\orgdiv{Department of Mathematics and Statistics}, \orgname{McMaster University},  \orgaddress{\street{1280 Main Street West}, \city{Hamilton}, \postcode{L8S 4L8}, \state{Ontario}, \country{Canada}}}

%


\abstract{
\normalfont\unboldmath

In this investigation, we conduct a systematic computational search for 
potential singularities in 3D Navier-Stokes flows on a periodic domain 
$\Omega$ based on the Ladyzhenskaya-Prodi-Serrin conditions. They 
assert that for a solution $\u(t)$ of the Navier-Stokes system to be 
regular on an interval $[0,T]$, the integral $\int_{0}^T 
\|\u(t)\|_{L^q}^p\,dt$, where $2/p+3/q=1,\;q>3$, and the 
expression $\sup_{t \in [0,T]} \|\u(t)\|_{L^3}$ must be 
bounded. Flows which might become singular and violate these conditions 
are sought by solving a family of variational PDE optimization problems 
where we identify initial conditions $\u_{0}$ with the corresponding 
flows $\u(t)$ locally maximizing the integral $\int_{0}^T 
\|\u(t)\|_{L^q}^p\,dt$ for a range of different values of $q$ 
and $p$ or the norm $\|\u(T)\|_{L^3}$ for different time 
windows $T$ and increasing sizes $\| \u_0 \|_{L^q}$ of the 
initial data. We consider two formulations where these expressions are 
maximized over appropriate Lebesgue spaces $L^q(\Omega)$ or the largest 
Hilbert-Sobolev spaces $H^s(\Omega)$ embedded in them. The lack of 
Hilbert-space structure in the first case necessitates development of a 
novel computational approach to solve the problem. While no evidence of 
unbounded growth of the quantities of interest, and hence also for 
singularity formation, was detected, we were able to quantify how 
``close" the flows realizing such worst-case scenarios come to forming a 
singularity. A comparison of these results with estimates on the rate 
of growth of the norms $||\u(t)||_{L^q}$ and of the enstrophy $\E(t)$ 
indicates that the extreme flows do enter a regime where these 
quantities are amplified at a rate consistent with singularity 
formation in finite time, but this growth is not sustained long enough 
for singularities to form. This transient nonlinear amplification is 
also found to increase with $q$.

}




\maketitle


\section{Introduction}
\label{sec:into}

We consider the incompressible Navier-Stokes system defined on the 
three-dimensional (3D) torus $\Omega=\TT^{3}:= \RR^{3}/\ZZ^{3},$ where 
``:='' means ``equal to by definition'', and 
 equipped with periodic boundary conditions 
\begin{subequations}\label{eq:NS}
\begin{alignat}{2}
\partial_t\u + \left(\u\cdot\bnabla_{\x}\right)\u + {\bnabla_{\x}}\, p - \nu\laplacian\u & = 0 & &\qquad\mbox{in} \,\,\Omega\times(0,T], \label{eq:NSa} \\
\bnabla_{\x}\cdot\u & = 0 & & \qquad\mbox{in} \,\,[0,T]\times\Omega, \\
\u(0) & = \u_0, &   &\label{eq:NSc}
\end{alignat}
\end{subequations}
where $\u = \u(t,\x) = {[u_1, u_2, u_3]^T}$ is the velocity field, $p$ 
the pressure, $\nu>0$ the coefficient of kinematic viscosity, 
$\bnabla_{\x}$ denotes the gradient with respect to the space variable 
$\x=[x_1,x_2,x_3]^{T}$ (this notation is needed to distinguish it from 
another notion of a gradient that will be introduced later) and $\u_0$ 
is the initial condition. The velocity gradient $\bnabla_{\x}\,\u$ is a 
tensor with components $[\bnabla_{\x}\,\u]_{ij} = \partial_j u_i$, 
$i,j=1,2,3$.  For simplicity when there is no risk of 
confusion we will denote $\u(t) = \u(t,\cdot)$.

The study of the existence, uniqueness and regularity of solutions of 
the Navier-Stokes system \eqref{eq:NS} dates back to the beginning of 
the 20th century. Despite this long effort, the question about the 
global-in-time existence of classical solutions of \eqref{eq:NS} 
corresponding to smooth initial data $\u_0$ still remains open and has 
been recognized as one of its ``Millennium Problems" by the Clay 
Mathematics Institute \cite{f00} (this problem involves flows either on 
the unbounded or periodic domain, $\RR^3$ or $\TT^{3}$). On the other 
hand, classical solutions of \eqref{eq:NS} are known to exist globally 
in time in two dimensions (2D) \cite{dg95}.  We add that the question 
of the global existence of classical solutions of the inviscid Euler 
system, obtained by setting $\nu = 0$ in \eqref{eq:NS}, and 
corresponding to smooth initial data $\u_{0}$ is also open in 3D 
\cite{mb02,gbk08}. Major progress concerning the existence of weak 
solutions on $\RR^3$ was made by Leray in 1934 \cite{l34}. Later, in 
1951, Hopf \cite{h51} established global existence of weak solutions, 
albeit without uniqueness, on bounded domains. The question of the 
uniqueness of Leray-Hopf weak solutions remains open and in 2022, 
Albitron, Bru\'e \& Colombo \cite{abc22} constructed nonunique 
Leray-Hopf weak solutions with a time-dependent force and borderline 
regularity. 

Classical solutions may cease to exist as a result of a spontaneous 
formation of a singularity, such that system \eqref{eq:NS} may no 
longer be satisfied in a pointwise sense. One approach to study the 
existence of classical solutions relies on the so-called ``conditional 
regularity results". These are global conditions that need to be 
satisfied by the weak Leray-Hopf solutions of \eqref{eq:NS} in order  
for these solutions to also satisfy this system in the classical sense. 
The best known results of this type are the enstrophy conditions 
\cite{ft89} and the family of the Ladyzhenskaya-Prodi-Serrin conditions 
\cite{KisLad57,Prodi1959,Serrin1962,Escauriaza2003}, both of this will 
be considered in this study.

While the question about singularity formation in the Navier-Stokes and 
Euler flows is one in the mathematical analysis of partial differential 
equations (PDEs), efforts to shed light on this possibility based on 
numerical computations have had a long history 
\cite{bmonmu83,ps90,b91,k93,p01,bk08,oc08,ghdg08, 
gbk08,h09,opc12,bb12,k13,opmc14,CampolinaMailybaev2018}. While most 
  of these investigations were rather inconclusive, the studies 
  \cite{lh14a,lh14b,Hou:22:Euler} stand out providing strong evidence 
  for finite-time singularity formation in Euler flows on 3D bounded 
  axisymmetric domains. In the context of Navier-Stokes flow we 
  highlight \cite{Hou:22:NS} which also suggested singularity formation 
  in a flow in the same geometry. All of these studies relied on 
  ad-hoc, although motivated by careful physical arguments, or 
  numerical experiments, choice of the initial data $\u_{0}$ that could 
  lead to singularity formation. There is thus a need to find initial 
  data that could trigger potential singularity formation in a 
  systematic manner.
  
In the present investigation we follow a fundamentally different 
approach, initiated in the seminar study \cite{ld08}, where extreme, 
possibly singular flows are sought methodically as solutions of 
suitable variational optimization problems. More specifically, this 
allows one to find initial data such that the resulting flows maximize 
the growth of a chosen quantity measuring the regularity of the 
solution related to one of the conditional regularity results and whose 
unbounded growth in finite time would signal singularity formation. Our 
earlier investigations based on this idea focused on maximizing the 
quantities involved in the enstrophy condition and in one member of the 
Ladyzhenskaya-Prodi-Serrin family of conditional regularity results 
\cite{KangYunProtas2020,KangProtas2021}. These studies provided no 
evidence for an unbounded growth of these quantities which, however, 
did show significant transient growth exhibiting certain universal 
behavior. In the present study we extend the results from 
\cite{KangProtas2021} by considering a wide range of the 
Ladyzhenskaya-Prodi-Serrin conditions while also overcoming some 
technical limitations that were encountered in \cite{KangProtas2021}. 
This optimization-based framework has also proven useful in the study 
of globally well-posed systems, such as the one-dimensional (1D) 
viscous Burgers and the 2D Navier-Stokes systems 
\cite{ap11a,ap13a,ayala_doering_simon_2018,Yun2018}, where rigorous a 
priori bounds are available on the growth of various energy-type 
quantities. By finding flows extremizing the growth of these quantities 
it is possible to verify the sharpness of these a priori estimates. In 
the context of 1D viscous Burgers flows, the results from \cite{ap11a} 
inspired efforts to find improved rigorous upper bounds on the growth 
of enstrophy and the gap was closed only recently 
\cite{AlbrittonDeNitti2023}. A survey of this research program is 
offered in the review paper \cite{p21a}.

The main contribution of the present study is to construct, by 
numerically solving suitable optimization problems, families of initial 
conditions such that the resulting Navier-Stokes flows locally maximize 
the quantities in the different Ladyzhenskaya-Prodi-Serrin regularity 
criteria. Although we find no evidence for singularity formation, by 
recalling known estimates on the rate of growth of the corresponding 
solution norms, we are able to assess how ``close" these extreme flows 
come to producing a singularity. The structure of the paper is as 
follows: after introducing key definitions in Section 
\ref{sec:KeyQuant}, in Section \ref{sec:RegulResults} we recall the 
enstrophy and the Ladyzhenskaya-Prodi-Serrin conditions together with 
the associated inequalities; in Section \ref{sec:OptFormulations} we 
state the optimization problems and the solution approach is described 
in Section \ref{sec:SolApr}; the computational results are presented in 
Section \ref{sec:results} whereas final conclusions are deferred to 
Section \ref{sec:final}. Some more technical results are collected in 
three appendices.

\section{Definitions of Key Quantities}
\label{sec:KeyQuant}

In this section we collect definitions of the key quantities. While 
they are quite standard, they are given here for completeness.

\begin{definition}[Lebesgue space]
We denote 
$$L^q(\Omega):=\left\{{\bf{f}}:\Omega\rightarrow\RR^{3}: {\bf{f}}\mbox{ is measurable and }
\|{\bf{f}}\|_{L^q}  <\infty\right\}$$
where $q\in[1,\,\infty]$ and the norm is defined as
\begin{subequations}
\begin{alignat}{2}
\|{\bf{f}}\|_{L^q} := & 
\left( \int_\Omega |{\bf{f}}(\x)|^q \,d\x \right)^{\frac{1}{q}}, & \qquad & 1\leq q <\infty,
\label{eq:normLq} \\
\|{\bf{f}}\|_{L^q}:= & \, \esssup_{\x\in\Omega}\; |{\bf{f}}(\x)| && \nonumber \\
 = & \, \inf \left\{a\in\RR : \mu\left(\left\{\x\in\Omega: |{\bf{f}}(\x)|\geq a\right\}\right)=0\right\},  &  & q =\infty \label{eq:normL3}
\end{alignat}
\end{subequations}
in which $\mu$ denotes the Lebesgue measure.
\end{definition}

\begin{definition}
For $s \geq 0,$ the Sobolev space $H^s(\Omega)$ is defined as
\begin{equation*}
H^{s}\left(\Omega\right):=\left\{{\bf{f}}\in L^2\left(\Omega\right):\|{\bf{f}}\|^{2}_{H^s}
<\infty\right\},
\end{equation*}
where the norm is defined in terms of the Fourier series of $\f$ with
the Fourier coefficients
\begin{equation}
\label{eq:fc}
{\widehat{{\bf{f}}}_{\k} := \int_{\Omega} e^{-2\pi i\k\cdot\x} \, {\bf{f}}(\x)\,d\x}
\end{equation}
as
\begin{equation}
\|{\bf{f}}\|^{2}_{H^s}:=\sum_{\k\in\ZZ^3}{\left[1+(2\pi k)^{2}\right]}^s |\widehat{{\bf{f}}}_{\k}|^2,
\label{normHs}
\end{equation}  
where $k:=|\k|$. The corresponding inner product is denoted $\langle\cdot,\cdot\rangle_{H^s}$, such that $\langle\f,\f\rangle_{H^s} = \| \f \|^{2}_{H^s}$.
\end{definition}

\begin{definition}
The kinetic energy of a velocity field flow $\u(t)$ is defined as
\begin{equation}
\K(\u(t))  :=  \frac{1}{2}\int_\Omega | \u(t,\x) |^2 \,d\x=\frac{1}{2}\|\u(t)\|_{L^2}^2.
\label{eq:K}
\end{equation}
\end{definition}

\begin{definition}
The enstrophy\footnote{The enstrophy is often defined
    without the prefactor of 1/2. However, for consistency with earlier
    studies belonging to this research program
    \cite{ap11a,ap13a,ap13b,ap16,Yun2018,KangYunProtas2020}, we choose
    to retain this factor here.}
of a divergence-free velocity field $\u(t)$ is defined as
\begin{equation}
\E(\u(t))  :=  \frac{1}{2}\int_\Omega | \boldsymbol{\omega}(t,\x) |^2 \,d\x=\frac{1}{2}\|\boldsymbol{\omega(t)}\|_{L^{2}} = \frac{1}{2}\|\bnabla_{\x}\,\u(t)\|_{L^2},
\label{eq:E1}
\end{equation}
where $\boldsymbol{\omega}(t,\x):=\bnabla_{\x}\times\u(t,\x)$ is the vorticity of $\u$ and the last equality is a consequence of the incompressibility of the field $\u(t)$ and the periodic boundary conditions \cite{dg95}. 
\end{definition}

\noindent For simplicity we will denote $\K(t) := \K(\u(t))$ and $\E(t) 
:= \E(\u(t))$ for $t > 0$, and $\K_{0} := \K(\u(0))$ and $\E_{0} := 
\E(\u(0))$. We also add that for smooth solutions of the Navier-Stokes 
system $\K(t)$ and $\E(t)$ are related via the energy equation
\begin{equation}
\frac{d\K}{dt}=-2\nu\E.
\label{eq:kinetic_enstrophy}
\end{equation}

\section{Conditional Regularity Results}
\label{sec:RegulResults}

Conditional regularity results are statements about Leray-Hopf weak 
solutions $\u(t)$ of the Navier-Stokes system expressing verifiable 
conditions that need to be fulfilled by these solutions if they are to 
also satisfy system \eqref{eq:NS} in the classical (pointwise) sense 
\cite{RobinsonRodrigoSadowski2016}. While in general it is not known if 
they hold, they represent necessary and sufficient conditions for 
$\u(t)$ to be smooth on the time interval of interest $[0,T]$. Here we 
consider two such conditions, namely, the enstrophy and the family of 
the Ladyzhenskaya-Prodi-Serrin conditions together with associated 
inequalities for the rate of growth of the quantities they involve.

\subsection{The Enstrophy Condition}
\label{sec:E}

Leray-Hopf weak solutions of the Navier-Stokes system \eqref{eq:NS} 
defined in the interval $[0,T]$ satisfy the inequality, 
cf.~\eqref{eq:kinetic_enstrophy},
\begin{equation}
\int_{0}^{T} \E(\u(t))\,dt<\infty.{}
\end{equation}
In order for such solutions to be strong (classical), they must also 
satisfy the sufficient and necessary condition \cite{ft89}
\begin{equation}
\label{eq:RegCrit_FoiasTemam}
\mathop{\sup}_{0 \leq t \leq T} \E(\u(t))  < \infty.
\end{equation}
Thus, the boundedness of enstrophy guarantees that the 
Navier-Stokes flow remains smooth and should a singularity form at 
$T^*$, then necessarily 
\begin{equation}
\lim_{t\rightarrow {T^{*}}^{-}}\E(\u(t))=\infty.
\label{eq:Ensblowup}
\end{equation}

\subsection{Bounds on the Rate of Growth of the Enstrophy}
\label{sec:dEdt}

Whether or nor blow-up \eqref{eq:Ensblowup} can occur is dictated by 
the maximum sustained rate at which enstrophy can be amplified in 
classical solutions of the Navier-Stokes system. This problem was 
considered in \cite{ld08}, where the following upper bound was obtained
\begin{equation}
\frac{d\E}{dt}\leq\frac{27}{8\pi^4 \nu^3}\E^3.
\label{eq:RateGrowthEbound}
\end{equation}
An analogous estimate involving a norm of the symmetric part of the 
velocity gradient $\bnabla_{\x}\u$ was obtained in \cite{Miller2020}. 
The numerical evidence provided in \cite{ld08,ap16} suggests that 
estimate \eqref{eq:RateGrowthEbound} is in fact sharp, in the sense 
that there exist families of incompressible vector fields $\tuE$ with 
the enstrophy $\E(\tuE) = \E_{0}$ such that $(d/dt)\E(\tuE) = 
\OO(\E_{0}^{3})$ as $\E_{0} \rightarrow \infty$. Direct integration of 
\eqref{eq:RateGrowthEbound} leads to 
\begin{equation}
\E(\u(t))\leq\frac{\E_{0}}{\sqrt{1-\frac{27}{4\pi^4\nu^3}\E_{0}^2\,t}},
\label{eq:boundE}
\end{equation}
where the right-hand side (RHS) becomes unbounded as $t \rightarrow 
T^{*}=4\pi^4 \nu^{3}/(27\E_{0}^2)$. Thus, together with 
\eqref{eq:RegCrit_FoiasTemam}, relation \eqref{eq:boundE} is equivalent 
to a local existence result.

In order for a singularity to form in finite time, the enstrophy must 
be amplified at a sustained rate that is sufficiently rapid. Assuming 
that 
\begin{equation}
\frac{d\E}{dt} < C\E^{\alpha}, \qquad C > 0, \quad 0 < \alpha \le 3
\label{eq:dEdt}
\end{equation}
(hereafter $C$ will denote a generic positive constant whose 
numerical value may change from instance to instance) and using  
Gr\"onwall's lemma together with the energy equality 
\eqref{eq:kinetic_enstrophy} in an integrated form
\begin{equation}
\int_{0}^{t} \E(s)\,ds={\frac{1}{2\nu}\left[\K_0-\K(\u(t))\right]}\leq\frac{1}{2\nu}\K_0,
\label{eq:boundEwithK}
\end{equation}
yields the upper bound 
$\E(t)\leq\E_{0}\,\exp\left[\frac{\K_0}{2\nu}\right]$ valid when 
$\alpha \in (0,2]$ in \eqref{eq:dEdt}. Therefore, in order for a 
singularity to form the enstrophy of the flow must be amplified at a 
sustained rate proportional to $\E^{\alpha}$ with $\alpha \in (2,3]$. 
This observation offers a convenient way to assess how close various 
extreme flows come to producing a singularity. We add that when the 
states $\tuE$ saturating the upper bound in \eqref{eq:RateGrowthEbound} 
are used as the initial data $\u_{0}$ in \eqref{eq:NS}, the growth of 
enstrophy is immediately depleted in the resulting Navier-Stokes flows 
\cite{ap16}.

\subsection{Ladyzhenskaya-Prodi-Serrin Conditions}
\label{sec:LPS}

Another well-known conditional regularity result is given by the family 
of the Ladyzhenskaya-Prodi-Serrin (LPS) conditions 
\cite{KisLad57,Prodi1959,Serrin1962}. They assert that $\u(t)$ is a 
smooth classical solution of the Navier-Stokes system \eqref{eq:NS} on 
the time interval $[0,T]$ if and only if 
\begin{align}
& \u \in L^p([0,T];L^q(\Omega)), \qquad 2/p+3/q \le 1, \quad q > 3,   \label{eq:LPS} \\
& \u \in L^{\infty}([0,T];L^3(\Omega)), \label{eq:LPS3}
\end{align}
with the latter condition corresponding to the limiting case with $q = 
3$ \cite{Escauriaza2003}. Another blow-up criterion related to 
\eqref{eq:LPS3} was recently established in \cite{Tao2020}. As regards 
the values of the indices $p$ and $q$ in \eqref{eq:LPS}, we are 
interested in the borderline (equality) case with $p = 2q / (q - 3)$. 
This condition implies that if a singularity is formed in a classical 
solution $\u(t)$ of the Navier-Stokes system \eqref{eq:NS} at the time 
$T^{*}\in(0,  \infty)$, then necessarily
\begin{equation}
\lim_{t \rightarrow T^{*}} \int_0^t \| \u(\tau) \|_{L^q}^p \, d\tau  = \infty, \qquad 2/p+3/q = 1, \quad q > 3.
\label{eq:LPSblowup}
\end{equation}
Leray-Hopf weak solutions are known to satisfy the following estimates 
\cite{Constantin1991,KangProtas2021}
\begin{subequations}
\begin{alignat}{2}
& \int_{0}^{T}\|\u(t)\|_{L^q}^{\frac{4q}{3(q-2)}}\,dt\leq C\K_{0}^{\frac{2q}{3(q-2)}},&  \qquad  & 2\leq q\leq6,
\label{eq:ub1} \\
& \int_{0}^{T}\|\u(t)\|_{L^q}^{\frac{q}{q-3}}\,dt\leq C\K_{0}^{3},& \quad & q>6.
\label{eq:ub2}
\end{alignat}
\end{subequations}
with exponents in the integrant expressions smaller than $p = 2q / (q - 
3)$ in \eqref{eq:LPS}. Thus, if a singularity is to form in a 
Navier-Stokes flow, the norm $\| \u(t) \|_{L^q}$ must grow at 
such a rate that the integral in \eqref{eq:LPSblowup} becomes unbounded 
as $t \rightarrow T^{*}$ while the integrals in 
\eqref{eq:ub1}--\eqref{eq:ub2} remain finite. This issue is considered 
in more detail in the subsection below.

A preliminary attempt to search for singular Navier-Stokes flows that 
would realize \eqref{eq:LPSblowup} based on an optimization formulation 
was conducted in \cite{KangProtas2021} where the values $q = 4$ and $p 
= 8$ were used. The goal of the present investigation is to complete 
this program by considering a wider range of parameters $q$ and $p$ 
in addition to overcoming certain technical limitations encountered in 
\cite{KangProtas2021}.

\subsection{A Priori Bounds on $\frac{d}{dt}\|\u(t)\|_{L^q}$}
\label{sec:dLPSdt}

Here we analyze conditions on the rate of growth of 
$\|\u(t)\|_{L^q}$ necessary for blow-up \eqref{eq:LPSblowup}. 
A well-known upper bound on the rate of growth of this quantity is 
\cite{RobinsonSadowski2014}
\begin{equation}
\label{eq:aprioriLq}
\frac{d}{dt}\|\u(t)\|_{L^q}\leq C\|\u(t)\|_{L^q}^{\frac{3(q-1)}{q-3}}{=C\|\u(t)\|_{L^q}^{p+1}}, \quad q>3,
\end{equation}
where $C = C(q) > 0$. By integrating this inequality with respect 
to time, we obtain
\begin{equation}
\|\u(t)\|_{L^q}\leq \frac{1}{\left(\|\u_{0}\|_{L^q}^{-p}-p\,C\,t\right)^{1/p}}
\end{equation}
which becomes unbounded as $t \rightarrow T^{*} = 
1 / (p\,C\,\|\u_{0}\|_{L^q}^{p})$. Together with 
\eqref{eq:LPS}, this relation thus represents a statement of a local 
existence result for smooth solutions, akin to \eqref{eq:boundE}. To 
the best of our knowledge, estimates analogous to \eqref{eq:aprioriLq} 
are not available in the limiting case with $q = 3$ where the ``best" 
upper bound is \cite{RobinsonRodrigoSadowski2016}
\begin{equation}
\frac{1}{3}\frac{d}{dt}\|\u(t)\|_{L^3}^3\leq C\|\u(t)\|_{L^3}^2\|\u(t)\|_{L^9}^{3},
\label{eq:du3dt}
\end{equation}
which however does not have the form $\frac{d}{dt}Y\leq CY^{\alpha}$. 
The question of the sharpness of bounds \eqref{eq:aprioriLq} was 
considered in \cite{BleitnerProtas2026} where numerical evidence is 
provided suggesting these estimates are indeed sharp (up to 
prefactors) for $q > 3$.

\begin{figure}
\centering
\includegraphics[width=0.5\textwidth]{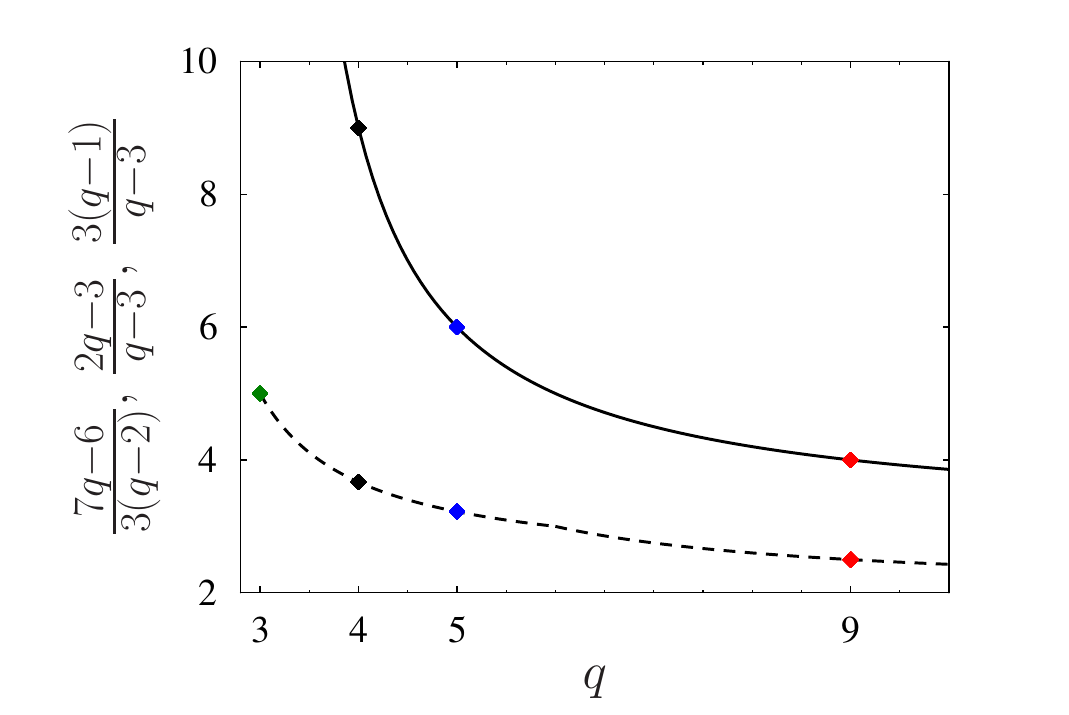}
\caption{Dependence of the exponents in the upper bounds  
in inequalities \eqref{eq:aprioriLq} (solid line) and \eqref{eq:noblowup_ub}
(dashed line) on $q$. The symbols represent the values of parameter $q$ considered 
in the present study.}
\label{fig:rateofgrowthLq}
\end{figure}

Estimates \eqref{eq:ub1}--\eqref{eq:ub2} can be used in conjunction 
with Gr\"onwall's lemma  to determine the largest rate of growth 
of the norm $\|\u(t)\|_{L^q}$ ensuring there is no blow-up in 
finite time, which gives
\begin{equation}
\frac{d}{dt}\|\u(t)\|_{L^q}\leq
\begin{cases}
C\|\u(t)\|_{L^q}^{\frac{7q-6}{3(q-2)}}, \qquad & 2\leq q \leq6, \\
\\
C\|\u(t)\|_{L^q}^{\frac{2q-3}{q-3}}, & q>6. \\   
\end{cases} 
\label{eq:noblowup_ub}
\end{equation}
To show this, let us first consider the case with $2\leq q \leq6$, where
\begin{equation*}
\frac{d}{dt}\|\u(t)\|_{L^q}\leq
C\|\u(t)\|_{L^q}^{\frac{7q-6}{3(q-2)}}=C\|\u(t)\|_{L^q}\|\u(t)\|_{L^q}^{\frac{4q}{3(q-2)}}.
\end{equation*}
By Gr\"onwall's lemma, we can then bound $\|\u(t)\|_{L^q}$ on any interval $[0,T]$, $T>0,$ by
\begin{equation*}
\begin{aligned}
\|\u(t)\|_{L^q} & \leq  \|\u_{0}\|_{L^q}\,\mbox{exp}\left(C\int_{0}^{T}\|\u(t)\|_{L^q}^{\frac{4q}{3(q-2)}}dt\right)\\
& \leq\|\u_{0}\|_{L^q}\mbox{exp}\left(C\K_{0}^{\frac{2q}{3(q-2)}}\right) \ < \infty.
\end{aligned} 
\end{equation*}
Similarly, for $q>6,$ we have
\begin{equation*}
\frac{d}{dt}\|\u(t)\|_{L^q}\leq
C\|\u(t)\|_{L^q}^{\frac{2q-3}{q-3}}=C\|\u(t)\|_{L^q}\|\u(t)\|_{L^q}^{\frac{q}{q-3}}.
\end{equation*}
Again, by Gr\"onwall's lemma, we can obtain a bound on $\|\u(t)\|_{L^q}$
on any interval $[0,T]$ 
\begin{equation*}
\begin{aligned}
\|\u(t)\|_{L^q} & \leq \|\u_{0}\|_{L^q}\,\mbox{exp}\left(C\int_{0}^{T}\|\u(t)\|_{L^q}^{\frac{q}{q-3}}dt\right)\\
& \leq\|\u_{0}\|_{L^q}\mbox{exp}\left(C\K_{0}^{3}\right) \ < \infty.
\end{aligned} 
\end{equation*}
On the other hand, assuming there is a blow-up at $t=T^*$, we 
can suppose that $\|\u(t)\|_{L^q}= C / (T^*-t)^{\alpha}$, 
$\alpha>0$. Plugging this ansatz in expression \eqref{eq:LPSblowup}, we 
then have that $(\alpha\,p) \geq 1$ and
\begin{equation}
\frac{d}{dt}\|\u(t)\|_{L^q} =\frac{C}{(T^{*}-t)^{1+\alpha}}  
 = C\left[\frac{1}{(T^*-t)^{\alpha}}\right]^{\frac{1+\alpha}{\alpha}} 
 = C\|\u(t)\|_{L^q}^{\beta} \label{eq:duqdt},
\end{equation}
where $\beta=(1+\alpha)/\alpha \leq p+1=3(q-1)/(q-3)$, cf.~\eqref{eq:aprioriLq}.
We therefore conclude that, somewhat counter-intuitively,
for a blow-up to occur in Navier-Stokes flows,
 $\frac{d}{dt}\|\u(t)\|_{L^q}$ cannot 
be arbitrarily large. In fact, it must  be at most proportional to 
$\|\u(t)\|_{L^q}^{3(q-1)/(q-3)}$
as otherwise the integrand expression in \eqref{eq:LPSblowup} would 
blow-up before the integral becomes unbounded. In other words, if 
blow-up is to occur in a Navier-Stokes flow where the norm 
$\|\u(t)\|_{L^q}$ grows following relation \eqref{eq:duqdt}, 
then the exponent $\beta$ must lie within a $q$-dependent range with
the bottom bracket given by the exponents in \eqref{eq:noblowup_ub} 
and the upper bracket given by the exponent in \eqref{eq:aprioriLq}. These 
relations are illustrated in Figure \ref{fig:rateofgrowthLq}. Needless 
to say, for blow-up to occur, the growth rate in \eqref{eq:duqdt} must 
be sustained sufficiently long and this time grows without bound as the 
exponent $\beta$ drops towards the bottom bracket in 
\eqref{eq:noblowup_ub}. These diagnostics will be helpful in assessing 
how ``close" the extreme flows constructed in the present study come 
to becoming singular.

\section{Optimization Formulations}
\label{sec:OptFormulations}

In this section we formulate a number of optimization problems we 
consider in this study. They are motivated by the conditional 
regularity results discussed in Section \ref{sec:RegulResults} and aim 
at finding optimal initial conditions $\u_0$ that will produce 
``extreme'' flows, in the sense of maximizing the different regularity 
indicators involved in the conditional regularity results. The goal is 
to see whether a singularity might occur in such flows. With the domain 
$\Omega$ and the  viscosity $\nu$ fixed, each problem is defined by 
three parameters, namely, the function space where the initial data is 
sought, the size of the optimal initial data and the length of the time 
window. The perspective adopted here is that for a given size of the 
initial data we begin with a time window $[0,T]$ that is assumed 
sufficiently short such that classical solutions exist. The 
maximization is then performed on time windows with gradually 
increasing length $T$ to see whether the quantities controlling the 
regularity of the solution might reveal an unbounded growth that would 
signal approach to a singularity, cf.~\eqref{eq:Ensblowup} and 
\eqref{eq:LPSblowup}, as $T$ becomes longer than the interval of local 
existence. Thus, one has to work with smooth classical solutions only 
which facilitates numerical computations (cf.~Appendix 
\ref{sec:numer}).

To provide context for the present investigation, we first recall the 
optimization problem motivated by the enstrophy condition 
\eqref{eq:RegCrit_FoiasTemam} that was investigated in 
\cite{KangYunProtas2020}, where  initial conditions $\tuET$ with a 
fixed enstrophy $\E(\tuET) = \E_{0}$ were sought to maximize the 
enstrophy at prescribed times $T$
\setcounter{problem}{-1}
\begin{problem}\label{pb:maxET}
  Given $\E_0, T \in\mathbb{R}_+$, find
\begin{align*}
\tuET & =  \mathop{\arg\max}_{\u_0 \in {\Q}_{\E_0}} \, \E\left(\u(T);\u_0\right), \quad \text{where} \\
{\Q}_{\E_0} & :=  \left\{\u_0\in H^1(\Omega)\,\colon\,\bnabla\cdot\u_0 = 0, \; \E(\u_0) = \E_0 \right\}.
\end{align*} 
\end{problem}
\noindent A key element of the statement of this and subsequent 
optimization problems is the definition of the constraint manifold, 
here ${\Q}_{\E_0}$. While no evidence of unbounded growth of enstrophy 
\eqref{eq:Ensblowup} was detected \cite{KangYunProtas2020}, the main 
finding from that study was the observation that the maximum transient 
growth of enstrophy scales as 
\begin{equation}
\max_{\u_{0}} \max_{t>0} \E(\u(t;\u_{0})) = \OO\left(\E(\u_{0})^{3/2}\right) \qquad \text{as} \ \E(\u_{0}) \rightarrow \infty.
\label{eq:maxEt}
\end{equation}

\begin{figure}[t]
\mbox{
\subfigure[]{\includegraphics[width=0.49\textwidth]{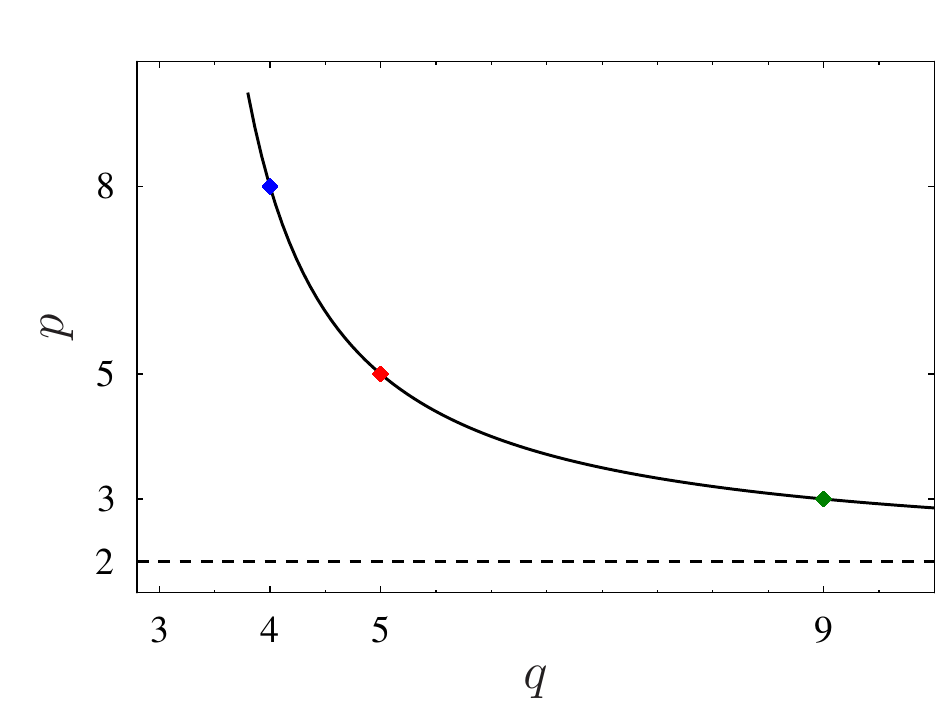}\label{fig:LPSq}}
\subfigure[]{\includegraphics[width=0.49\textwidth]{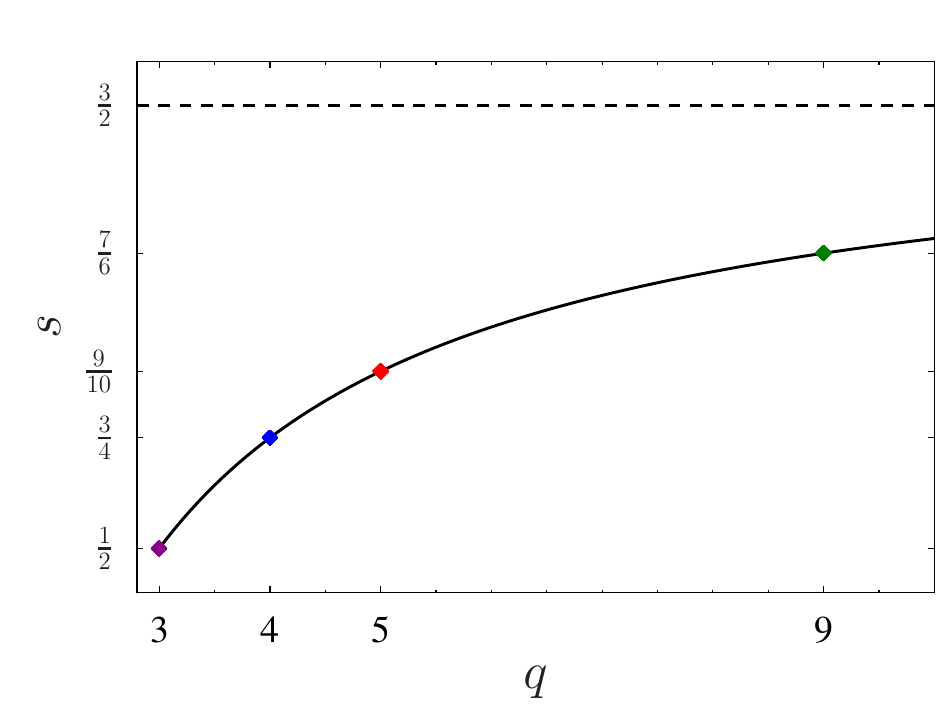}\label{fig:LPSs}}
}
\caption{
Dependence of \subref{fig:LPSq} $p$ in \eqref{eq:Phi} and 
\subref{fig:LPSs} $s$ in \eqref{eq:SobEmbb} on the index $q$ with solid 
symbols representing the values considered in the present study. The dashed 
horizontal lines correspond to the limiting values of $p$ and $s$ obtained 
when $q \rightarrow \infty$.
}
\label{fig:LPSq_s}
\end{figure}

Motivated by the Ladyzhenskaya-Prodi-Serrin conditions 
\eqref{eq:LPS}--\eqref{eq:LPS3}, we define the following two objective 
functionals
\begin{align}
\Phi_T^{q}(\u_0)  & := \frac{1}{T} \int_0^T \| \u(\tau) \|_{L^q}^{p} \, d\tau, 
\qquad p = \frac{2q}{q-3}, \quad q > 3, \label{eq:Phi} \\
\Psi_T(\u_0)  & :=  \| \u(T) \|_{L^3}^{3} \label{eq:Psi}
\end{align}
for a given $T > 0$.
The first one coincides with the integral in \eqref{eq:LPSblowup}, 
except for the prefactor $T^{-1}$ which ensures expression 
\eqref{eq:Phi}  remains bounded even if the integral grows without 
bound when $T$ tends to infinity. The values of the exponent $p$ are 
shown as a function of $q$ in Figure \ref{fig:LPSq}. The definition 
\eqref{eq:Psi} is motivated  by the observation that the norm $\| \cdot 
\|_{L^{\infty}(\Omega)}$ is not a differentiable function of its 
argument, cf.~\eqref{eq:normL3}, hence its use in the objective 
functional would result in a non-smooth optimization problem that would 
be much harder to solve \cite{r06}. This difficulty is circumvented by 
considering $\| \u(t) \|_{L^3}$ at a fixed time $t = T$, 
rather than over the interval $[0,T]$, in the objective functional 
\eqref{eq:Psi}.

An important aspect is the definition of the constraint manifold over 
which the objective functionals are to be maximized. The form of the 
functionals \eqref{eq:Phi}--\eqref{eq:Psi} involving Lebesgue norms $\| 
\u(t) \|_{L^q}$, $q \ge 3$, suggests that optimization should be 
performed over these spaces. However, the difficulty is that these 
spaces are not endowed with the Hilbert structure which, as will be 
discussed in the next section, is the preferred setting for the 
formulation of optimization algorithms. The reason is that an inner 
product is required for the usual definition of the gradient of the 
objective functional. Therefore, for each of the two objective 
functionals \eqref{eq:Phi}--\eqref{eq:Psi} we will consider two distinct 
formulations of the optimization problem:
\begin{itemize}
\item in the first, optimization is performed in the largest Sobolev 
space with the Hilbert structure $H^{s}(\Omega)$ embedded in the given 
Lebesgue space $L^q(\Omega)$, which will lend itself to solution 
using standard methods of numerical PDE optimization,
\item in the second, maximization is performed directly in the 
Lebesgue space $L^q(\Omega)$ leading to variational optimization  
problems with a nonstandard structure; in particular, in the absence of 
an inner product, this will necessitate the construction of a different 
notion of the gradient of the objective functional, namely, the metric 
gradient \cite{gt72}.
\end{itemize}
In regard to the first formulation, the relevant Hilbert-Sobolev space 
$H^{s}(\Omega)$  is determined based on the Sobolev embedding theorem 
in 3D \cite{af05}
\begin{equation}
 H^{s}(\Omega)\hookrightarrow L^{q}(\Omega), \quad \mbox{if}\quad s\ge\frac{3}{2}-\frac{3}{q}. 
 \label{eq:SobEmbb}
\end{equation}
In other words,  $H^{s}(\Omega)$ with $s = 3 / 2 - 3 / q$ is 
the largest Sobolev space embedded in $L^{q}(\Omega)$ and hence will 
serve as a proxy for the letter in the first formulation. The values of 
the index $s$ in \eqref{eq:SobEmbb} are shown as function of $q$ in 
Figure \ref{fig:LPSs}.

This thus leads to the following four optimization problems:
\begin{problem}\label{pb:PhiHs}
 Given $B, T \in\mathbb{R}_+$, $q>3$, $s=3/2-3/q$ and the objective functional $\Phi_T^q(\u_0)$ from
 \eqref{eq:Phi}, find
\begin{align}
\tuBT & =  \mathop{\arg\max}_{\u_0 \in {\M}_{B}} \, \Phi_T^q(\u_0), \quad \text{where} \\
 {\M}_{B} & :=  \left\{\u_0\in H^{s}(\Omega)\,\colon\,\bnabla_{\x}\cdot\u_0 = 0, \; \int_{\Omega} \u_0 \, d\x = 0,  \; \|\u_0\|_{L^q} = B \right\}.\label{manifold:pb1}
\end{align} 
\end{problem}
\begin{problem}\label{pb:PhiLq}
Given $B, T \in\mathbb{R}_+$, $q>3$ and the objective functional $\Phi_T^q(\u_0)$ from
\eqref{eq:Phi}, find
\begin{align}
\tuBT & =  \mathop{\arg\max}_{\u_0 \in \mathcal{L}_{B}} \, \Phi_T^q(\u_0), \quad \text{where} \\
 {\mathcal{L}}_{B} & :=  \left\{\u_0\in L^q\,\colon\,\bnabla_{\x}\cdot\u_0 = 0, \; \int_{\Omega} \u_0 \, d\x = 0,  \; \|\u_0\|_{L^q} = B \right\}. \label{manifold:pb2}
\end{align} 
\end{problem}
\begin{problem}\label{pb:PsiHs}
Given $B, T \in\mathbb{R}_+$ and the objective functional $\Psi_T(\u_0)$ from
\eqref{eq:Psi}, find
\begin{align}
\tuBT & =  \mathop{\arg\max}_{\u_0 \in {\N}_{B}} \, \Psi_T(\u_0), \quad \text{where} \\
 {\N}_{B} & :=  \left\{\u_0\in H^{1/2}(\Omega)\,\colon\,\bnabla_{\x}\cdot\u_0 = 0, \; \int_{\Omega} \u_0 \, d\x = 0,  \; \|\u_0\|_{L^3} = B \right\}.
\end{align} 
\end{problem}
\begin{problem}\label{pb:PsiLq}
Given $B, T \in\mathbb{R}_+$ and the objective functional $\Psi_T(\u_0)$ from
\eqref{eq:Psi}, find
\begin{align}
\tuBT & =  \mathop{\arg\max}_{\u_0 \in {\mathcal{S}}_{B}} \, \Psi_T(\u_0), \quad \text{where} \\
 {\mathcal{S}}_{B} & :=  \left\{\u_0\in L^3(\Omega)\,\colon\,\bnabla_{\x}\cdot\u_0 = 0, \; \int_{\Omega} \u_0 \, d\x = 0,  \; \|\u_0\|_{L^3} = B \right\}.
\end{align} 
\end{problem}
\noindent Problems \ref{pb:PhiHs} and \ref{pb:PsiHs} represent the 
first formulation mentioned above, whereas  Problems \ref{pb:PhiLq} and 
\ref{pb:PsiLq} represent the second.

In an earlier work \cite{KangProtas2021}, Problem \ref{pb:PhiHs} was 
studied for a single set only of parameters $q = 4, p = 8$, and  in 
the present investigation we consider Problems  \ref{pb:PhiHs} and 
\ref{pb:PhiLq}  for $q = 4,5,9$, in addition to studying Problems 
\ref{pb:PsiHs} and \ref{pb:PsiLq} (which correspond to $q = 3$, 
cf.~\eqref{eq:Psi}). The reason we focus on these particular values is 
because these are the only cases where both $q$ and $p$ are integers, 
which facilitates interpretation of the results. These values of $q$ 
are marked together with the corresponding values of $p$ and $s$ in 
Figures \ref{fig:LPSq} and \ref{fig:LPSs}, respectively. Computational 
approaches to find families of local maximizers of Problems 
\ref{pb:PhiHs}--\ref{pb:PsiLq} are presented next.

\section{Solution Approach}
\label{sec:SolApr}

First, we describe the Riemannian gradient method used to solve 
constrained optimization problems. A key element of this approach is 
evaluation of the gradient of the objective functional 
\eqref{eq:Phi}--\eqref{eq:Psi} with respect to the control variable 
which is the initial condition $\u_{0}$ in \eqref{eq:NS}. This needs to 
be done in a suitable topology and is described below. Finally, we 
introduce additional tools needed to implement the Riemannian gradient 
approach. 

We follow the ``optimize-then-discretize" strategy where the 
optimization approach is first formulated in the infinite-dimensional 
(continuous) setting and only then the resulting equations and 
expressions are discretized for the purpose of numerical solution 
\cite{g03}. Some details of the numerical implementation are presented 
in Appendix \ref{sec:numer}. While the general approach is essentially 
the same as used in \cite{KangProtas2021}, the main novelty is that 
here we also compute Lebesgue gradients, i.e, gradients of the 
objective functional \eqref{eq:Phi}--\eqref{eq:Psi} defined in Lebesgue 
spaces $L^{q}(\Omega)$, $q \ge 3$, without Hilbert structure, which are 
needed to solve Problems \ref{pb:PhiLq} and \ref{pb:PsiLq}. As will be 
evident below, this is technically more complicated than in Hilbert 
spaces \cite{gt72,protas2008}. To focus attention, we first describe 
solution of Problems \ref{pb:PhiHs} and \ref{pb:PhiLq}, and then 
mention adaptations needed to solve Problems \ref{pb:PsiHs} and 
\ref{pb:PsiLq} which have a similar structure.

\subsection{Riemannian Gradient Method}

Problem \ref{pb:PhiHs} is Riemannian {in the sense that} local 
maximizers $\tuBT$ belong to a constraint manifold $\mathcal{M}_B$ 
\cite{ams08}. To locally characterize this manifold, we define the 
tangent space $\T_{\z}\mathcal{M}_B$ at a point $\z \in \mathcal{M}_B$. 
To do so, the fixed-norm constraint can be expressed in terms of the 
function ${G}_{q} \; :\; H^s(\Omega) \rightarrow \RR_+$, where 
${G}_{q}(\z) := \|\z\|^q_{L^q}$. Computing the G\^{a}teaux 
differential of ${G}_{q}(\z) = B$ and using the Riesz representation 
theorem \cite{b77}, we obtain
\begin{equation}
\forall \z' \in H^{s}(\Omega) \qquad {G'}_{q}(\z;\z') := \frac{d}{d\epsilon} {G}_{q}(\z + \epsilon \z')\big|_{\epsilon = 0} = \left\langle \bnabla {G}_{q}(\z), \z' \right\rangle_{H^s} = 0,
\label{eq:dGHs}
\end{equation}
where $\bnabla {G}_{q}(\z) \in H^{s}(\Omega)$ is the gradient of the 
function ${G}_{q}$ at $\z \in \M_{B}$ and can be interpreted as an 
element orthogonal to the subspace $\T_{\z}\M_B$ in $H^{s}(\Omega)$. 
Thus, the tangent subspace is given by
\begin{equation}
\T_{\z}\M_B  := \left\{ \vvec \in H^s(\Omega) \, : \, \bnabla_{\x}\cdot\vvec = 0, \; \int_{\Omega} \vvec \, d\x = \0,  \; \left\langle \bnabla {G}_{q}(\z), \vvec \right\rangle_{H^s} = 0 \right\}.
\label{eq:TzHS} 
\end{equation}
A local maximizer $\tuBT$ will then be found by constructing a sequence
of divergence-free zero-mean vector fields with a fixed $L^q(\Omega)$ norm, 
$\left\{\uBT^{(n)}\right\}_{n\in\mathbb{N}}$, such that
\begin{equation}
\tuBT = \lim_{n\rightarrow\infty} \uBT^{(n)}.
\label{eq:tuBT}
\end{equation}
This sequence is defined using  the  following iterative procedure 
representing a discretization of a gradient 
flow projected onto the manifold $\M_B$
\begin{equation}
\begin{aligned}
\uBT^{(n+1)} & =  \R_{\M_B}\left(\;\uBT^{(n)} + \uptau_n \, \P_{\T_n\mathcal{M}_B}\bnabla^{H^s}\Phi_T^q\left(\uBT^{(n)}\right)\;\right), \\ 
\uBT^{(1)} & =  \u^0,
\end{aligned}
\label{eq:descHs}
\end{equation}
where $\u^0$ is an initial guess, $\P_{\T_n\mathcal{M}_B} \; : \; 
L^{q}(\Omega) \rightarrow \T_n \mathcal{M}_B:= 
\T_{\uBT^{(n)}}\mathcal{M}_B$ is an operator representing projection 
onto the tangent subspace \eqref{eq:TzHS} at the $n$th iteration, 
$\uptau_n$ is the step size, $\bnabla^{H^s}\Phi_T^q$ is the gradient of 
the functional $\Phi_T^q$ in the Sobolev space $H^s(\Omega)$, whereas 
$\R_{\mathcal{M}_B} \: : \; \T_n\mathcal{M}_B \rightarrow 
\mathcal{M}_B$ is a retraction from the tangent subspace to the 
constraint manifold \cite{ams08}.  Precise definitions of $\P_{\T_n 
\mathcal{M}_B},$ $\uptau_{n}$ and $\R_{\mathcal{M}_B}$ will be given in 
Section \ref{sec:projection} while computation of the gradient 
$\bnabla^{H^s}\Phi_T^q$ will be described in Section 
\ref{sec:eval_gradient}.

As regards solution of Problem \ref{pb:PhiLq}, the approach is quite 
similar to that in Problem \ref{pb:PhiHs}, with the difference that now 
we solve the optimization problem in the Banach space $L^q(\Omega)$ 
which lacks Hilbert structure. The maximizer $\tuBT$ belongs to the 
constraint manifold $\mathcal{L}_B$ and to locally characterize this 
manifold, we define the tangent subspace $\T_{\z}\mathcal{L}_B$ at a 
point $\z \in \mathcal{L}_B$. As above, the fixed-norm constraint is 
expressed in terms of the function $F_{q} \; : \; L^{q}(\Omega) 
\rightarrow \RR_+$, $F_{q}(\z) := \|\z\|^q_{L^q}$ which is 
now defined on the Lebesgue space $L^q(\Omega)$. Computing the 
G\^{a}teaux differential of ${F}_{q}(\z) = B$, we obtain
\begin{equation}
\forall \z' \in L^{q}(\Omega) \qquad {F'}_{q}(\z;\z') := \frac{d}{d\epsilon} {F}_{q}(\z + \epsilon \z')\big|_{\epsilon = 0} = \left\langle \bnabla {F}_{q}(\z), \z' \right\rangle_{(L^{q})^{*}\times L^{q}} = 0,
\label{eq:dFLq}
\end{equation}
where $\langle \cdot, \cdot \rangle_{(L^{q})^{*}\times 
L^{q}}$ is the duality pairing between the Lebesgue space 
$L^q(\Omega)$  and its dual $(L^{q}(\Omega))^{*}$ which can be 
identified with the space $L^p(\Omega)$  where $1/p+1/q=1$ \cite{af05}. 
Thus, the subspace is defined as 
\begin{equation}
\T_{\z}\mathcal{L}_B  := \left\{ \vvec \in L^{q}(\Omega) \, : \, \bnabla_{\x}\cdot\vvec = 0, \int_{\Omega} \vvec \, d\x = \0,  \left\langle \bnabla F_{q}(\z), \vvec \right\rangle_{(L^{q})^{*}\times L^{q}} = 0 \right\},
\label{eq:TzLB} 
\end{equation}
which is a natural extension of \eqref{eq:TzHS} to the case of a general Banach space
with the inner product replaced  by a duality pairing to define the tangency condition.

Thus, local maximizers $\tuBT$  in Problem \ref{pb:PhiLq} are found as 
in \eqref{eq:tuBT}--\eqref{eq:descHs}, except that the Sobolev gradient 
$\bnabla^{H^s}\Phi_T^q$ of the objective functional needs to be 
replaced with the Lebesgue gradient $\bnabla^{L^q}\Phi_T^q$, whereas 
the projection and retraction operators have to be defined as 
$\P_{\T_n\mathcal{L}_B} \; : \; L^{q}(\Omega) \rightarrow \T_n 
\mathcal{L}_B:= \T_{\uBT^{(n)}}\mathcal{L}_B$  and 
$\R_{\mathcal{L}_B}\;:\; \T_n\mathcal{L}_B \rightarrow \mathcal{L}_B$, 
such that now the discretized gradient flow on the manifold 
$\mathcal{L}_B$ takes the form
\begin{equation}
\begin{aligned}
\uBT^{(n+1)} & =  \R_{\mathcal{L}_B}\left(\;\uBT^{(n)} + \uptau_n \, \P_{\T_n\mathcal{L}_B}\bnabla^{L^q}\Phi_T^q\left(\uBT^{(n)}\right)\;\right), \\ 
\uBT^{(1)} & =  \u^0.
\end{aligned}
\label{eq:descLq}
\end{equation}

\subsection{Evaluation of the Gradient}
\label{sec:eval_gradient}

A key element of the iterative procedure \eqref{eq:descHs} is 
evaluation of the gradients $\bnabla^{H^s}\Phi_T^q$ and 
$\bnabla^{L^q}\Phi_T^q$ of the objective functional $\Phi_T^q$, 
cf.~\eqref{eq:Phi}. The first step in determining these gradients is to 
find the gradient of \eqref{eq:Phi} with respect to the $L^2$ topology. 
We begin by considering the G\^{a}teaux differential of the objective 
functional \eqref{eq:Phi} 
\begin{equation}
(\Phi_T^q)'(\u_0;\u_0') := \lim_{\epsilon \rightarrow 0}
\epsilon^{-1} \left[ \Phi_T^q(\u_0+\epsilon \u_0') - \Phi_T^q \right]
\label{eq:dPhi}
\end{equation}
which with $\u_{0}$ fixed is a bounded linear functional of its second 
argument in $L^{2}(\Omega)$. Then, invoking the Riesz representation 
theorem \cite{b77}, we have 
\begin{equation}
(\Phi_T^q)'(\u_0;\u_0')
= \Big\langle \bnabla^{L^2}\Phi_T^q(\u_0), \u_0' \Big\rangle_{L^2}.
\label{eq:rieszL2}
\end{equation}
As shown in Appendix \ref{sec:L2grad}, an expression for the G\^{a}teaux differential can be conveniently obtained by solving the following adjoint system
\begin{subequations}
\label{eq:aNSE3D}
\begin{align}
 \mathcal{L}^*\begin{bmatrix} \u^* \\ p^* \end{bmatrix} := 
& \begin{bmatrix}
-\partial_{t}\u^*-\left[\bnabla_{\x}\,\u^*+\left(\bnabla_{\x}\,\u^{*}\right)^T\right]\u-\bnabla_{\x}\, p^*-\nu\laplacian\u^* \\
-\bnabla_{\x}\cdot\u^*
\end{bmatrix}  = \begin{bmatrix} \f  \\ 0\end{bmatrix}, \label{eq:aNSE3Da} \\
\quad \f(t,\x)  := & \frac{2q}{(q-3)T}  \|\u(t)\|_{L^q}^{\frac{q(5-q)}{q-3}} \,  |\u(t,\x)|^{q-2} \u(t,\x), \qquad \x\in \Omega, \, t \in [0,T], \label{eq:aNSE3Df} \\
 \u^*(T)= & \0  \label{eq:aNSE3Db}
\end{align}
\end{subequations}
subject to periodic boundary conditions, where $\u^* \, : \,{}
[0,T]\times\Omega \rightarrow \RR^3$ and $p^* \, : \, [0,T]\times\Omega 
\rightarrow \RR$ are the {\em adjoint states} associated to the 
velocity field $\u$ and the pressure $p$, respectively, in 
\eqref{eq:NS}.  The $L^2$ gradient is then obtained from 
\eqref{eq:rieszL2} as
\begin{equation}
\bnabla^{L^2}\Phi_T^q = \u^*(0).
\label{eq:gradL2}
\end{equation}
We note that the adjoint system \eqref{eq:aNSE3D} is a linear 
terminal-value problem and hence needs to be integrated backwards in 
time. The coefficients in \eqref{eq:aNSE3Da} and the source term 
\eqref{eq:aNSE3Df} are determined by the state $\u = \u(t,\x)$ corresponding to 
the initial condition $\u_{0}$ around which linearization is performed. 
When Problems \ref{pb:PsiHs} and \ref{pb:PsiLq} are solved, the $L^{2}$ 
gradient of the objective functional \eqref{eq:Psi} is given by an 
expression analogous to \eqref{eq:gradL2}, namely
\begin{equation}
\bnabla^{L^2}\Psi_T^q = \u^*(0)
\label{eq:gradL3L2}
\end{equation}
where $\u^{*}$ is the solution of the adjoint problem \eqref{eq:aNSE3Da} with
\begin{subequations}
\label{eq:aL3NSE3D}
\begin{align}
 \f(t,\x)  = & \, \0, \qquad t \in [0,T],\ \x\in \Omega, \label{eq:aL3NSE3Df} \\
 \u^*(T,\x)= & \, 3 \, |\u(T,\x)|\u(T,\x).  \label{eq:aL3NSE3Db}
\end{align}
\end{subequations}
The difference between the structure of the source term and the 
terminal condition in \eqref{eq:aNSE3Df}--\eqref{eq:aNSE3Db} and in 
\eqref{eq:aL3NSE3Df}--\eqref{eq:aL3NSE3Db} reflects the difference in 
the form of the objective functionals \eqref{eq:Phi} and \eqref{eq:Psi} 
defined, respectively, as an integral over $[0,T]$ and an expression 
evaluated at $t = T$. Once the $L^{2}$ gradient is found as in 
\eqref{eq:gradL2}--\eqref{eq:gradL3L2}, it is possible to construct the 
Sobolev and Lebesgue gradients $\bnabla^{H^s}\Phi_T^q$ and 
$\bnabla^{L^q}\Phi_T^q$ as described in the following subsections.

\subsubsection{Gradient in the Sobolev Space $H^s$}

The Sobolev gradient is obtained by reinterpreting the G\^{a}teaux 
differential \eqref{eq:rieszL2}, viewed as a function of its second 
argument  with $\u_{0}$ fixed, as a bounded linear functional on the 
Sobolev space $H^{s}(\Omega)$ \cite{Neuberger2010book,pbh04}. 
Therefore, the differential also admits the Riesz representation
\begin{equation}
(\Phi_T^q)'(\u_0;\u_0')
= \Big\langle \bnabla^{H^s}\Phi_T^q(\u_0), \u_0' \Big\rangle_{H^s}.
\label{eq:rieszHs}
\end{equation}
Motivated by computational consideration mentioned below, here we use 
an equivalent form of the inner product in $H^{s}(\Omega)$, namely, 
$\big\langle \f,\g\big\rangle_{H^s}:=\int_\Omega \f\cdot\g \,d\x 
+\ell^{2s}\int_\Omega \Delta^{s/2}\f\cdot\Delta^{s/2}\g\, d\x$, where 
$0 < \ell < \infty$ is an adjustable parameter. Identifying the two 
Riesz representations \eqref{eq:rieszL2} and \eqref{eq:rieszHs} of the 
G\^{a}teaux differential, using \eqref{eq:gradL2} together with this 
definition of the inner product, performing integration by parts and 
noting the arbitrariness of the perturbation  $\u_{0}'$, we obtain 
the elliptic boundary-value problem \cite{pbh04}
\begin{equation}
\left[\mbox{Id}-\ell^{2s}\Delta^{s}\right]\bnabla^{H^s}\Phi_T^q=\bnabla^{L^2}\Phi_T^q \quad \mbox{in}\;\Omega,
\label{eq:HsBVP}
\end{equation}
subject to periodic boundary conditions which allows us to determine 
the Sobolev gradient $\bnabla^{H^s}\Phi_T^q$ when the $L^{2}$ gradient 
is available from \eqref{eq:gradL2}--\eqref{eq:gradL3L2}. Clearly, 
\eqref{eq:HsBVP} preserves the divergence-free property of the $L^{2}$ 
gradient. Transforming equation \eqref{eq:HsBVP} to the Fourier space, 
where $\left[\widehat{\f}\right]_{\k}\in\CC^3$ represents the Fourier 
coefficient of the vector field $\f$ with wavenumber $\k$, we obtain 
\begin{subequations}
\label{eq:SobvgradientFourier}
\begin{align}
\left[\widehat{\bnabla^{H^s}\Phi_T^q}\right]_{\k}&= \overbrace{\frac{1}{1+\ell^{2s}|\k|^{s}}}^{\mathcal{F}(k)} \left[\bnabla^{L^2}\Phi_T^q\right]_{\k}, \qquad \k\in\mathbb{Z}^3\setminus \{\0\},\label{eq:SobvgradientFouriera}\\
\left[\widehat{\bnabla^{H^s}\Phi_T^q}\right]_{\0}&=\0,\label{eq:SobvgradientFourierb}
\end{align}
\end{subequations}
where $k = |\k|$. This demonstrates that the computation of the Sobolev 
gradient $\bnabla^{H^s}\Phi_T^q$ can be regarded as application of the 
low-pass filter $\mathcal{F}(k)$, which is a smoothing operation, to 
the $L^{2}$ gradient \eqref{eq:gradL2}--\eqref{eq:gradL3L2} with the 
cut-off wavenumber given by $\ell^{-1}$ \cite{pbh04}. While problems 
with different values of $\ell \in (0, \infty)$ are mathematically 
equivalent, adjusting the value of this parameter can in practice have 
a significant effect on the rate of convergence of iterations 
\eqref{eq:descHs}.

\subsubsection{Gradient in the Lebesgue Space $L^q$}
\label{sec:gradLq}

Determination of the gradient in the Lebesgue space $L^{q}(\Omega)$, $q 
\ge 3$, needed in Problems \ref{pb:PhiLq} and \ref{pb:PsiLq} is more 
involved as we do not have an inner product in this space and hence 
there is no Riesz identity. We thus need  to invoke the concept of a 
{\em metric} gradient \cite{gt72}, which is a generalization of the 
notion of the gradient to normed spaces. It relies on the observation 
that the gradient is the element maximizing the directional derivative 
of the objective functional under certain constraints.  We can thus 
define the gradient of the objective functional \eqref{eq:Phi} 
variationally as a solution of the following optimization subproblem 
subject to the suitable constraints \cite{Neuberger2010book,protas2008}
\begin{equation}
\bnabla^{L_{q}} \Phi^q_{T} = \mathop{\arg\max}_{\|\boldsymbol{\Theta}\|_{L^q}=1,\,\bnabla_{\x}\cdot\boldsymbol{\Theta}=0} 
\left[ (\Phi_T^q)'(\u_0;\u_0') = 
 \left\langle \bnabla^{L^{2}} \Phi^q_{T}, \boldsymbol{\Theta}\right\rangle_{(L^{q})^{*}\times L^{q}} \right].
\label{eq:OptGradLq}
\end{equation}
Optimization problem \eqref{eq:OptGradLq} can be converted to an 
unconstrained form by introducing the Lagrange multipliers $\mu\in\RR$ 
and $\eta\in L^q(\Omega)$ associated with each of these two constraints
\begin{multline}
\bnabla^{L^{q}} \Phi^q_{T} = 
\mathop{\arg\min}_{\mu\in\RR,\,\eta\in L^{q}(\Omega)}\mathop{\arg\max}_{\boldsymbol{\Theta}\in L^{q}(\Omega)}  \left[\left\langle \bnabla^{L^2} \Phi^q_{T},\boldsymbol{\Theta}\right\rangle_{(L^{q})^{*}\times L^{q}}+\right. \\ 
\left.+\frac{\mu}{q}\left(\|\boldsymbol{\Theta}\|^{q}_{L^q}-1\right)+\int_{\Omega} \eta\left(\bnabla_{\x}\cdot\boldsymbol{\Theta}\right)\,d\x\right],
\end{multline}
where the expression in the bracket can be interpreted as the Lagrangian
corresponding to the right-hand side of \eqref{eq:OptGradLq}. After 
performing integration by parts, we arrive at
\begin{equation}
\label{eq:grad_lagrange}
\bnabla^{L^{q}} \Phi^q_{T} = 
{\mathop{\arg\min}_{\mu\in\RR,\,\eta\in L^{q}(\Omega)}}\mathop{\arg\max}_{\boldsymbol{\Theta}\in L^{q}(\Omega)} \int_{\Omega} \left( \u^{*}{(0)}\cdot\boldsymbol{\Theta}  
+\frac{\mu}{q}|\boldsymbol{\Theta}|^{q}-\frac{\mu}{q|\Omega|}-\boldsymbol{\Theta}\cdot \bnabla_{\x}\,\eta\right)\,d\x.
\end{equation}
Now, using the optimality conditions requiring the vanishing of the 
G\^ateaux differentials of the integral expression in 
\eqref{eq:grad_lagrange} with respect to $\boldsymbol{\Theta}$, $\mu$ 
and $\eta$, we obtain the following expressions
\begin{subequations}
\begin{alignat}{2}
\int_{\Omega} \left( \u^{*}{(0)}  
+\mu \bnabla^{L^{q}}\Phi^q_{T}|\bnabla^{L^{q}}\Phi^q_{T}|^{q-2}- \bnabla_{\x}\,\eta\right)\cdot \boldsymbol{\Theta}'\,d\x&=0,\;\qquad\forall\; \boldsymbol{\Theta}'\in L^{q}(\Omega){,} \label{eq:OptCond1}\\
\int_{\Omega} |\boldsymbol{\Theta}|^q\,d\x&=1,\qquad\forall\; \boldsymbol{\Theta}\in L^{q}(\Omega){,} \label{eq:OptCond2}\\
\int_{\Omega} \left(\bnabla_{\x}\cdot\boldsymbol{\Theta}\right)\,\eta\,d\x&=0,\qquad\forall\; \eta\in L^{q}(\Omega).\label{eq:OptCond3}
\end{alignat}
\end{subequations}
Conditions \eqref{eq:OptCond2} and \eqref{eq:OptCond3} trivially say 
that solutions of the optimization problem \eqref{eq:OptGradLq} have a 
unit $L^{q}(\Omega)$ norm and are divergence-free. Since 
$\boldsymbol{\Theta}'$ is arbitrary, condition \eqref{eq:OptCond1} is 
equivalent to the following relation
\begin{equation}
\label{eq:gradLq1}
\bnabla^{L^{q}}\Phi^q_{T}|\bnabla^{L^{q}}\Phi^q_{T}|^{q-2}=\frac{1}{\mu}\left(-\u^{*}{(0)}+\bnabla_{\x}\,\eta\right),  \quad \x\in\Omega.
\end{equation}
As observed in \cite{protas2008}, the RHS in expression 
\eqref{eq:gradLq1} represents the Helmholtz-Weyl decomposition of a 
vector field in $L^q(\Omega)$ \cite{RobinsonRodrigoSadowski2016}, where 
such a vector field is written as the sum of a divergence-free 
vector field and the gradient of a certain scalar potential. Then, 
applying the divergence operator to  \eqref{eq:gradLq1} and noticing 
that both $\u^*(0)$ and $\bnabla^{L^q} \Phi_{T}^q$ are divergence-free, 
we obtain an elliptic boundary-value problem which can 
be used to determine $\eta$
\begin{equation}
\label{eq:gradLq2}
\Delta\eta=\mu\left(\bnabla_{\x}\,|\bnabla^{L^{q}}\Phi_{T}^q|^{q-2}\right)\cdot\bnabla^{L^{q}}\Phi_T^{q}, \quad \x\in\Omega
\end{equation}
Therefore, to determine  the metric gradient in the space 
$L^q(\Omega),$ it is necessary to solve the following non-linear system 
subject to periodic boundary conditions
\begin{subequations}\label{eq:LGsystem}
\begin{align}
\bnabla^{L^{q}}\Phi^q_{T}|\bnabla^{L^{q}}\Phi^q_{T}|^{q-2}=\frac{1}{\mu}\left(-\u^{*}{(0)}+\bnabla_{\x}\,\eta\right),&\quad \x\in\Omega,\label{eq:LGsystema}\\
\Delta\eta=\mu\left(\bnabla_{\x}\,|\bnabla^{L^{q}}\Phi_T^q|^{q-2}\right)\cdot\bnabla^{L^{q}}\Phi_T^q, &\quad\x\in\Omega{,} \label{eq:LGsystemb}\\
\|\bnabla^{L^{q}}\Phi_T^q\|_{L^q}=1{.}\label{eq:LGsystemc}&
\end{align}
\end{subequations}
We emphasize that, in contrast to the Hilbert-space setting where the 
map between the $L^{2}$ gradient and  the corresponding Sobolev 
gradient $\bnabla^{H^s} \Phi_{T}^q$ is {\em linear}, 
cf.~\eqref{eq:HsBVP}, in the present case of Lebesgue spaces this 
relation is nonlinear. We note that in the special case with $q = 2$, 
system \eqref{eq:LGsystem} becomes linear and the normalization 
condition \eqref{eq:LGsystemc} can be dropped, so that we have
\begin{subequations}
\begin{align}
\bnabla^{L^{2}}\Phi^q_{T}=\frac{1}{\mu}\left(-\u^{*}{(0)}+\bnabla_{\x}\,\eta\right),  &\qquad\x\in\Omega \label{eq:gradL2_2a}\\
\Delta\eta=0, & \qquad\x\in\Omega.
\end{align}
\end{subequations}
Since in this case $\eta$ is harmonic and equal to a constant, 
equation \eqref{eq:gradL2_2a} becomes 
$\bnabla^{L^{2}}\Phi^q_{T}=- (1 / \mu) \u^{*}{(0)}$ which, up to 
normalization, is the same as \eqref{eq:gradL2}. Finally, we add that a 
rigorous demonstration of the existence and possible uniqueness of 
solutions of system \eqref{eq:LGsystem} remains an open problem.

\subsection{Projection, Retraction and Arc-Maximization}
\label{sec:projection}

In Problems \ref{pb:PhiHs} and \ref{pb:PsiHs} the condition 
characterizing the subspaces tangent to the manifolds $\M_{B}$ and 
$\N_{B}$ has the form $\left\langle\bnabla 
G_{q}(\z),\z'\right\rangle_{H^s} =\left\langle 
q\,|\z|^{q-2}\z,\z'\right\rangle_{\dot{H}^s}$ for all $\z'\in 
H^s(\Omega)$. cf.~\eqref{eq:dGHs}. We note that given the nonlinearity 
of the term $|\z|^{q-2}\z$, the element $\bnabla G_{q}(\z)$ does not, 
in general, satisfy the divergence-free and zero-mean conditions, even 
if they are satisfied by $\z$. Therefore, the projection operator 
$\mathcal{P}_{\T_n\mathcal{M}_{B}}:\,H^{s}(\Omega)\rightarrow 
\T_{n}\mathcal{M}_{B}$ in iterations \eqref{eq:descHs} is defined as 
\cite{ams08,KangProtas2021}
\begin{align}
\mathcal{P}_{\T_n\mathcal{M}_{B}} \z & := 
\overline{\z} - \frac{\left\langle \z,\bnabla {G}_{q}\left(\uBT^{(n)}\right)\right\rangle_{H^s}}{
\left\langle\overline{\bnabla F_{q}\left(\uBT^{(n)}\right)},\bnabla {G}_{q}\left(\uBT^{(n)}\right)\right\rangle_{H^s}}\, \overline{\bnabla {G}_{q}\left( \uBT^{(n)}\right)}, \label{eq:PHs} \\
& \text{where} \quad  \overline{\v} := \v - \bnabla \Delta^{-1} (\bnabla \cdot \v) - \int_{\Omega} \v \, d\x \nonumber
\end{align}
and likewise for the projection operator 
$\mathcal{P}_{\T_n\mathcal{N}_{B}}:\,H^{1/2}(\Omega)\rightarrow 
\T_{n}\mathcal{N}_{B}$, which ensures the divergence-free and zero-mean 
properties are satisfied by construction. In Problems \ref{pb:PhiLq} 
and \ref{pb:PsiLq} the tangency condition has the form \eqref{eq:dFLq} 
such that the projection operator 
$\mathcal{P}_{\T_n\mathcal{L}_{B}}:\,L^{q}(\Omega)\rightarrow 
\T_{n}\mathcal{L}_{B}$ is 
\begin{equation}
\mathcal{P}_{\T_n{\mathcal{L}_{B}}} \z := 
\overline{\z} - \frac{\left\langle \bnabla F_{q}\left(\uBT^{(n)}\right),\overline{\z} \right\rangle_{(L^{q})^{*}\times L^{q}}}{\int_{\Omega}\bnabla F_{q}\left(\uBT^{(n)}\right)\cdot\overline{\bnabla F_{q}\left( \uBT^{(n)}\right)} \,d\x}\, \overline{\bnabla F_{q}\left( \uBT^{(n)}\right)}
\label{eq:PLq} 
\end{equation}
and likewise for the projection operator 
$\mathcal{P}_{\T_n\mathcal{S}_{B}}:\,L^{3}(\Omega)\rightarrow 
\T_{n}\mathcal{S}_{B}$. We observe that in this case, inner products 
have been replaced by duality parings due to the lack of Hilbert 
structure in Lebesgue spaces. 

\begin{figure}[t]
\centering
\includegraphics[width=0.5\textwidth]{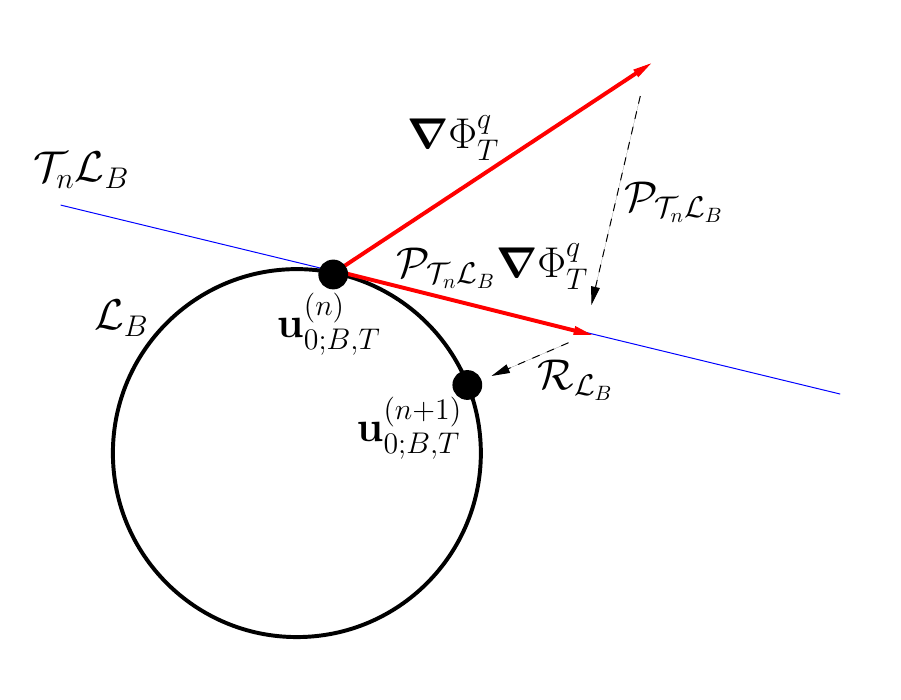}
\caption{Schematic representation of the operations performed at each 
iteration of Algorithm \ref{alg:optimAlg} with projection 
$\mathcal{P}_{\T_n\mathcal{L}_{B}}:\,L^{q}(\Omega)\rightarrow 
\T_{n}\mathcal{L}_{B}$  and retraction 
$\mathcal{R}_{\mathcal{L}_B}:\,\T_{n}\mathcal{L}_{B}\rightarrow\mathcal{L}_{B}$ 
\cite{KangProtas2021}.}
\label{fig:ProjLq}
\end{figure}

The retraction operator 
$\mathcal{R}_{\mathcal{M}_B}:\,\T_{n}\mathcal{M}_{B}\rightarrow\mathcal{M}_{B}$ 
in \eqref{eq:descHs} is defined as the normalization 
\cite{ams08,KangProtas2021}
\begin{equation}
\mathcal{R}_{\mathcal{M}_B}(\z) := \frac{B}{\|\z\|_{L^q}}\,\z \qquad \mbox{for all}\; \z \in \T_n\mathcal{M}_{B}
\label{eq:R_MB}
\end{equation}
and likewise for the retraction operators in the solution of Problems 
\ref{pb:PhiLq}, \ref{pb:PsiHs} and \ref{pb:PsiLq}. The step size 
$\uptau_{n}$ in the iterative algorithms  \eqref{eq:descHs} and \eqref{eq:descLq} 
is found by 
solving the arc-search problems
\begin{equation}\label{eq:tau_nHS}
\uptau_n = \mathop{\arg\max}_{\uptau>0} \Phi^q_T\left[ \mathcal{R}_{\M_B}\left(\;\uST^{(n)} + \uptau \, \mathcal{P}_{\T_n{\mathcal{M}_{B}}}\bnabla^{H^s}\Phi^q_T\left(\uBT^{(n)}\right)\;\right) \right]
\end{equation}
and 
\begin{equation}\label{eq:tau_nLq}
\uptau_n = \mathop{\arg\max}_{\uptau>0} \Phi^q_T\left[ \mathcal{R}_{\mathcal{L}_B}\left(\;\uST^{(n)} + \uptau \, \mathcal{P}_{\T_n{\mathcal{L}_{B}}}\bnabla^{L^q}\Phi^q_T\left(\uBT^{(n)}\right)\;\right) \right],
\end{equation}
respectively. Analogous formulations are employed when solving Problems 
\ref{pb:PsiHs} and \ref{pb:PsiLq}. These problems can be regarded as a 
generalization of the standard line-search problem with maximization 
performed following an arc (a geodesic in the limit of infinitesimal 
step sizes) lying on the constraint manifold $\M_B$ and 
$\mathcal{L}_B$, rather than along a straight line. Computations 
involved in the discrete gradient flow applied to solve Problem 
\ref{pb:PhiLq} are summarized as Algorithm \ref{alg:optimAlg} and 
illustrated schematically in Figure \ref{fig:ProjLq}. We have also 
implemented the Riemannian conjugate-gradient method \cite{ams08}, but 
since it did not lead to appreciable improvements in terms of 
convergence, we do not use it here (we add that development of this 
approach for problems defined on Lebesgue spaces not endowed with an 
inner product leads to some interesting theoretical questions 
\cite{s2005,r2025}).

\bigskip 

\begin{algorithm}[H]
\begin{algorithmic}
\State $\uBT^{(0)} = \u^{0}$
\State Compute $\Phi^q_{T}(\u_{0})$
\State $n = 0$

\Repeat 

\State  \{-------------------------- Optimization Iterations \eqref{eq:descHs}  -------------------------------- \}

\State Solve the Navier-Stokes system with initial condition $\uBT^{(n)}$, see equation \eqref{eq:NS}

\State Solve the adjoint system to obtain $\u^*$ and $p^*$, see equation \eqref{eq:aNSE3D}

\State Compute the $L^2$ gradient $\bnabla^{L_2}\Phi^q_T\left(\uBT^{(n)}\right)$, see equation \eqref{eq:gradL2}

\State Compute the Lebesgue gradient $\bnabla^{L^q}\Phi^q_T\left(\uBT^{(n)}\right)$,
            see system \eqref{eq:LGsystem}

\State Compute the optimal step size $\uptau_n$, see equation \eqref{eq:tau_nLq}

\State Set $\uBT^{(n+1)} = \mathcal{R}_{\mathcal{L}_{B}}\left(\;\uBT^{(n)} + \uptau_n \mathcal{P}_{\mathcal{T}_{n}{\mathcal{L}_{B}}}\bnabla^{L^q}\Phi_T^q\left(\uBT^{(n)}\right)\;\right)$

\State Evaluate the termination condition \texttt{relative\_change $ =\frac{\Phi_T^q\left(\uBT^{(n+1)}\right)-\Phi_T^q\left(\uBT^{(n)}\right)}{\Phi_T^q\left(\uBT^{(n)}\right)}$}

\State Set $n=n+1$

\Until{ \ \texttt{relative\_change} $<$ $\epsilon$ or $n>N_{\max}$}
\State $\tuBT = \uBT^{(n+1)}$

\end{algorithmic}
\caption{
Solution of Problem \ref{pb:PhiLq} for fixed $T$ and $B$.  \newline
     \textbf{Input:} \newline
 \hspace*{0.22cm} $B$ --- size of the initial data \newline
 \hspace*{0.22cm} $T$ --- length of the time window \newline
 \hspace*{0.22cm} $\u^{0}$ --- initial guess \newline
 \hspace*{0.22cm}    $\epsilon$ --- tolerance in the termination condition for iterations \eqref{eq:descLq} \newline
  \hspace*{0.22cm}  $N_{\max}$ --- maximum number of iterations \newline
 \textbf{Output:} \newline
 \hspace*{0.22cm}  Optimal initial data {$\tuBT\in\mathcal{L}_B.$}
}
\label{alg:optimAlg}
\end{algorithm}

\section{Numerical Results}
\label{sec:results}

In this section, we solve Problems \ref{pb:PhiHs}, \ref{pb:PhiLq}, 
\ref{pb:PsiHs} and \ref{pb:PsiLq} for different values of $q,$ $B,$ $T$ 
and using different initial guesses $\u^{0},$ cf.~Algorithm 
\ref{alg:optimAlg}. These problems were tackled using the two variants 
of the gradient descent method introduced in Section \ref{sec:SolApr}, 
namely, \eqref{eq:descHs} and \eqref{eq:descLq}. Our analysis focuses 
on all integer values of $(q,p) \in \{ (4,8), (5,5), (9,3) \}$ in 
addition to $q = 3$, cf.~Figure \ref{fig:LPSq_s}, whereas the parameter 
$s$ in Problems \ref{pb:PhiHs} and \ref{pb:PsiHs} depends on $q$ via 
the Sobolev embedding \eqref{eq:SobEmbb}. This represents a large 
number of parameters, and while we have explored their many different 
combinations, detailed analysis will be provided for the extreme values 
$q = 3$ and $q = 9$ in subsections \ref{sec:L3} and \ref{sec:L9} below 
with a summary information for all considered values of $q$ collected 
in Section \ref{sec:Lqdiagn}. We note that the case $q = 4$  was 
extensively investigated, albeit in the context of Problem 
\ref{pb:PhiHs} only, in \cite{KangProtas2021}.

Without loss of generality, we set $\nu = 1$ in \eqref{eq:NS} and, 
unless stated otherwise, the numerical resolution used in computations 
in Algorithm \ref{alg:optimAlg} is $N^{3} = 256^{3}$. All computations 
reported below are well resolved which was verified by carefully 
examining the energy spectra at different flow stages (details are 
provided in \cite{r2025}). We employed a variety of initial guesses 
$\u^{0}$, including the (suitably rescaled) optimal initial conditions 
found for $q = 4$ in \cite{KangProtas2021} and various states with a 
random structure. The constraint manifolds in Problems 
\ref{pb:PhiHs}--\ref{pb:PsiLq} are defined in terms of the parameter 
$B$ fixing the ``size" of the initial data $\u_{0}$ measured in the 
$L^q$ norm. When solving each of Problems 
\ref{pb:PhiHs}--\ref{pb:PsiLq} we consider three distinct values of the 
parameter $B$ and in order for these values to be consistent across 
different $q = 3,4,5,9$, they are chosen such that the values of the 
Taylor-scale Reynolds number \cite{frisch1995turbulence} characterizing 
the corresponding initial conditions 
\begin{equation} 
Re:=\sqrt{\frac{10}{3\,|\Omega|}}\frac{\K(\u_{0})}{\nu\sqrt{\E(\u_{0})}} 
\label{eq:Re} 
\end{equation} 
are approximately the same for the smallest, intermediate and largest 
values of $B$. This ensures that solutions of Problems 
\ref{pb:PhiHs}--\ref{pb:PsiLq}  obtained with different values of $q$ 
and $B$ selected in this way are characterized by a similar balance of 
linear viscous and nonlinear advection effects at $t = 0$. We add that 
since the values of $\K\left( \uBT^{(n)} \right)$ and $\E\left( 
\uBT^{(n)} \right)$ change during iterations \eqref{eq:descHs} and 
\eqref{eq:descLq}, the Reynolds numbers \eqref{eq:Re} characterizing 
the optimal initial conditions $\tuBT$ found for different $q$ and $B$ 
are only approximately equal. Finally, we state that in our solutions 
of Problems \ref{pb:PhiHs}--\ref{pb:PsiLq} we have found no evidence 
for an unbounded growth of quantities \eqref{eq:Phi}--\eqref{eq:Psi} 
that would signal an incipient singularity. However, as documented 
below, these quantities were found to exhibit a significant transient 
growth.


\subsection{Extremal Flows in $L^{3}(\Omega)$}
\label{sec:L3}

\begin{figure}[H]
\mbox{
\subfigure[]{\includegraphics[width=0.49\textwidth]{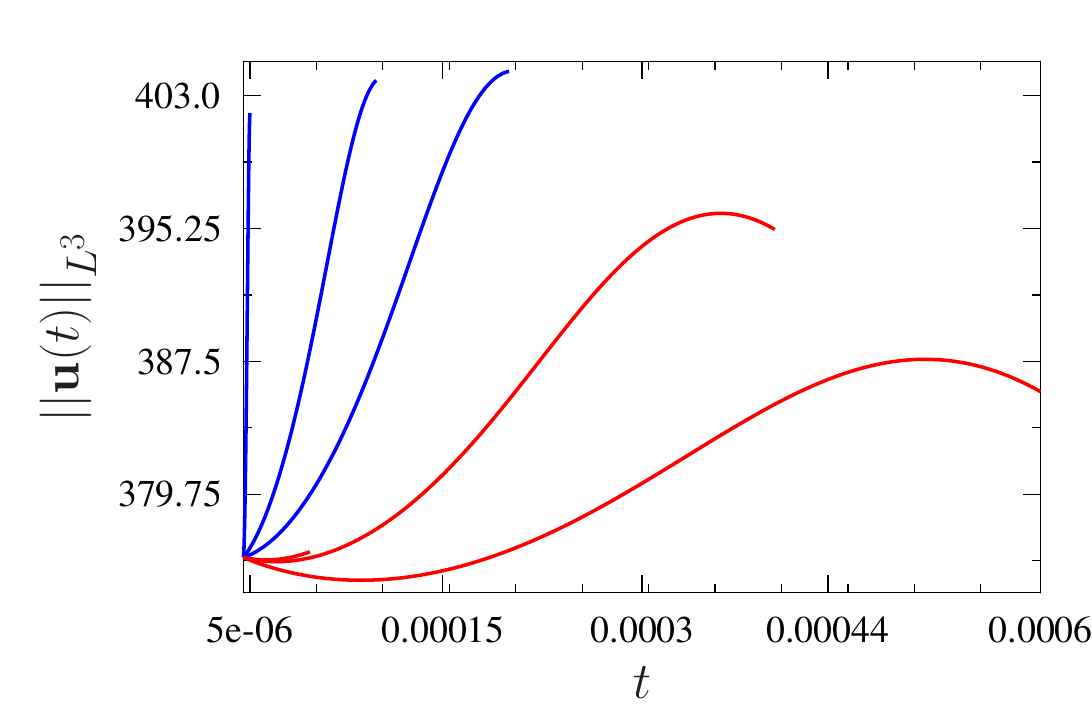}\label{fig:L3ta}}
\subfigure[]{\includegraphics[width=0.49\textwidth]{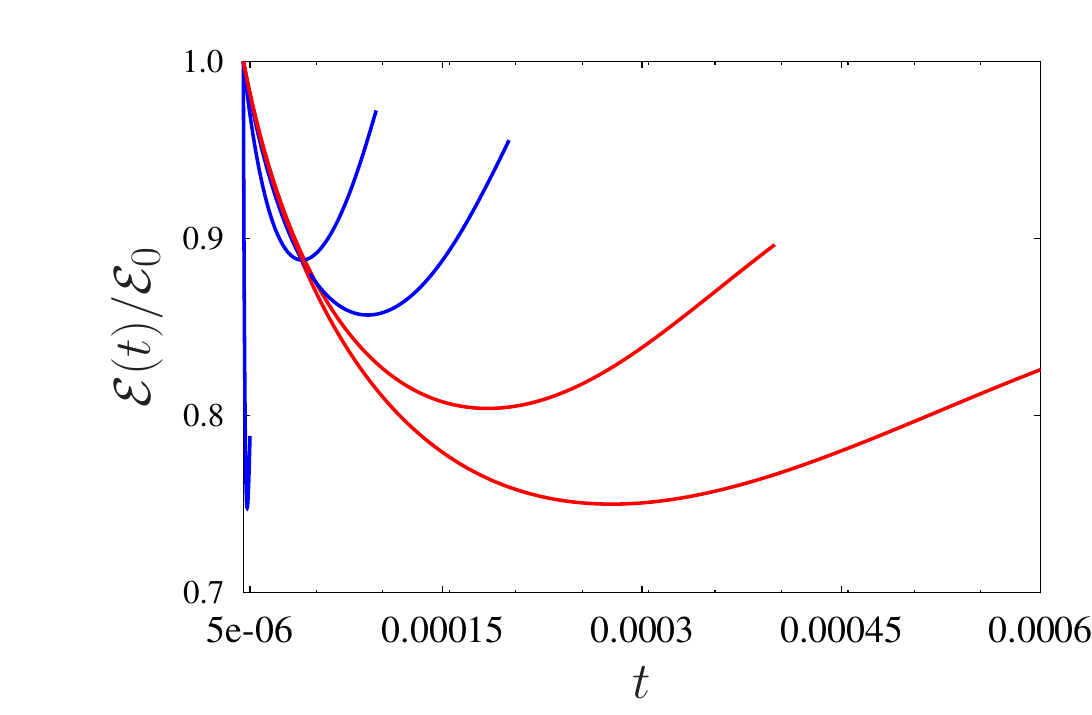}\label{fig:L3Ea}}
}
\mbox{
\subfigure[]{\includegraphics[width=0.49\textwidth]{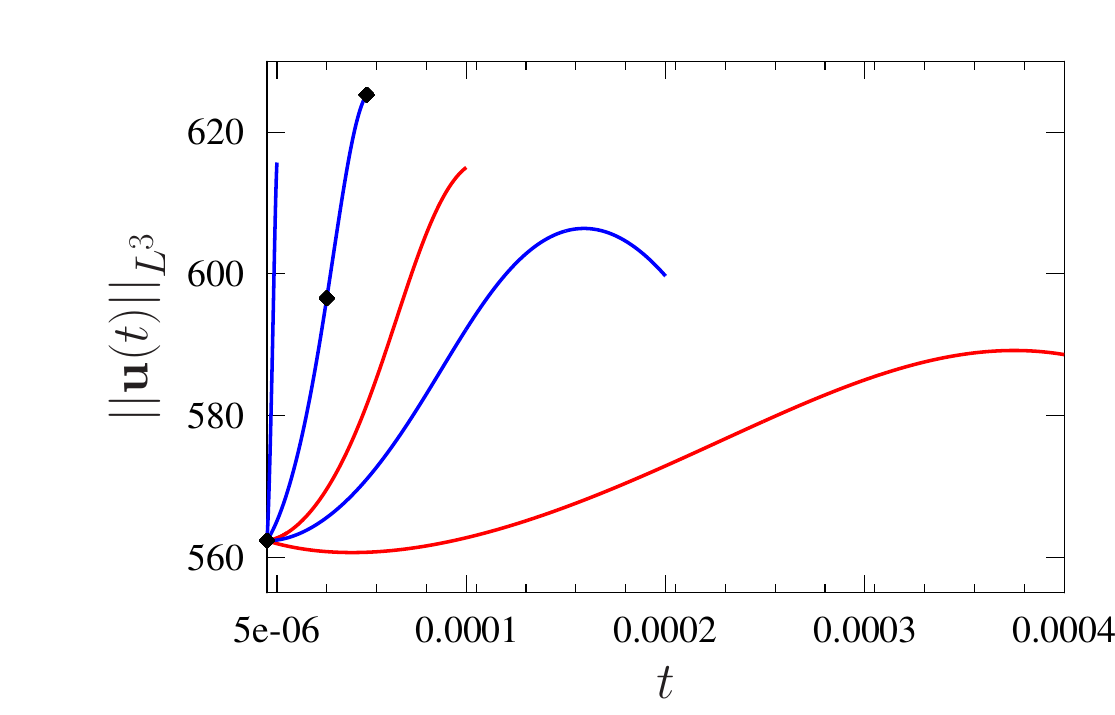}\label{fig:L3tb}}
\subfigure[]{\includegraphics[width=0.49\textwidth]{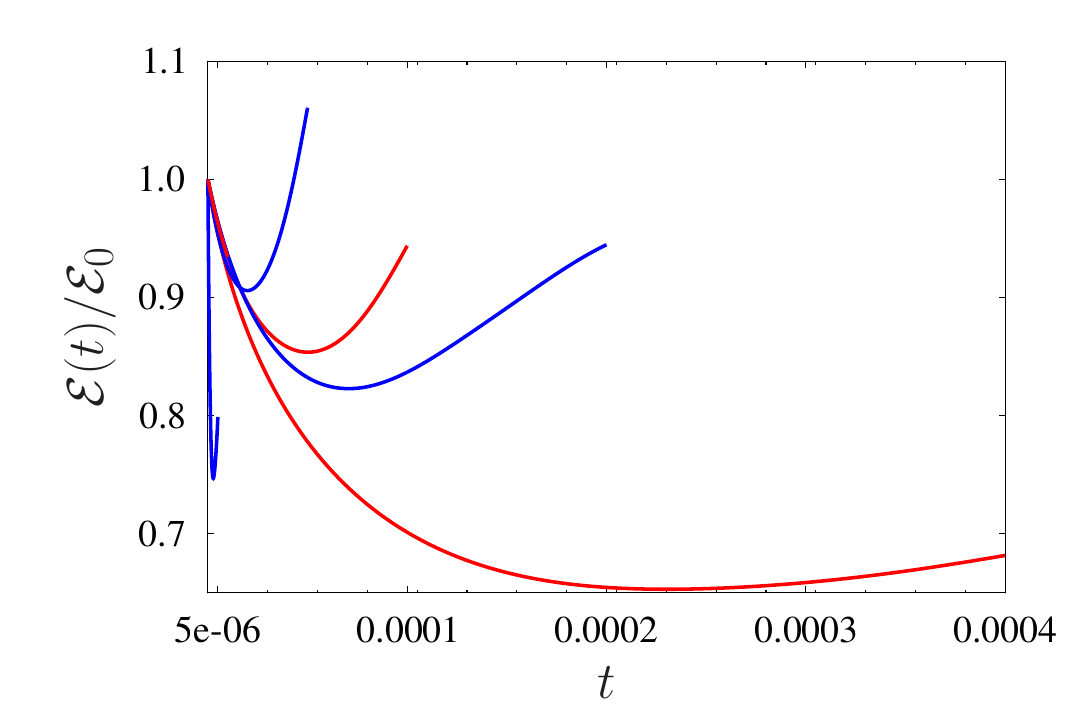}\label{fig:L3Eb}}
}
\mbox{
\subfigure[]{\includegraphics[width=0.49\textwidth]{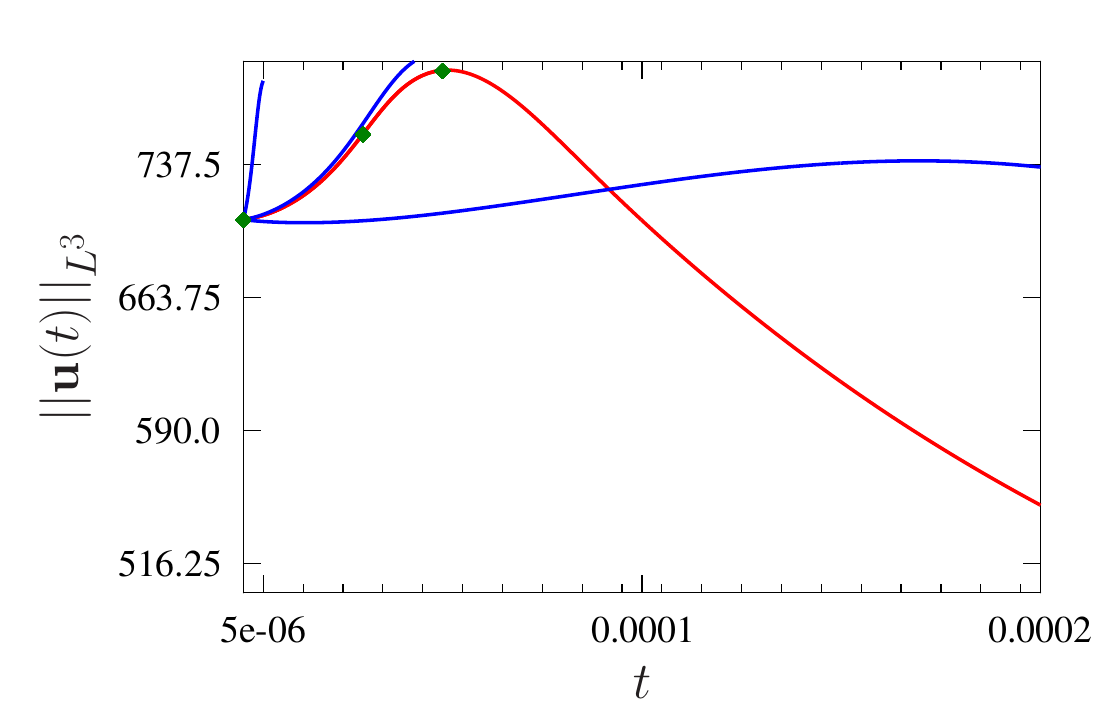}\label{fig:L3tc}}
\subfigure[]{\includegraphics[width=0.49\textwidth]{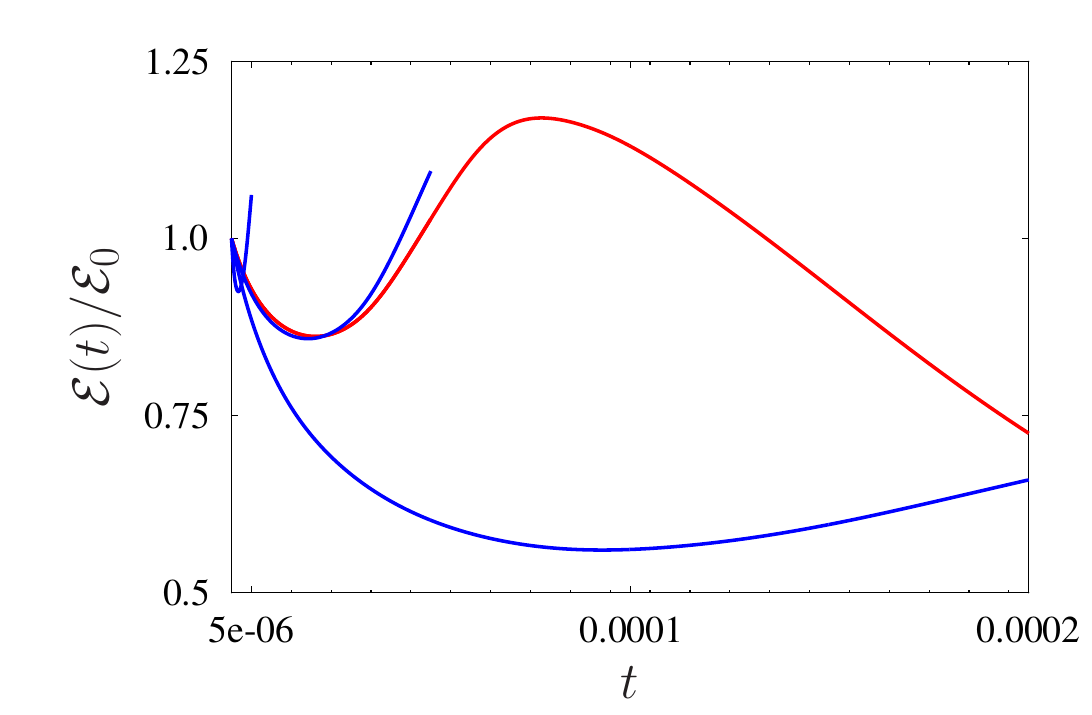}\label{fig:L3Ec}}
}
\caption{[$q=3$] Time evolution of (a,c,e) the norm $||\u(t)||_{L^3}$ 
and (b,d,f) the normalized enstrophy $\E(\u(t)) / \E(0)$ in the 
Navier-Stokes flows with optimal initial conditions obtained by solving 
Problem \ref{pb:PsiHs} (blue) and Problem \ref{pb:PsiLq} (red). The 
values of the constraint parameter are \subref{fig:L3ta}, 
\subref{fig:L3Ea} $B=376.06$, \subref{fig:L3tb}, \subref{fig:L3Eb} 
$B=562.34$ and \subref{fig:L3tc}, \subref{fig:L3Ec} $B=707.10$.  
Solutions are computed over the time window $[0,T]$ where $T$ takes 
three representative values, including the one which gives 
approximately the largest value of $\Psi_{T}(\tuBT)$ for a given value 
of $B$. The solid symbols in \subref{fig:L3tb} and \subref{fig:L3tc} mark 
the time instances at which the flows are visualized in Figures 
\ref{fig:VisualizationTimeEvolSbq3}--\ref{fig:VisualizationTimeEvolLGq3}.}
\label{fig:L3t}
\end{figure}

In this section we consider the limiting case with $q = 3$ and solve 
Problems \ref{pb:PsiHs} and \ref{pb:PsiLq} for three values of the 
constraint parameter $B = 376.06, 562.34, 707.10$ and several time 
windows $T$ chosen empirically to produce the largest growth of 
$||\u(T)||_{L^3}$. The time evolution of the norm $||\u(t)||_{L^3}$ and 
of the normalized enstrophy $\E(\u(t)) / \E(0)$ in the Navier-Stokes 
flows corresponding to the initial conditions $\tuBT$  found in this 
way is shown in Figures \ref{fig:L3t}a,c,e and \ref{fig:L3t}b,d,f, 
respectively, for three representative values of $T$, including the one 
which gives approximately the largest value of $\Psi_{T}(\tuBT)$ for a 
given value of $B$ (in some panels certain curves nearly overlap making 
the impression that some datasets may be missing). We see that the norm 
$||\u(t)||_{L^3}$ exhibits growth at least at some stage in all 
considered cases. However, this growth is rather weak and in relative 
terms is smaller than the growth of the enstrophy $\E(t)$ or the norm 
$||\u(t)||_{L^4}$ achieved in the extreme flows constructed in 
\cite{KangYunProtas2020} and \cite{KangProtas2021}, respectively. This 
is consistent with the observations made in \cite{Hou:22:NS} that, in 
comparison to other quantities characterizing the regularity of 
solutions, the norm $||\u(t)||_{L^3}$  tends to be a slowly growing 
quantity. Even though we maximize the norm at the final time $t = T$, 
cf.~\eqref{eq:Psi}, in some cases the maximum of $||\u(t)||_{L^3}$ is 
achieved at intermediate times $0 < t <T$. We also note that the flows 
obtained by solving Problem \ref{pb:PsiHs} generally lead to a larger 
growth of the norm $||\u(t)||_{L^3}$ than the flows obtained by solving 
Problem \ref{pb:PsiLq}. In all cases the evolution of the normalized 
enstrophy reveals an initial decay followed by some growth.

\subsubsection{Branches of Local Maximizers}
\label{sec:brmL3}

\begin{figure}[t]
\mbox{
\subfigure[]{\includegraphics[width=0.49\textwidth]{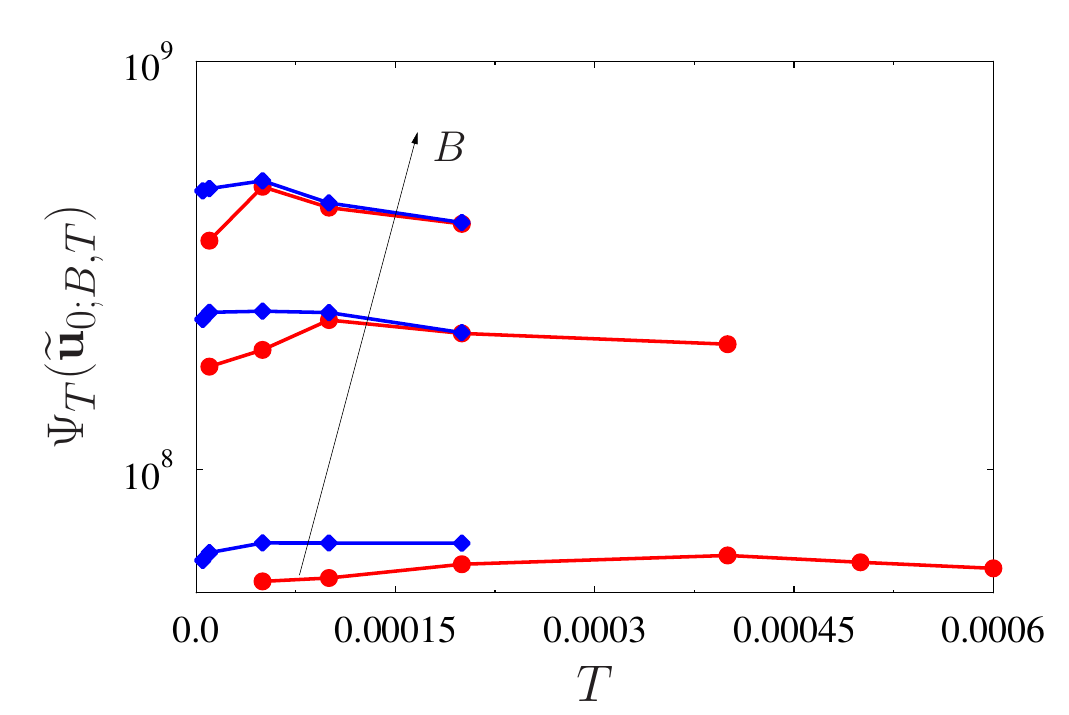}\label{fig:brL3a}}
\subfigure[]{\includegraphics[width=0.49\textwidth]{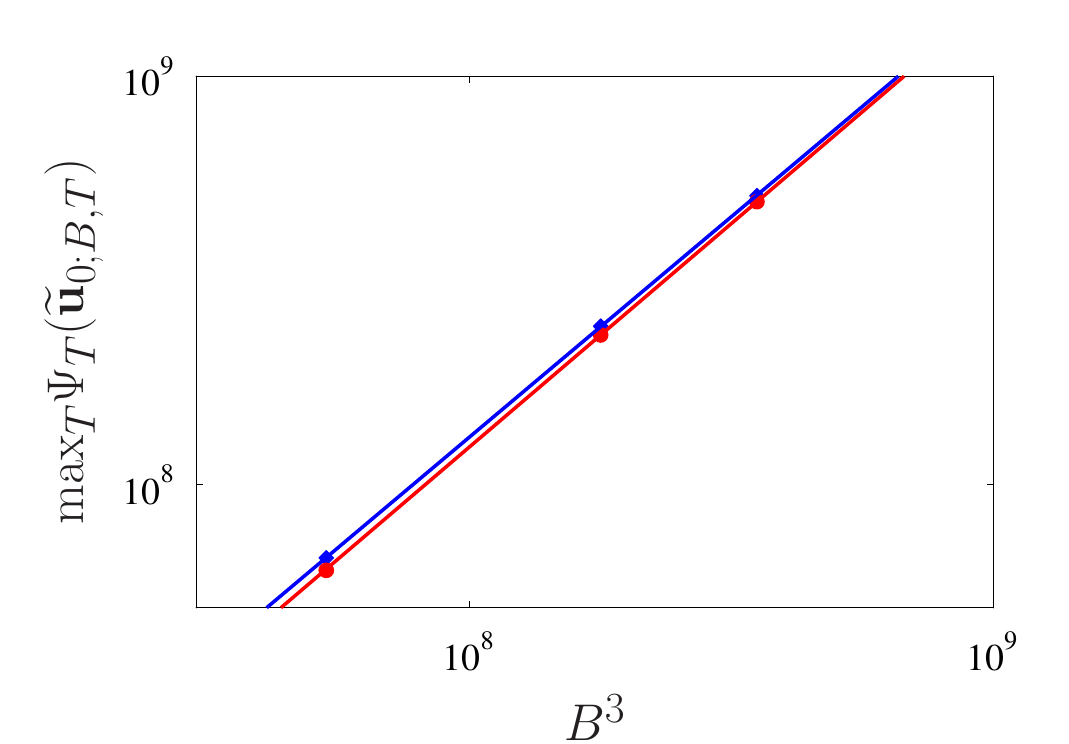}\label{fig:brL3b}}
}
\caption{[$q=3$]  \subref{fig:brL3a} Dependence of the local maxima of 
the objective functional \eqref{eq:Psi} on the length $T$ of the 
optimization window in Problems \ref{pb:PsiHs} (blue) and 
\ref{pb:PsiLq} (red) for different values of the constraint $B=376.06, 
562.34, 707.10$ with the arrow indicating the trend with the increase 
of $B$. \subref{fig:brL3b} Dependence of 
$\max_{T}\,\Psi_{T}\left(\widetilde{\u}_{0;B,T}\right)$ on 
$B^3=||\widetilde{\u}_{0;B,T}||_{L^{3}}^3$ with solid lines representing 
the power-law fits \eqref{eq:sqSb_q3}--\eqref{eq:sqLG_q3}. }
\label{fig:brL3}
\end{figure}

We define a ``branch'' as a family of local maximizers parametrized by 
the length of the time window T and obtained by solving Problem 
\ref{pb:PhiHs}, \ref{pb:PhiLq}, \ref{pb:PsiHs} or \ref{pb:PsiLq} with 
the parameter $B$ fixed. The branches found by solving Problem 
\ref{pb:PsiHs} and Problem \ref{pb:PsiLq} are illustrated by plotting 
the objective functional $\Psi_{T}(\tuBT)$  as a function of $T$ for 
different $B$ in Figure \ref{fig:brL3a}. We note that the dependence of 
$\Psi_{T}(\tuBT)$ on $T$ is rather weak, although a well-defined 
maximum is evident on each branch and shifts towards smaller values of 
$T$ as $B$ increases. The evolutions of $||\u(t)||_{L^3}$ and 
$\E(\u(t)) / \E(0)$ in the flows corresponding to these maxima are also 
included in Figure \ref{fig:L3t}. Figure \ref{fig:brL3a} confirms that 
the values of $||\u(T)||_{L^3}$ in the flows obtained by solving 
Problem \ref{pb:PsiHs}, with optimization performed in the space 
$H^{1/2}(\Omega)$, are generally larger than in the flows obtained by 
solving Problem \ref{pb:PsiLq}, where optimization is performed in the 
space $L^{3}(\Omega)$. It is interesting to quantify how the maximum 
values of $\Psi_{T}(\tuBT)$ attained on different branches scale with 
the constraint parameter $B$ and this relation is shown in Figure 
\ref{fig:brL3b} for flows obtained by solving Problems \ref{pb:PsiHs} 
and \ref{pb:PsiLq}. We see that in both cases there is a clear 
power-law relation whose parameters, the exponent and the prefactor, 
can be determined via a least-squares fit to give
\begin{subequations}
\begin{alignat}{2}
& \text{Problem \ref{pb:PsiHs}:}& \qquad
& \max_{T}\,\Psi_{T}\left(\widetilde{\u}_{0;B,T}\right)\sim 0.31\, 
\left(B^{3}\right)^{1.07},
\label{eq:sqSb_q3} \\
& \text{Problem \ref{pb:PsiLq}:}& \qquad
& \max_{T}\,\Psi_{T}\left(\widetilde{\u}_{0;B,T}\right)\sim 0.22\, 
\left(B^{3}\right)^{1.09}.
\label{eq:sqLG_q3}
\end{alignat}
\end{subequations}
We note that in both cases the exponents are only slightly larger than 
one, meaning that of the increase of $\Psi_{T}(\tuBT)$ comes primarily 
from the increase of the constraint parameter $B$ rather than from the 
nonlinear amplification.

\begin{figure}[t]
\mbox{
\subfigure[]{\includegraphics[width=0.49\textwidth]{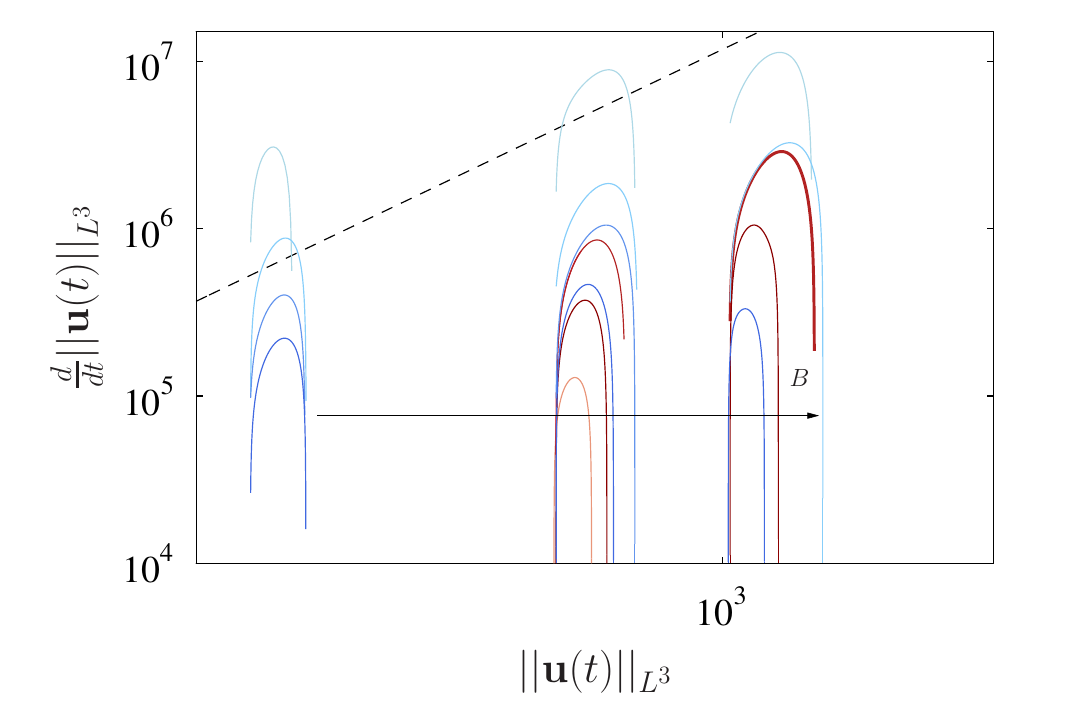}\label{fig:RoCL3a}}
\subfigure[]{\includegraphics[width=0.49\textwidth]{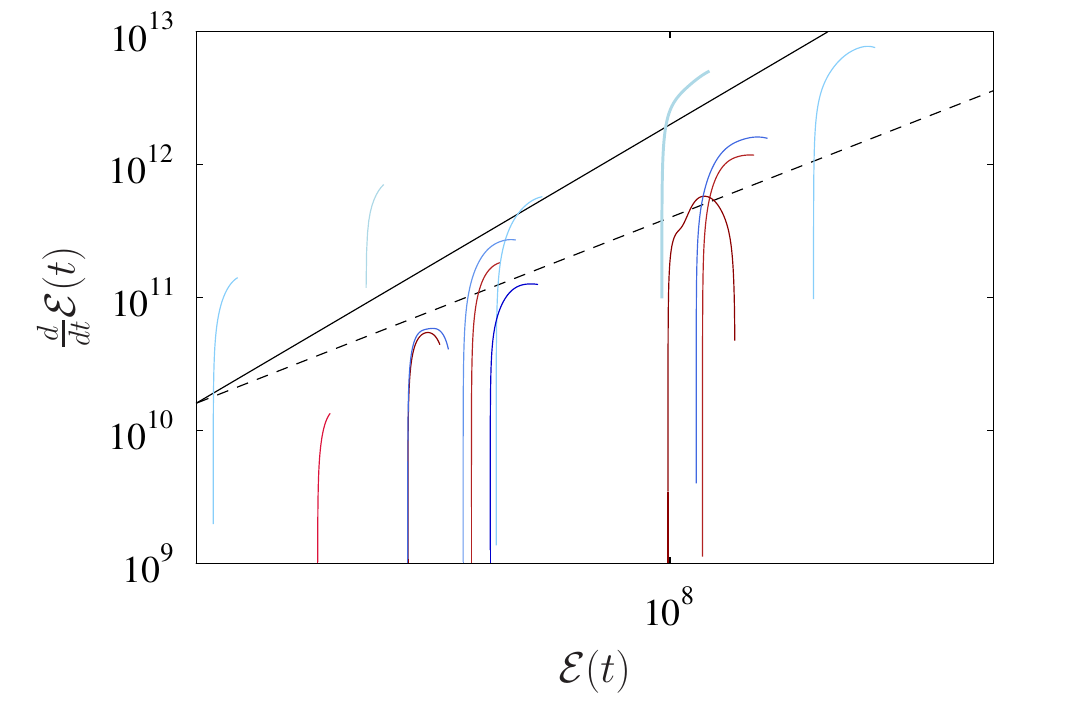}\label{fig:RoCL3b}}
}
\caption{[$q=3$] Navier-Stokes flows corresponding to the optimal 
initial conditions found by solving Problem \ref{pb:PsiHs} and 
\ref{pb:PsiLq} for different $B$ and $T$ shown using the coordinates 
\subref{fig:RoCL3a} 
$\left\{||\u(t)||_{L^{3}},\frac{d}{dt}\,||\u(t)||_{L^{3}}\right\}$ 
and \subref{fig:RoCL3b} $\{\E(t),d\E(t)/dt\}$. In \subref{fig:RoCL3b} 
the black solid line shows the upper bound $d\E/dt\sim\E^3$ on the rate of change of 
enstrophy, cf.~\eqref{eq:RateGrowthEbound}. The dashed lines show the relations 
\subref{fig:RoCL3a} 
$\frac{d}{dt}\,||\u(t)||_{L^{3}}\sim||\u(t)||_{L^{3}}^{5}$ 
from \eqref{eq:noblowup_ub} and \subref{fig:RoCL3b} $d\E/dt\sim\E^{2}$. 
Trajectories marked in blue and red correspond to solutions of Problems 
\ref{pb:PsiHs} and \ref{pb:PsiLq}, respectively. The intensity 
of the color is related to the length of the time window with darker 
colors corresponding to solutions obtained on longer time windows $T$. 
}
\label{fig:RoCL3}
\end{figure}

In order to assess how ``close" the extreme flows studied in Figure 
\ref{fig:L3t} come to forming a singularity, we plot this data using 
the coordinates 
$\left\{||\u(t)||_{L^{3}},\frac{d}{dt}\,||\u(t)||_{L^{3}}\right\}$ 
and $\{\E(t),d\E(t)/dt\}$ in Figures \ref{fig:RoCL3a} and 
\ref{fig:RoCL3b}, respectively. This makes it possible to compare the 
observed behavior with the a priori bounds on the rate of growth of 
$||\u(t)||_{L^{3}}$ and $\E(t)$ from Sections \ref{sec:dEdt} and 
\ref{sec:dLPSdt}. We note that since an upper bound on the rate of 
growth of the norm $||\u(t)||_{L^{3}}$ is not available, here we can 
only consider the ``lower" bound \eqref{eq:noblowup_ub} with $q = 3$. 
The three clusters of trajectories evident in Figures \ref{fig:RoCL3a} 
and \ref{fig:RoCL3b} correspond to the three values of the constraint 
parameter $B$ we use when solving Problems \ref{pb:PsiHs} and 
\ref{pb:PsiLq}. The results in these figures indicate that while in 
some of the flows considered here the norm $||\u(t)||_{L^{3}}$ and the 
enstrophy $\E(t)$ are for some time amplified at a rate consistent with 
singularity formation, this growth is not sustained long enough for 
such singularities to actually form.

\subsubsection{Structure of the Extremal Flows}
\label{sec:vizL3}

\begin{figure}[t]
\centering
\mbox{
\subfigure[${t_0}=0$]{\includegraphics[width=0.35\textwidth]{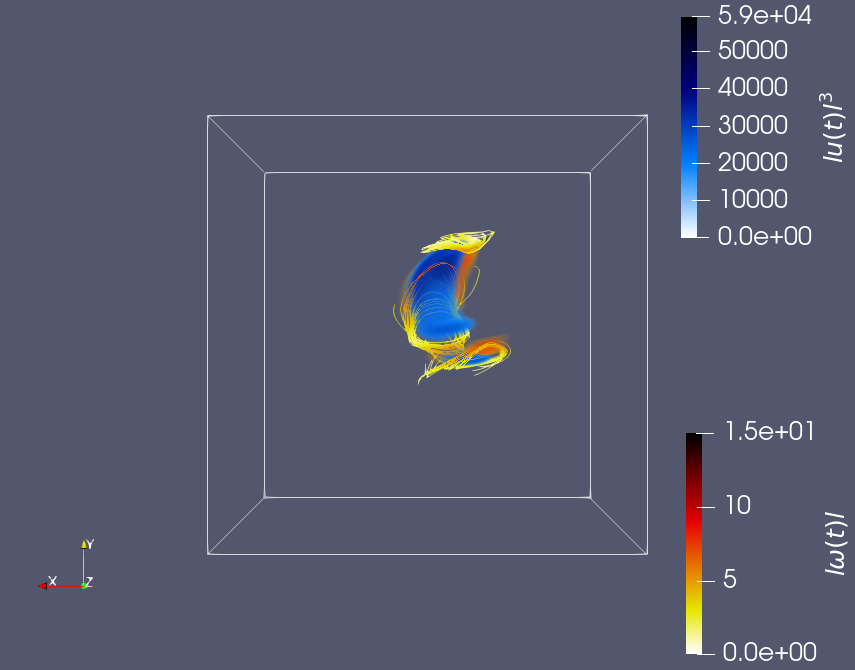}\label{L3SbL1e3T0001_0}}
\quad
\subfigure[${t_1}=6\times10^{-6}$]{\includegraphics[width=0.35\textwidth]{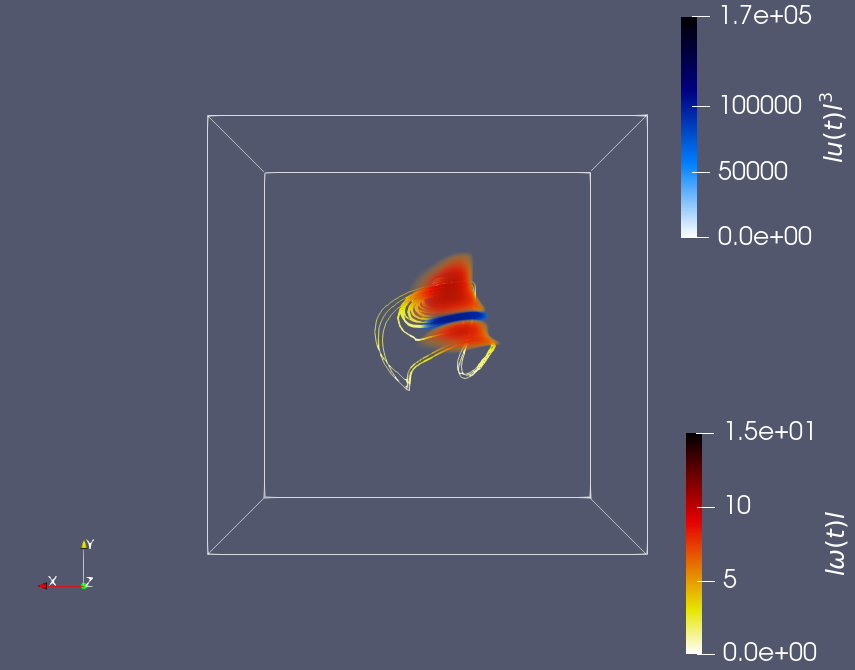}\label{L3SbL1e3T0001_1}}
}
\mbox{
\subfigure[{${t_2}=5\times10^{-5}$}]{\includegraphics[width=0.35\textwidth]{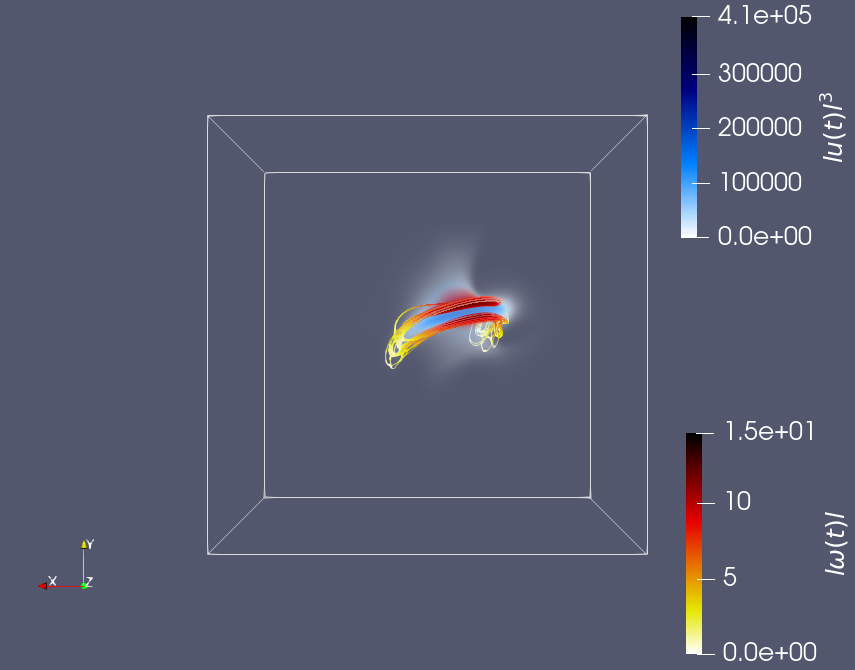}\label{L3SbL1e3T0001_2}}
}
\caption{[$q=3$, Problem \ref{pb:PsiHs}, $B=562.34$ and {$T=5\times10^{-5}$}] 
Snapshots of the magnitude of the vorticity {$|\bomega({t_i}\x)|$} (red 
color scale) in the solution of the Navier-Stokes system \eqref{eq:NS} 
along with vortex lines (red) and the quantity {$|\u(t_{i},\x)|^3$} 
(blue color scale) shown at the times $t_0,...,t_2$. These time 
instances are marked with solid black symbols in Figure \ref{fig:L3tb}.}
\label{fig:VisualizationTimeEvolSbq3}
\end{figure}

\begin{figure}[t]
\centering
\mbox{
\subfigure[${t_0}=0$]{\includegraphics[width=0.35\textwidth]{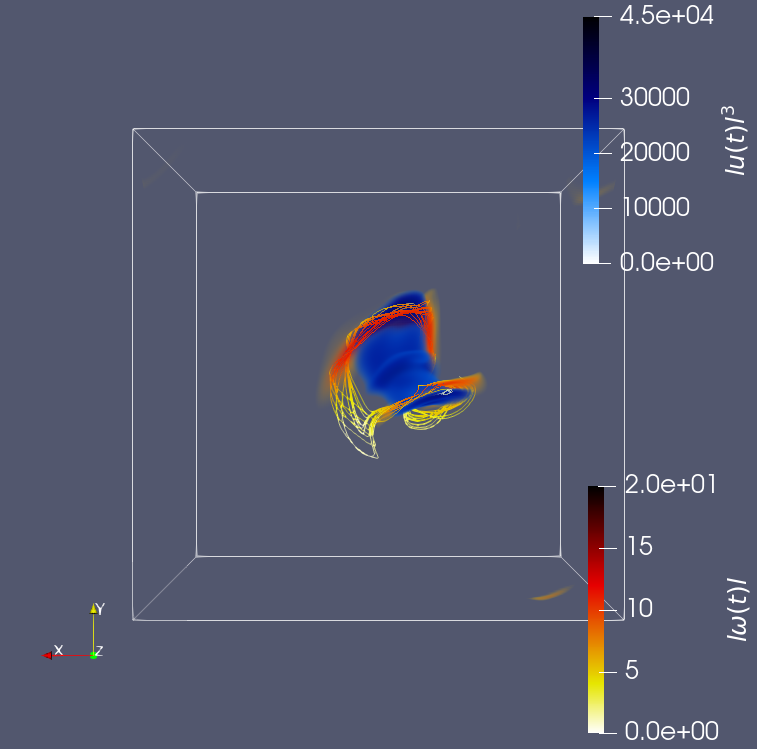}\label{L3LGL25e2T0005_0}}
\quad
\subfigure[${t_1}=3\times10^{-5}$]{\includegraphics[width=0.35\textwidth]{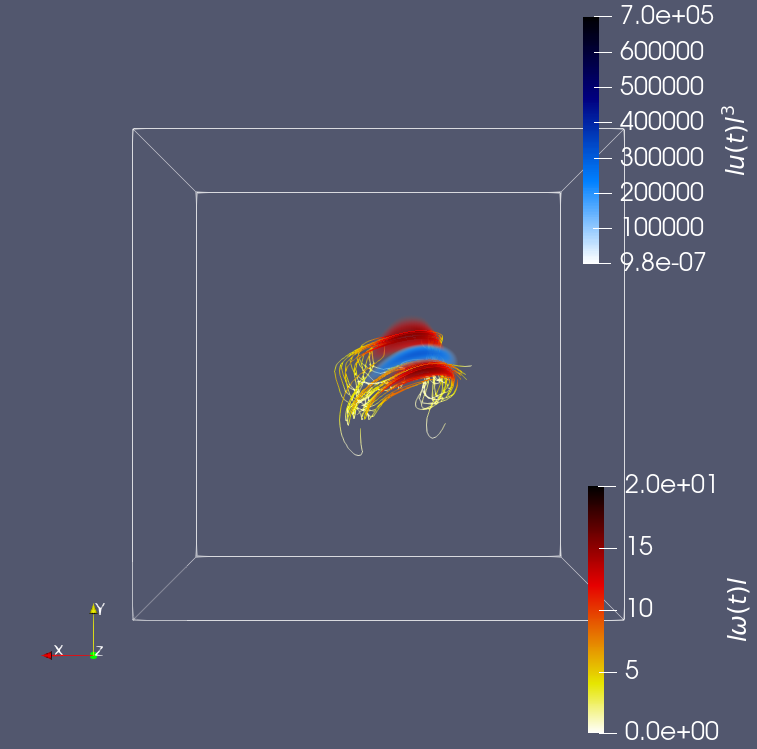}\label{L3LGL25e2T0005_1}}
}
\mbox{
\subfigure[${t_2}=5\times10^{-5}$]{\includegraphics[width=0.35\textwidth]{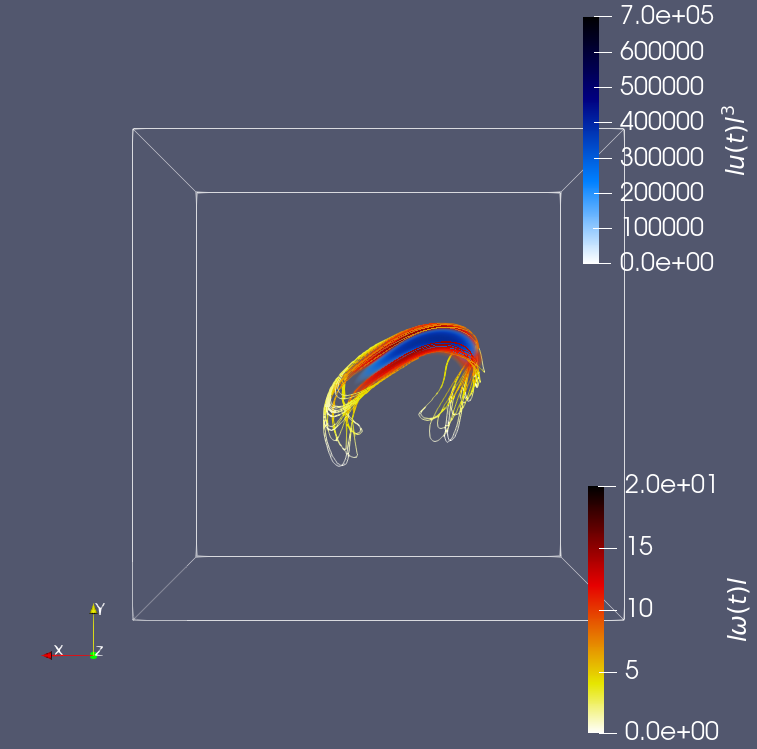}\label{L3LGL25e2T0005_2}}
}
\caption{[$q=3$, Problem \ref{pb:PsiLq}, $B=707.10$ and $T=5\times 
10^{-5}$] Snapshots of the magnitude of the vorticity 
{$|\bomega(t_i,\x)|$} (red color scale) in the solution of the 
Navier-Stokes system \eqref{eq:NS} along with vortex lines (red) and 
the quantity {$|\u(t_i,\x)|^3$} (blue color scale) shown at the times 
$t_0,...,t_2$. These time instances are marked with solid green symbols in 
Figure \ref{fig:L3tc}.}
\label{fig:VisualizationTimeEvolLGq3}
\end{figure}

Finally, we discuss the structure of the extremal flows belonging to 
the different branches shown in Figure \ref{fig:brL3a}. We do this by 
visualizing the vorticity magnitude $|\bomega(t_{i},\x)|$ and the 
quantity $|\u(t_{i},\x)|^3$ in space at different time instances $t_{i}$, 
$i=0,...,3$, defined as follows:
\begin{itemize}
\item $t_{0}$ is the initial time, i.e., $t_0=0,$
\item $t_{1}\approx\argmax_{t\in[0,T]}\frac{d}{dt}\|\u(t)\|_{L^q},$
\item $t_{2}\approx\argmax_{t\in[0,T]}\|\u(t)\|_{L^q},$
\item $t_{3}$ is the final time, i.e., $t_3=T.$
\end{itemize}
These plots are shown in Figures \ref{fig:VisualizationTimeEvolSbq3} 
and \ref{fig:VisualizationTimeEvolLGq3} for flows obtained with the 
optimal initial conditions $\tuBT$ found by solving Problems 
\ref{pb:PsiHs} and \ref{pb:PsiLq} with different values of the 
constraint parameter $B$ and $T$ corresponding to the maxima over the 
branches shown in Figure \ref{fig:brL3a}. We note that here $t_{2} 
\approx t_{3}$, so the flow state at $t = t_{3}$ is not separately 
visualized (however, this will be done for the flows corresponding to $q 
= 9$ in Section \ref{sec:vizL9}). We observe that the flow structures 
corresponding to solutions of Problems \ref{pb:PsiHs} and 
\ref{pb:PsiLq}  are distinct and are both strongly localized, in 
contrast to what was found when solving Problem \ref{pb:maxET} in 
\cite{KangYunProtas2020}. They resemble the flow structures found by 
solving Problem \ref{pb:PsiHs} with $q = 4$ in \cite{KangProtas2021} 
and involve bent vortex rings which undergo deformation as time goes on. The 
maxima of $|\u(T,\x)|^3$, which is what contributes the most to 
the value of the functional $\Psi_{T}(\tuBT)$, occur inside the opening 
in the vortex ring as it is stretched.


\FloatBarrier

\subsection{Extremal Flows in $L^9(\Omega)$}
\label{sec:L9}

\begin{figure}[H]
\mbox{
\subfigure[]{\includegraphics[width=0.49\textwidth]{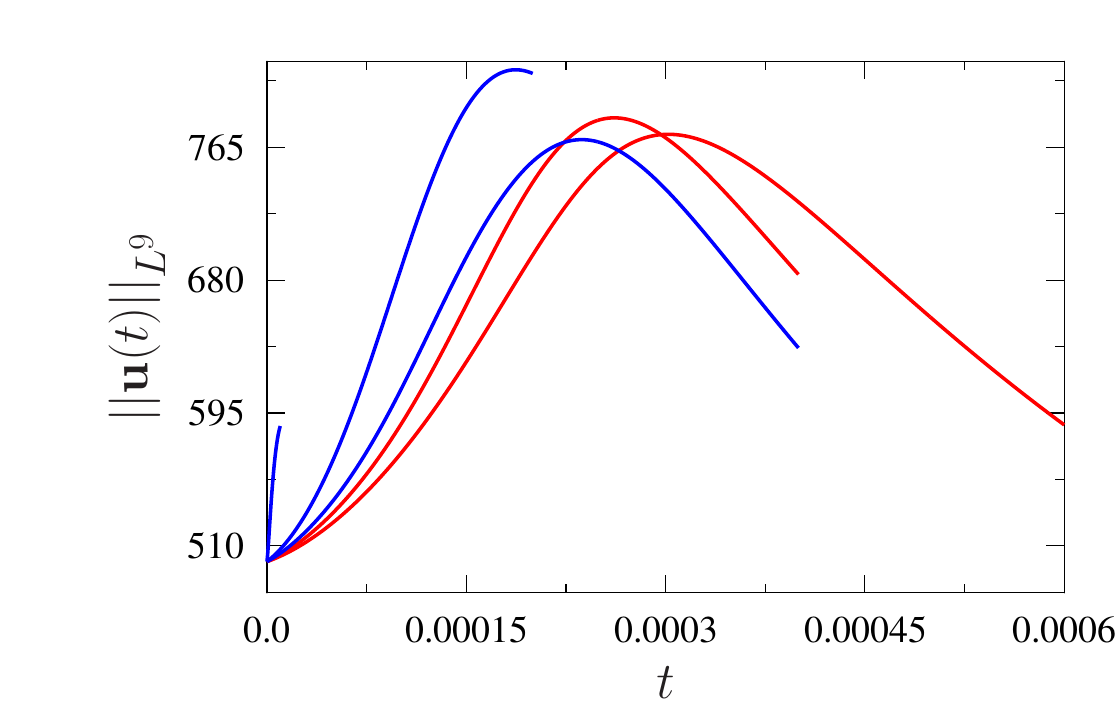}\label{fig:L9ta}}
\subfigure[]{\includegraphics[width=0.49\textwidth]{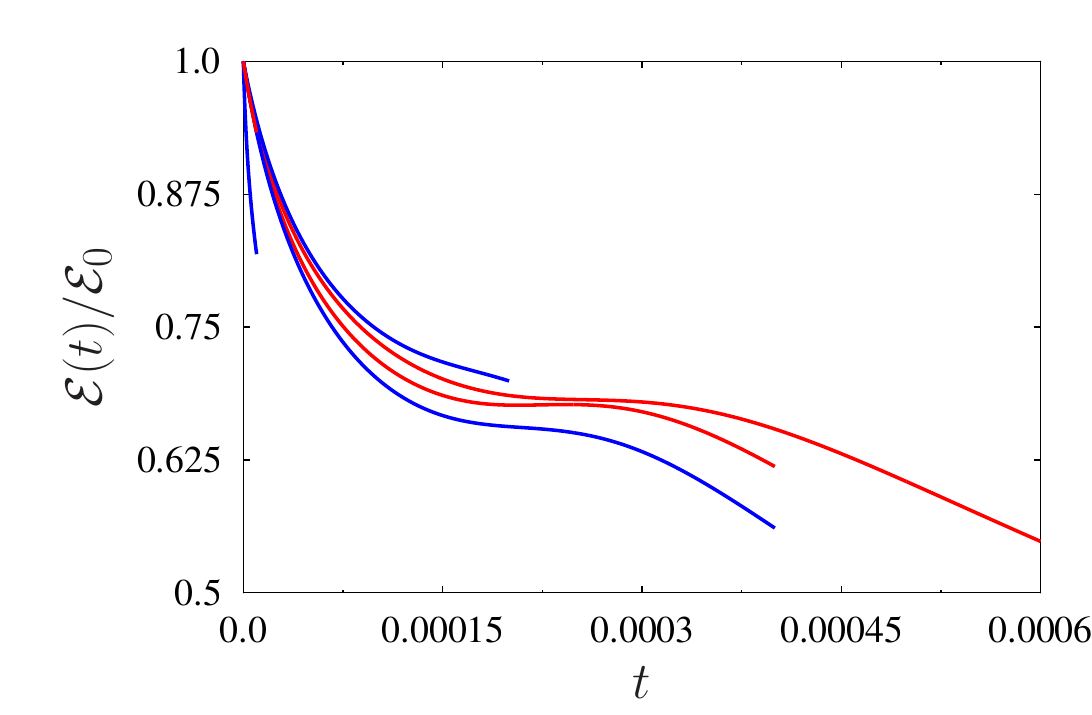}\label{fig:L9Ea}}
}
\mbox{
\subfigure[]{\includegraphics[width=0.49\textwidth]{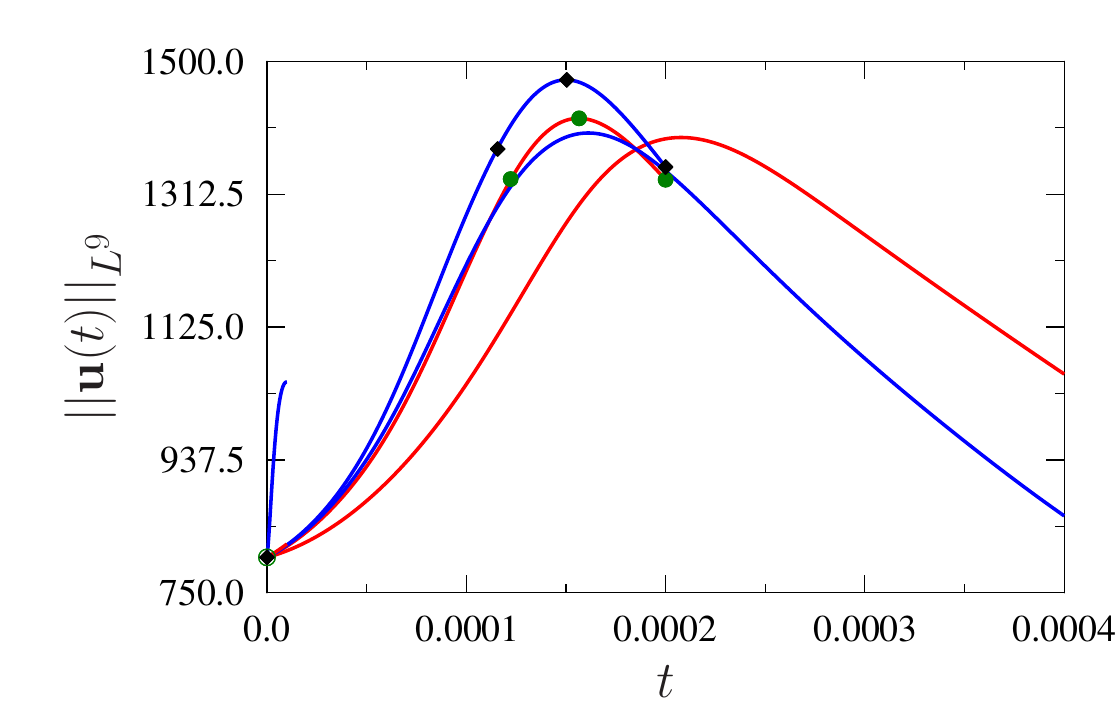}\label{fig:L9tb}}
\subfigure[]{\includegraphics[width=0.49\textwidth]{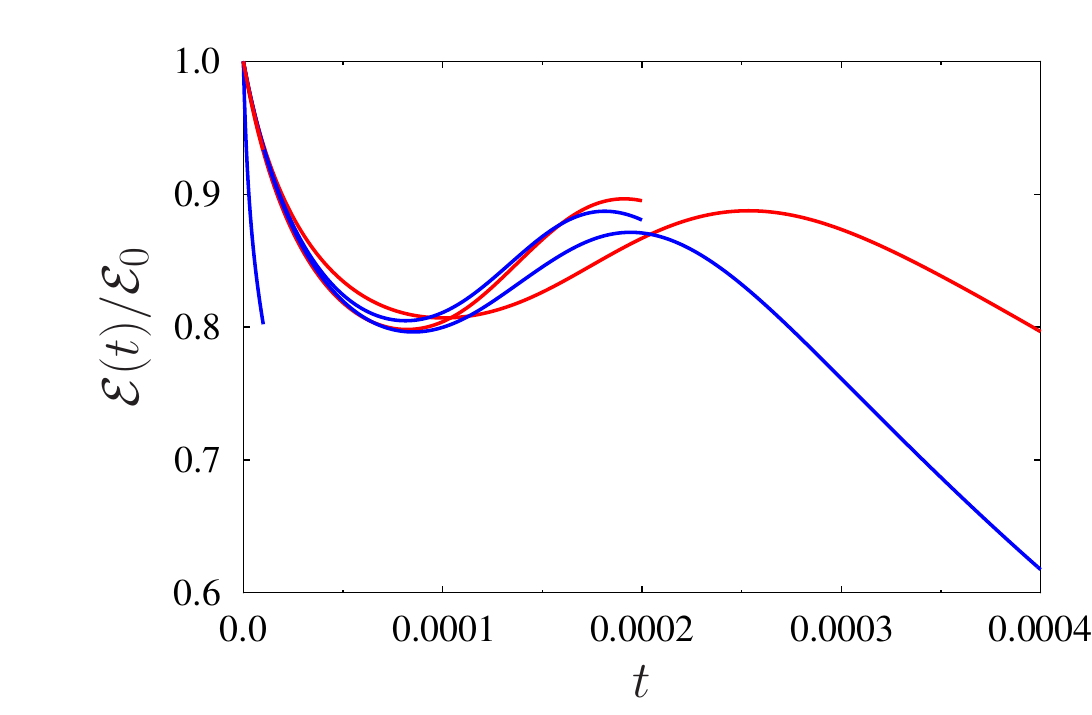}\label{fig:L9Eb}}
}
\mbox{
\subfigure[]{\includegraphics[width=0.49\textwidth]{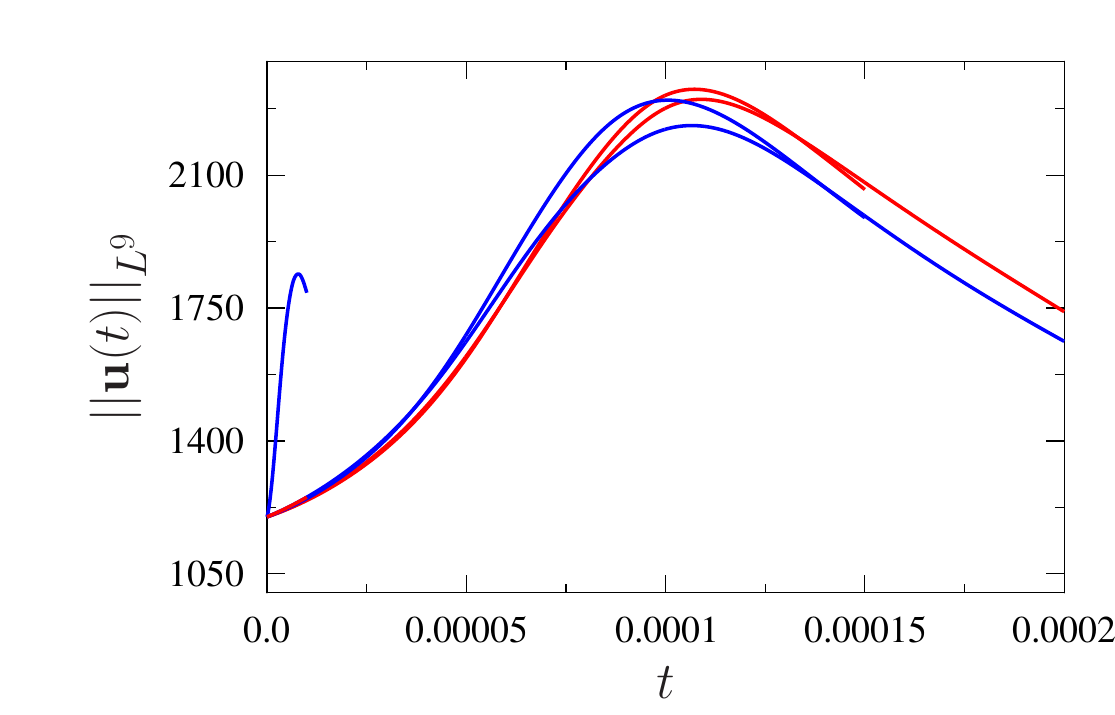}\label{fig:L9tc}}
\subfigure[]{\includegraphics[width=0.49\textwidth]{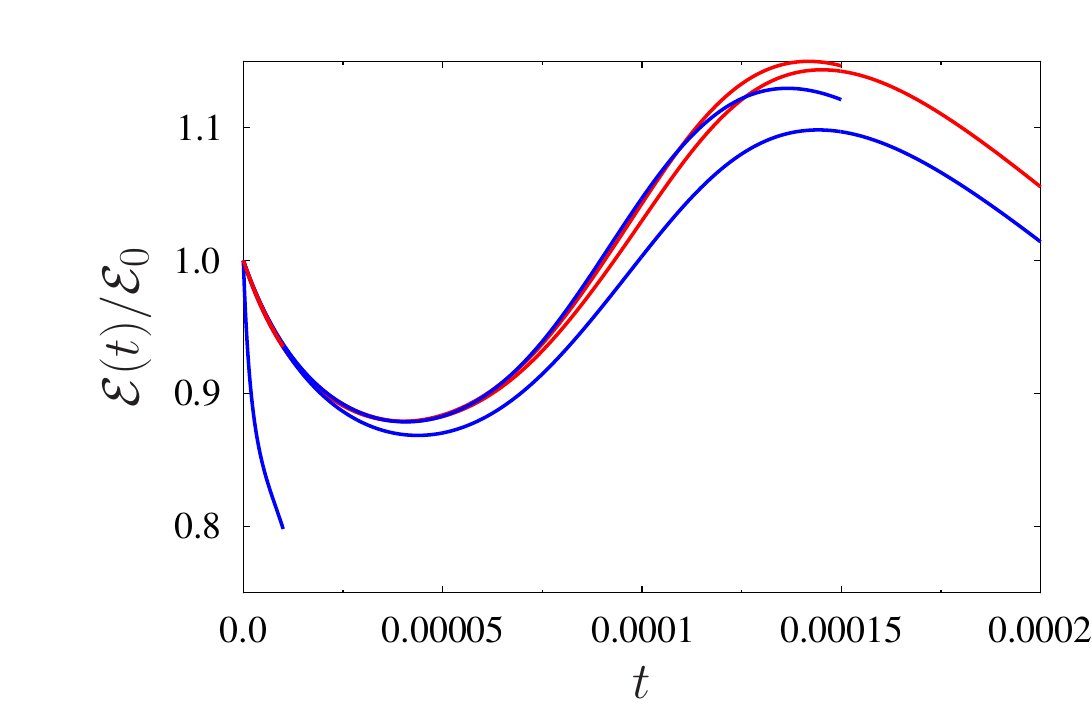}\label{fig:L9Ec}}
}
\caption{[$q=9$] Time evolution of (a,c,e) the norm $||\u(t)||_{L^9}$ 
and (b,d,f) the normalized enstrophy $\E(\u(t)) / \E(0)$ in the 
Navier-Stokes flows with optimal initial conditions obtained by solving 
Problem \ref{pb:PhiHs} (blue) and Problem \ref{pb:PhiLq} (red). The 
values of the constraint parameter are \subref{fig:L3ta}, 
\subref{fig:L3Ea} $B=500$, \subref{fig:L3tb} , \subref{fig:L3Eb} 
$B=800$ and \subref{fig:L3tc}, \subref{fig:L3Ec} $B=1200$.  Solutions 
are computed over the time window $[0,T]$ where $T$ takes three 
representative values, including the one which gives approximately the 
largest value of $\Phi^{9}_{T}(\tuBT)$ for a given value of $B$.
The solid symbols in \subref{fig:L9tb} mark 
the time instances at which the flows are visualized in Figures 
\ref{fig:VisualizationTimeEvolSbq9}--\ref{fig:VisualizationTimeEvolLGq9}}.
\label{fig:L9t}
\end{figure}

In this section we consider the case with $q = 9$ and solve Problems 
\ref{pb:PhiHs} and \ref{pb:PhiLq} for three values of the constraint 
parameter $B = 500, 800, 1200$ and 
several time windows $T$ chosen empirically to produce the largest 
optimized values of the objective functional $\Phi^{9}_{T}(\tuBT)$, 
cf.~\eqref{eq:Phi}. The time evolution of the norm $||\u(t)||_{L^9}$ 
and of the normalized enstrophy $\E(\u(t)) / \E(0)$ in the 
Navier-Stokes flows corresponding to the initial conditions $\tuBT$ 
found in this way is shown in Figures \ref{fig:L9t}a,c,e and 
\ref{fig:L9t}b,d,f, respectively, for three representative values of 
$T$, including the one which gives approximately the largest value of 
$\Phi^{9}_{T}(\tuBT)$ for a given value of $B$. In contrast to the 
results obtained by solving Problems \ref{pb:PsiHs} and \ref{pb:PsiLq}, 
cf.~Figure \ref{fig:L3t}, now the norm $||\u(t)||_{L^9}$ reveals a much 
larger relative growth in time which is comparable in flows obtained by 
solving the two optimization problems. On the other hand, the growth of 
the enstrophy $\E(t)$ is weaker, and even entirely absent in the case 
of the smallest value of $B$ (cf.~Figure \ref{fig:L9Ea}), as compared 
to the flows obtained by solving Problems \ref{pb:PsiHs} and 
\ref{pb:PsiLq}.

\subsubsection{Branches of Local Maximizers}
\label{sec:brmL9}

\begin{figure}[t]
\mbox{
\subfigure[]{\includegraphics[width=0.49\textwidth]{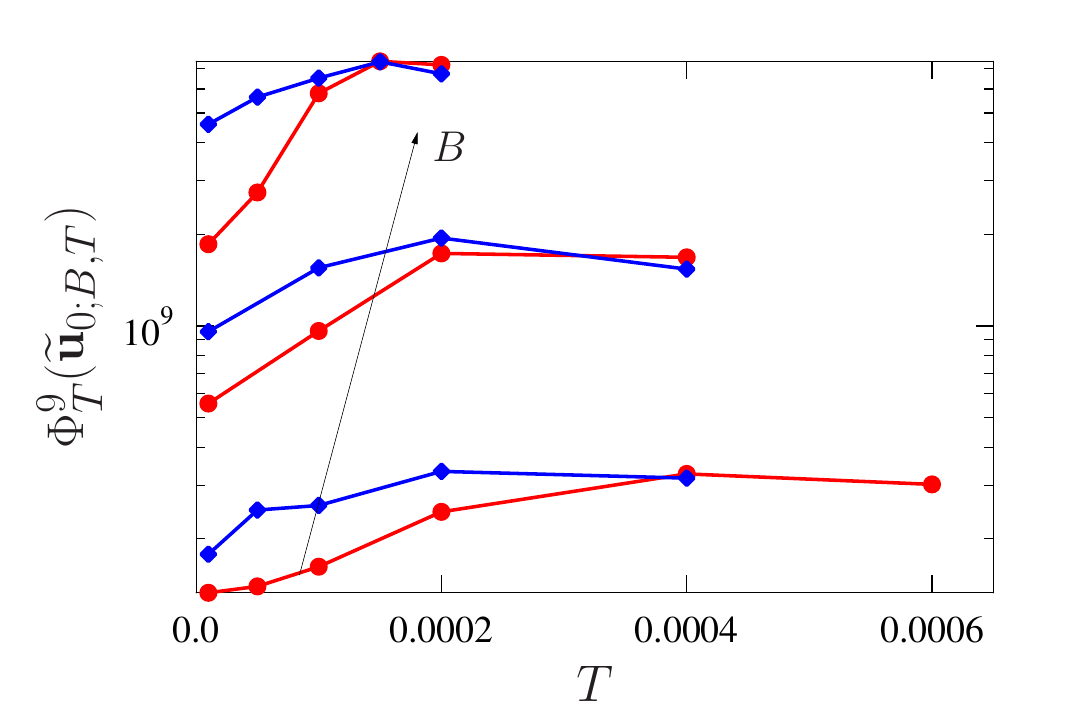}\label{fig:brL9a}}
\subfigure[]{\includegraphics[width=0.49\textwidth]{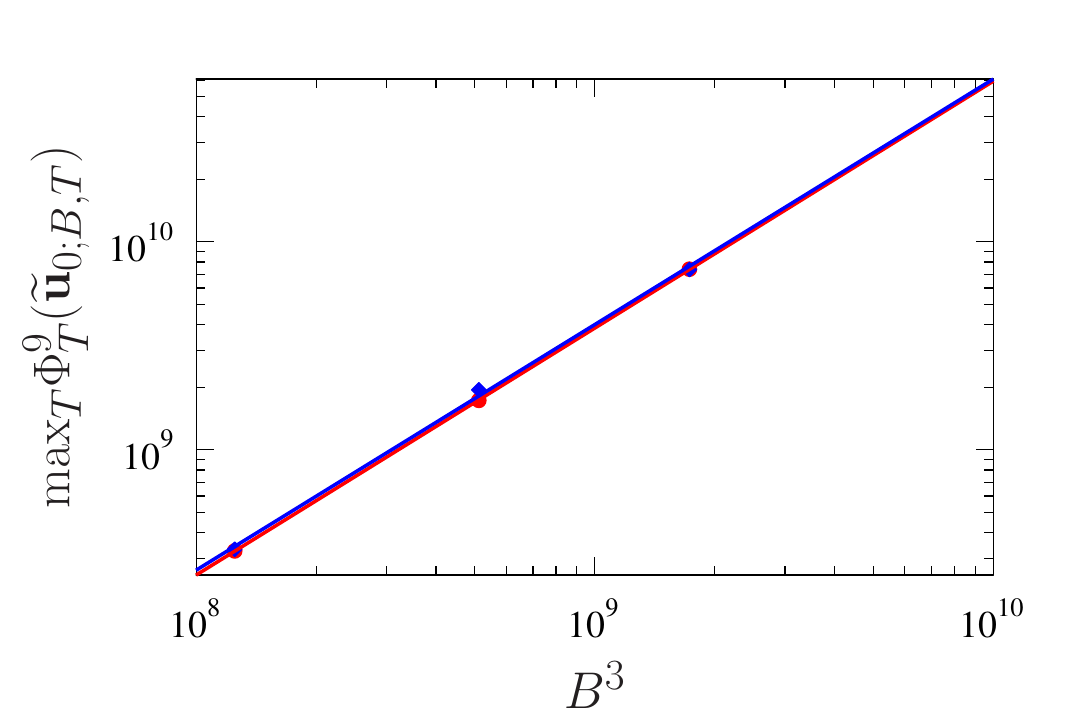}\label{fig:brL9b}}
}
\caption{[$q=9$]  \subref{fig:brL9a} Dependence of the local maxima of 
the objective functional $\Phi^{9}_{T}(\tuBT)$ on the length $T$ of the 
optimization window in Problems \ref{pb:PhiHs} (blue) and 
\ref{pb:PhiLq} (red) for different values of the constraint $B=500 
,800, 1200$ with the arrow indicating the trend with the increase of 
$B$. \subref{fig:brL9b} Dependence of 
$\max_{T}\,\Phi^{9}_{T}\left(\widetilde{\u}_{0;B,T}\right)$ on 
$B^3=||\widetilde{\u}_{0;B,T}||_{L^9}^3$ with solid lines 
representing the power-law fits \eqref{eq:sqSb_q9}--\eqref{eq:sqLG_q9}.}
\label{fig:brL9}
\end{figure}

The branches of local maximizers found by solving Problems 
\ref{pb:PhiHs} and \ref{pb:PhiLq} are illustrated by plotting the 
objective functional $\Phi^{9}_{T}(\tuBT)$  as a function of $T$ for 
different $B$ in Figure \ref{fig:brL9a}. In contrast to the results 
obtained for $q = 3$, cf.~Figure \ref{fig:brL3a}, now the dependence of 
$\Phi^{9}_{T}(\tuBT)$  on $T$ is more pronounced and we still observe a 
well-defined global maximum on each branch. The values of the objective 
functional $\Phi^{9}_{T}(\tuBT)$ corresponding to these maxima are 
comparable for the extreme flows obtained by solving Problems 
\ref{pb:PhiHs} and \ref{pb:PhiLq}. Figure \ref{fig:brL9b} shows that 
the dependence of the maximum 
values of $\Psi_{T}(\tuBT)$ attained on different branches on the constraint  
parameter $B$ has the form of power laws of the following form determined 
 via least-squares fits
\begin{subequations}
\begin{alignat}{2}
& \text{Problem \ref{pb:PhiHs}:}& \qquad
& \max_{T}\,\Phi_{T}^9\left(\widetilde{\u}_{0;B,T}\right)\sim 0.09\, 
\left(B^{3}\right)^{1.18},
\label{eq:sqSb_q9} \\
& \text{Problem \ref{pb:PhiLq}:}& \qquad
& \max_{T}\,\Phi_{T}^9\left(\widetilde{\u}_{0;B,T}\right)\sim 0.07\, 
\left(B^{3}\right)^{1.19}.
\label{eq:sqLG_q9}
\end{alignat}
\end{subequations}
The quantity $B^{3}$ was chosen as the independent variable in these 
fits recognizing the fact that with $q = 9$ the objective functional 
\eqref{eq:Phi} scales with the 3rd power of the velocity. We observe 
that as compared to the fits obtained for $q = 3$, cf.~Figure 
\ref{fig:brL3b} and relations \eqref{eq:sqSb_q3}--\eqref{eq:sqLG_q3}, 
now the exponents are noticeably larger reflecting stronger effects of 
nonlinear amplification.

\begin{figure}[t]
\mbox{
\subfigure[]{\includegraphics[width=0.49\textwidth]{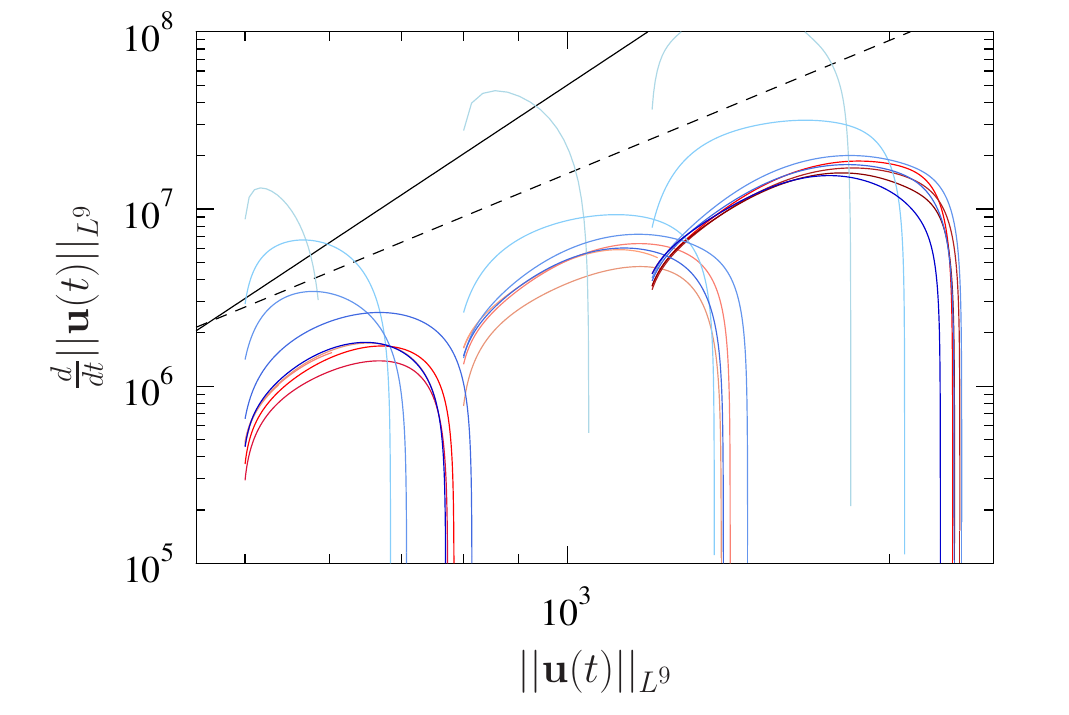}\label{fig:RoCL9a}}
\subfigure[]{\includegraphics[width=0.49\textwidth]{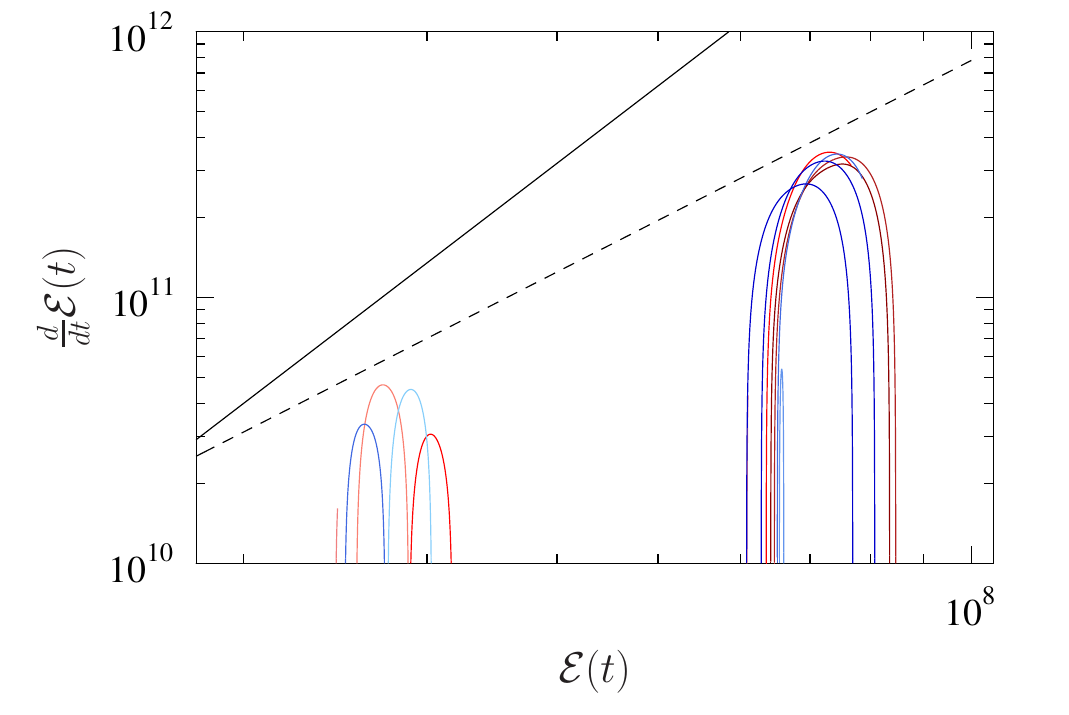}\label{fig:RoCL9b}}
}
\caption{[$q=9$] Navier-Stokes flows corresponding to the optimal 
initial conditions found by solving Problem \ref{pb:PhiHs} and 
\ref{pb:PhiLq} for different $B$ and $T$ shown using the coordinates 
\subref{fig:RoCL9a} 
$\left\{||\u(t)||_{L^{9}},\frac{d}{dt}\,||\u(t)||_{L^{9}}\right\}$ 
and \subref{fig:RoCL9b} $\{\E(t),d\E(t)/dt\}$. The black solid lines 
represent the upper bounds  \subref{fig:RoCL9a} 
$\frac{d}{dt}\,||\u(t)||_{L^{9}}\sim||\u(t)||_{L^{9}}^4$ 
from \eqref{eq:aprioriLq} and \subref{fig:RoCL9b}  $d\E/dt\sim\E^3$ 
from \eqref{eq:RateGrowthEbound}. The dashed lines show the relations 
\subref{fig:RoCL9a} 
$\frac{d}{dt}\,||\u(t)||_{L^{9}}\sim||\u(t)||_{L^{9}}^{5/2}$ 
from \eqref{eq:noblowup_ub} and \subref{fig:RoCL9b} $d\E/dt\sim\E^{2}$. 
Trajectories marked in blue and red correspond to solutions of Problems 
\ref{pb:PhiHs} and \ref{pb:PhiLq}, respectively. The intensity of the 
color is related to the length of the time window with darker colors 
corresponding to solutions obtained on longer time windows $T$.}
\label{fig:RoCL9}
\end{figure}

In order to assess how ``close" the extreme flows studied in Figure 
\ref{fig:L9t} come to forming a singularity, we plot this data using 
the coordinates 
$\left\{||\u(t)||_{L^9},\frac{d}{dt}\,||\u(t)||_{L^9}\right\}$ 
and $\{\E(t),d\E(t)/dt\}$ in Figures \ref{fig:RoCL9a} and 
\ref{fig:RoCL9b}, respectively. This makes it possible to compare the 
observed behavior with the a priori bounds on the rate of growth of 
$||\u(t)||_{L^{9}}$ and $\E(t)$ from Sections \ref{sec:dEdt} and 
\ref{sec:dLPSdt}. Since the enstrophy $\E(t)$ is a decaying function of 
time for all flows corresponding to the smallest value of the 
constraint parameter $B = 500$, cf.~Figure \ref{fig:L9Ea}, there are 
only two clusters of trajectories in Figure \ref{fig:RoCL9b}. Figure 
\ref{fig:RoCL9a} shows that in the flows found by solving Problems 
\ref{pb:PhiHs} and \ref{pb:PhiLq} the the norm $||\u(t)||_{L^{9}}$ is 
amplified at a rate consistent with singularity formation for a longer 
time than was the case for the norm $||\u(t)||_{L^{3}}$ in the flows 
found by solving Problems \ref{pb:PsiHs} and \ref{pb:PsiLq}, cf.~Figure 
\ref{fig:RoCL3a}. However, this process is again not sustained nearly 
long enough  for a singularity to form.

\subsubsection{Structure of the Extremal Flows}
\label{sec:vizL9}

\begin{figure}[H]
\centering
\mbox{
\subfigure[${t_0}=0$]{\includegraphics[width=0.35\textwidth]{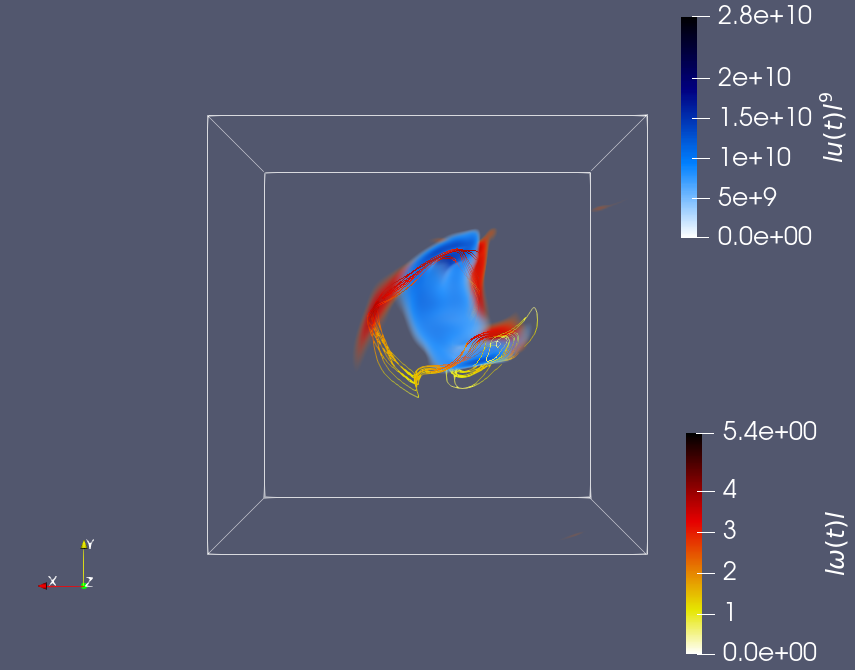}\label{L9SbL1e3T002_0}}
\quad
\subfigure[{${t_1}=1.1\times10^{-4}$}]{\includegraphics[width=0.35\textwidth]{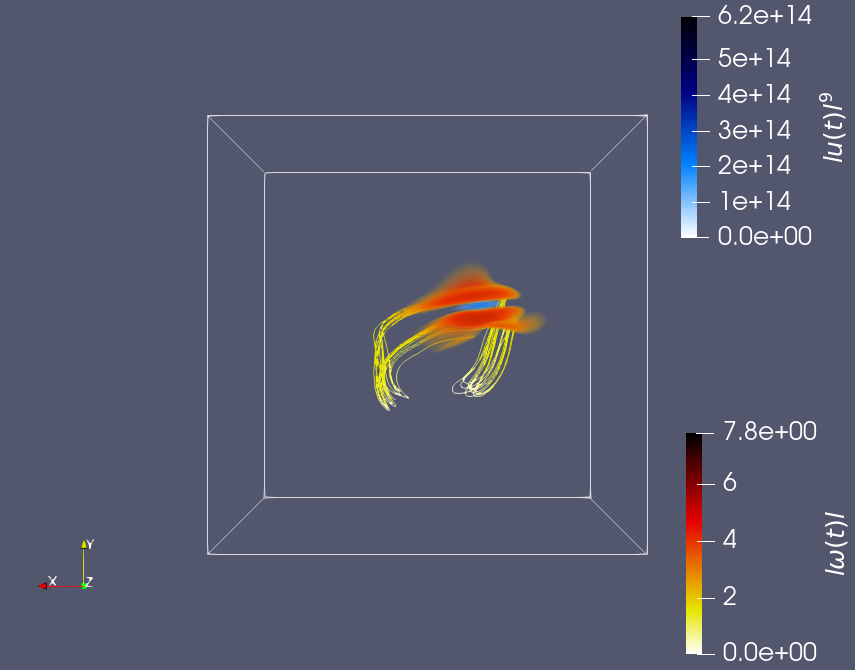}\label{L9SbL1e3T002_1}}
}
\mbox{
\subfigure[{${t_2}=1.5\times10^{-4}$}]{\includegraphics[width=0.35\textwidth]{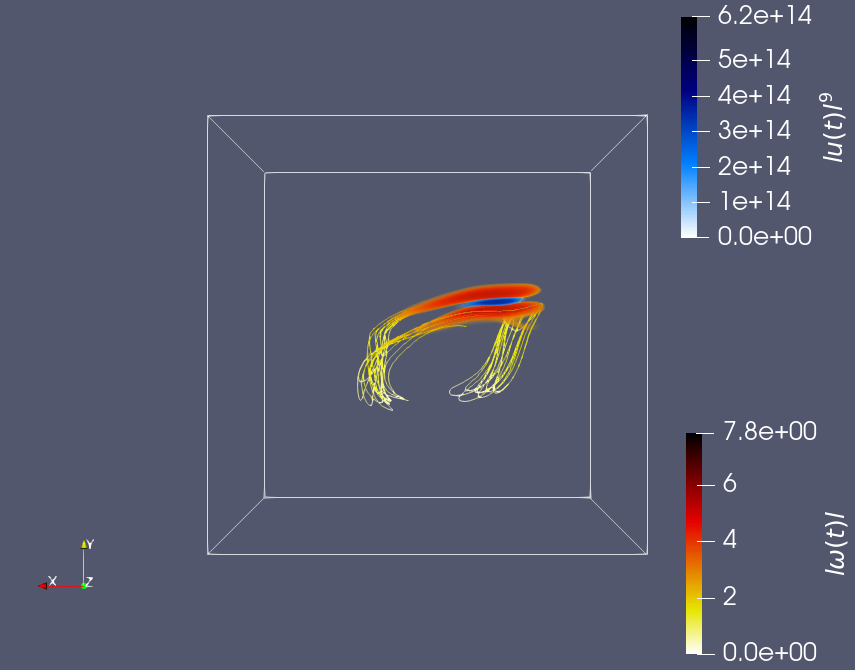}\label{L9SbL1e3T002_2}}
\quad
\subfigure[{${t_3}=2\times10^{-4}$}]{\includegraphics[width=0.35\textwidth]{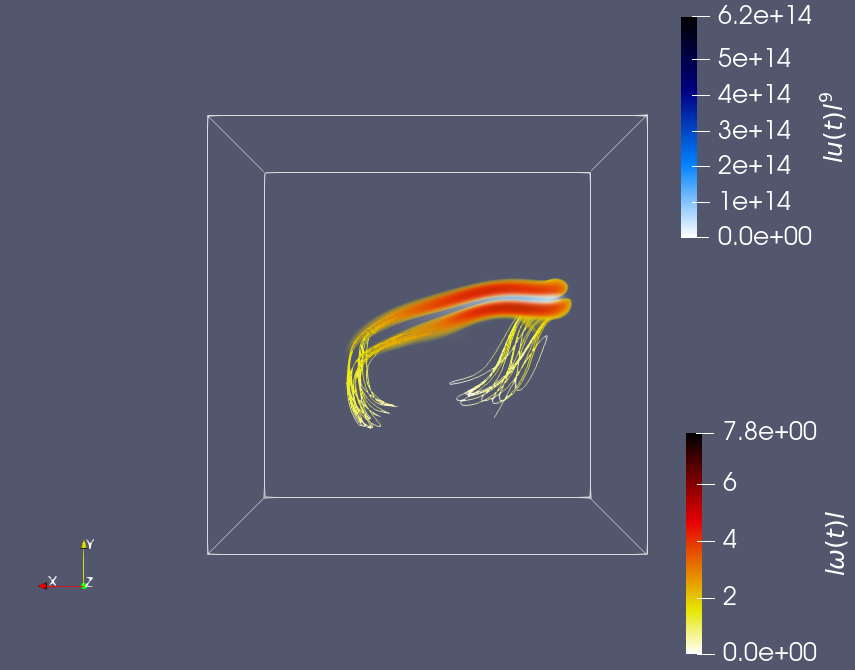}\label{L9SbL1e3T002_3}}
}
\caption{[$q=9$, Problem \ref{pb:PhiHs}, {$B=800$} and $T=2\times 10^{-4}$] 
Snapshots of the magnitude of the vorticity {$|\bomega(t_i,\x)|$} (red 
color scale) in the solution of the Navier-Stokes system \eqref{eq:NS} 
along with vortex lines (red) and the quantity {$|\u(t_i,\x)|^9$} 
(blue color scale) shown at the times $t_0,...,t_3$. These time 
instances are marked with solid black symbols in Figure \ref{fig:L9tb}.}
\label{fig:VisualizationTimeEvolSbq9}
\end{figure}

\begin{figure}[H]
\centering
\mbox{
\subfigure[${t_0}=0$]{\includegraphics[width=0.35\textwidth]{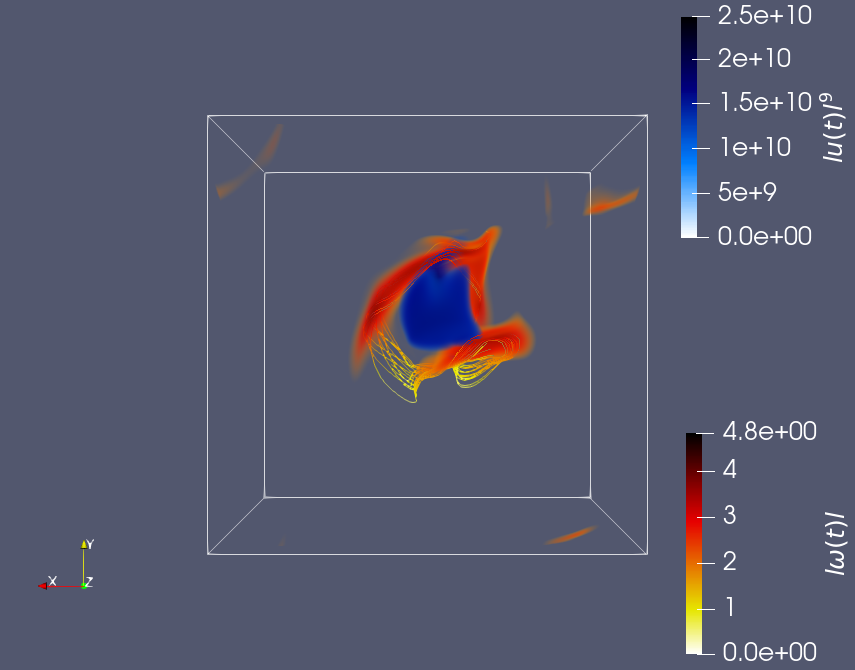}\label{L9LGL1e3T002_0}}
\quad
\subfigure[{${t_1}=1.2\times10^{-4}$}]{\includegraphics[width=0.35\textwidth]{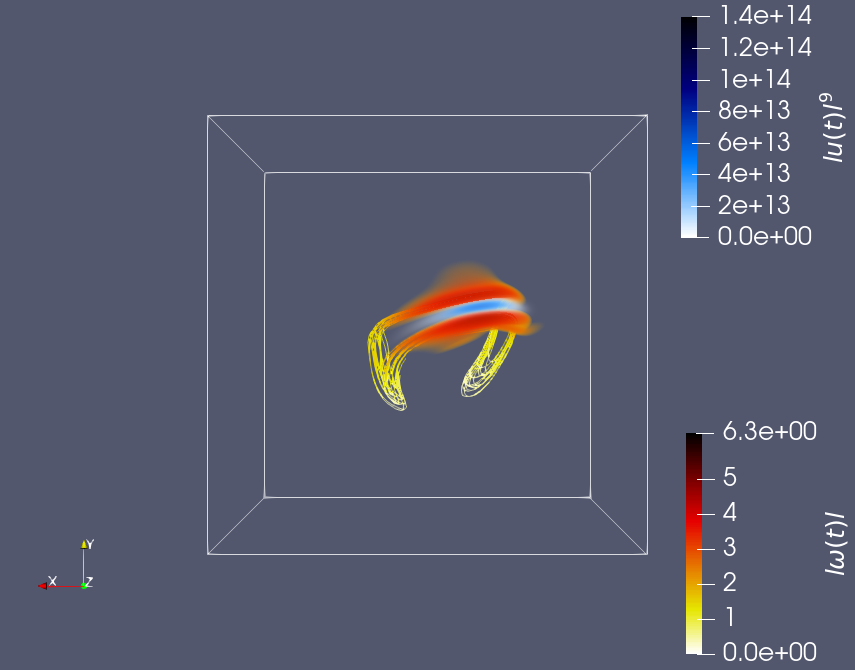}\label{L9LGL1e3T002_1}}
}
\mbox{
\subfigure[{${t_2}=1.5\times10^{-4}$}]{\includegraphics[width=0.35\textwidth]{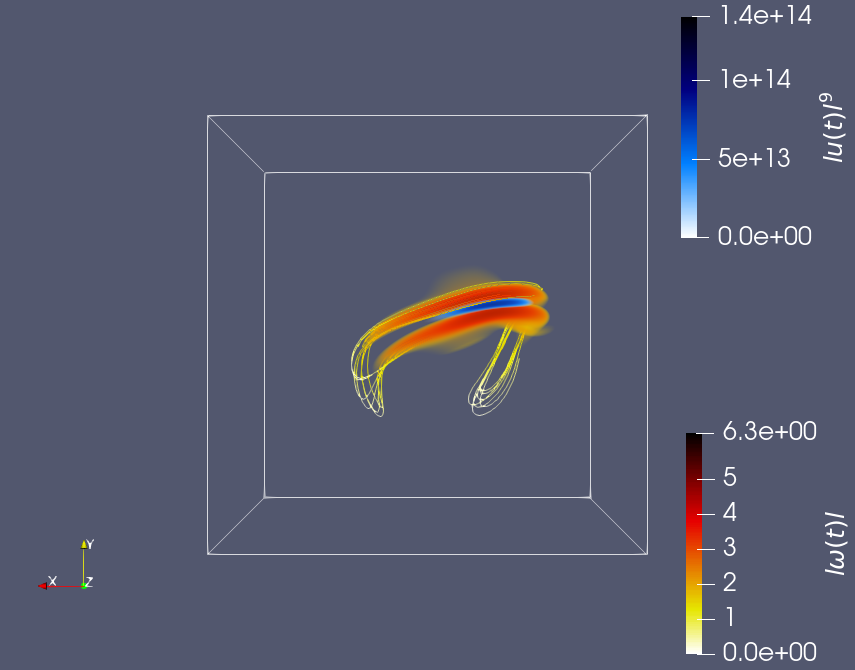}\label{L9LGL1e3T002_2}}
\quad
\subfigure[{${t_3}=2\times 10^{-4}$}]{\includegraphics[width=0.35\textwidth]{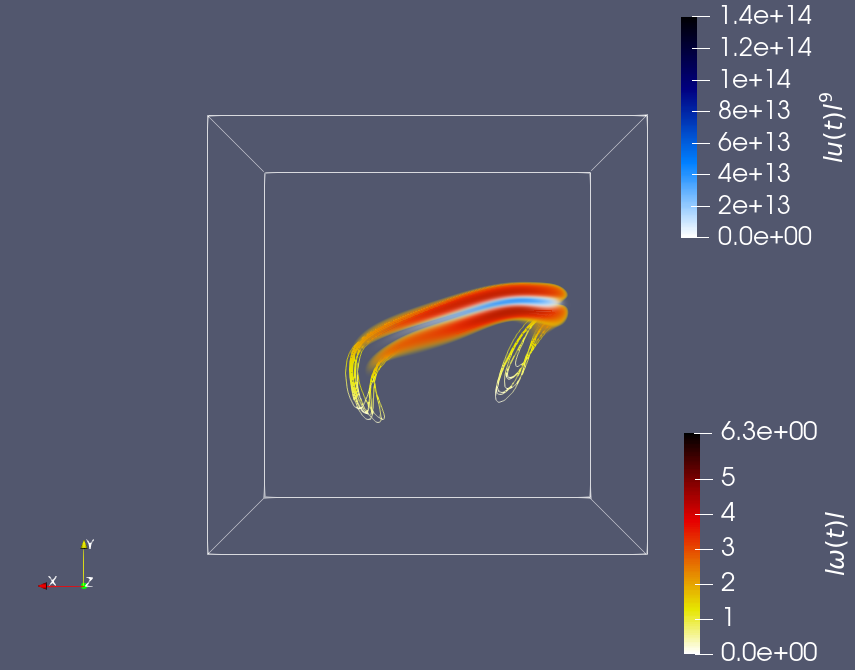}\label{L9LGL1e3T002_3}}
}
\caption{[$q=9$, Problem \ref{pb:PhiLq}, {$B=800$ and $T=2\times 10^{-4}$}] 
Snapshots of the magnitude of the vorticity {$|\bomega(t_i,\x)|$} (red 
color scale) in the solution of the Navier-Stokes system \eqref{eq:NS} 
along with vortex lines (red) and the quantity {$|\u(t_i,\x)|^9$} 
(blue color scale) shown at the times $t_0,...,t_3$. These time 
instances are marked with solid green symbols in Figure \ref{fig:L9tb}.}
\label{fig:VisualizationTimeEvolLGq9}
\end{figure}

Finally, we discuss the structure of the extremal flows belonging to 
the different branches shown in Figure \ref{fig:brL9a}. We do this by 
visualizing the vorticity magnitude $|\bomega(t_{i},\x)|$ and the 
quantity $|\u(t_{i},\x)|^9$ in space at the time instances $t_{i}$, 
$i=0,...,3$, defined in Section \ref{sec:vizL3}. These plots are shown 
in Figures \ref{fig:VisualizationTimeEvolSbq9} and 
\ref{fig:VisualizationTimeEvolLGq9} for flows obtained with the optimal 
initial conditions $\tuBT$ found by solving Problems \ref{pb:PhiHs} and 
\ref{pb:PhiLq} with the constraint parameter $B = 800$ and $T$ 
corresponding to the maxima over the branches shown in Figure 
\ref{fig:brL9a}. We observe that the flow structures corresponding to 
solutions of Problems \ref{pb:PhiHs} and \ref{pb:PhiLq} are quite 
similar, and they both feature flattened and bent vortex rings that 
stretch as the flows evolve. There are also similarities to the flow 
structures obtained for $q = 3$ in Section \ref{sec:L3}, cf.~Figures 
\ref{fig:VisualizationTimeEvolSbq3} and 
\ref{fig:VisualizationTimeEvolLGq3}.


\FloatBarrier
\subsection{Comparison of the Extreme Flows Obtained for Different Values of $q$}
\label{sec:Lqdiagn}

\begin{figure}[t]
\centering
\mbox{
\subfigure[]{\includegraphics[width=0.49\textwidth]{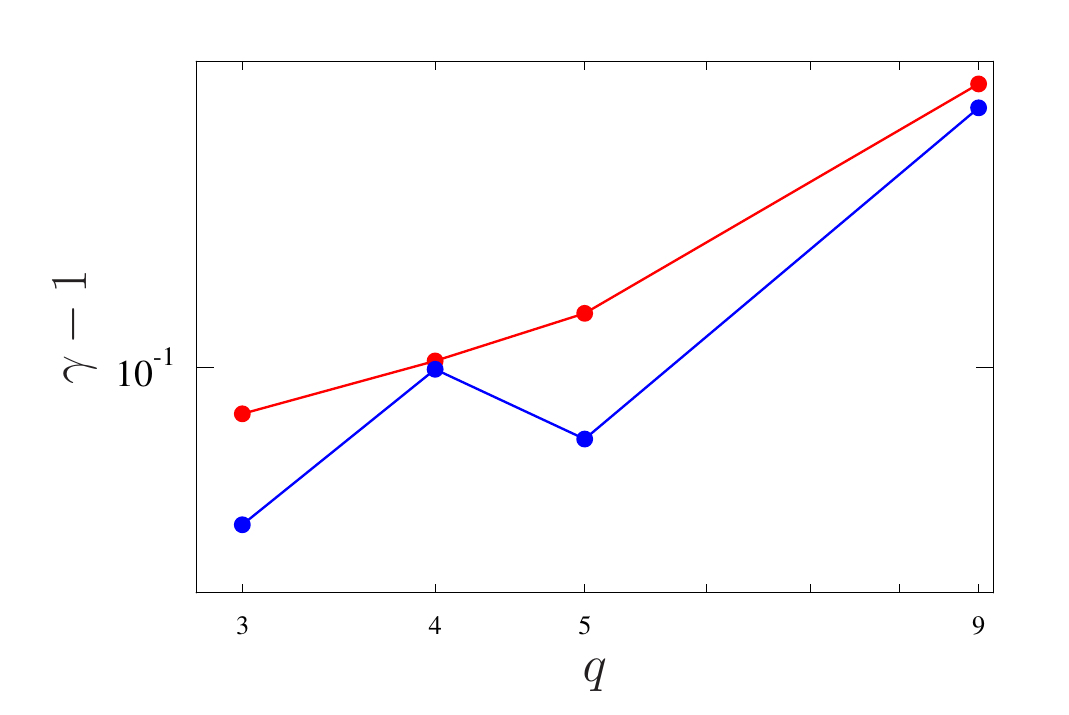}\label{fig:gamma}}
\quad
\subfigure[]{\includegraphics[width=0.49\textwidth]{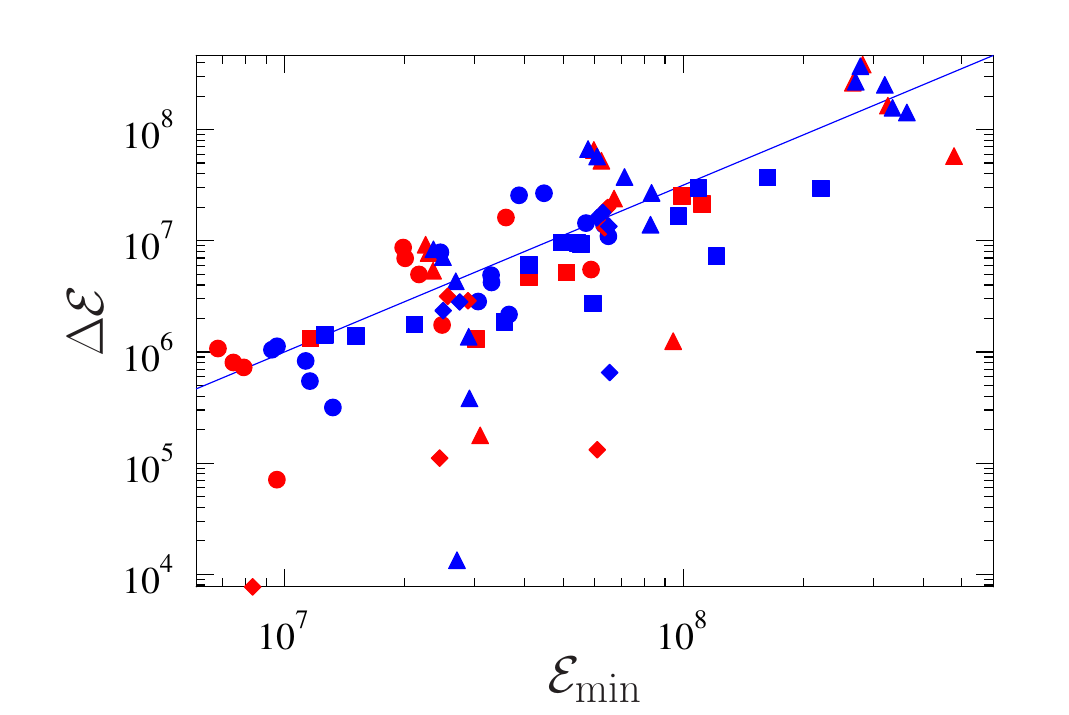}\label{fig:deltaE}}
}
\caption{\subref{fig:gamma} Dependence of the exponent $\gamma$ in 
expressions \eqref{eq:gamma} on  $q$  for flows with optimal initial 
conditions found by solving (blue) Problems \ref{pb:PhiHs} and 
\ref{pb:PsiHs} and (red) Problems \ref{pb:PhiLq} and \ref{pb:PsiLq}. 
\subref{fig:deltaE} Dependence of $\Delta\E$, cf.~\eqref{eq:deltaE}, on 
the minimum enstrophy  $\E_{\text{min}}$ in flows with the optimal 
initial conditions found by solving (blue) Problems \ref{pb:PhiHs} and 
\ref{pb:PsiHs}, and (red) Problems \ref{pb:PhiLq} and \ref{pb:PsiLq}. 
In \subref{fig:deltaE}, each symbol corresponds to a different value of 
$q$: (squares) $q=3$, (circles) $q=4$, (triangles) $q=5$ and (diamonds) 
$q=9$. The solid blue line represents the relation $\Delta \E \sim C 
\E_{\text{min}}^{3/2}$ for some $C>0$. }
\label{fig:compq}
\end{figure}

In this section we consolidate the results concerning the extreme flows 
obtained for different values of $q$ in order to identify certain 
general trends. In addition to the results for $q = 3$ and $q = 9$ 
already discussed above, we also include the data for $q = 4$ and $q = 
5$ which was obtained in a similar manner and for brevity is not 
presented here (see \cite{r2025} for all details). First, by comparing 
the plots in the left columns of Figures \ref{fig:L3t} and  
\ref{fig:L9t}, we observe that the relative, with respect to the 
constraint parameter $B$, growth of the norm $||\u(t)||_{L^q}$ is 
more pronounced for larger $q$ (this trend is also supported by the 
results for $q = 4,5$ \cite{r2025}). Another observation is that while 
for $q = 3,4$ solutions of optimization problems defined on the Sobolev 
spaces $H^{s}(\Omega)$ (Problems \ref{pb:PhiHs} and \ref{pb:PsiHs}) 
generally give rise to a larger growth of $||\u(t)||_{L^q}$ than 
solutions of problems defined on the Lebesgue spaces $L^{q}(\Omega)$ 
(Problems \ref{pb:PhiLq} and \ref{pb:PsiLq}), this growth is comparable 
in the two cases  for larger values of $q$. For a fixed $q$ this 
growth is quantified by the exponent $\gamma$ in the empirical 
relation, cf.~\eqref{eq:sqSb_q3}--\eqref{eq:sqLG_q3} and 
\eqref{eq:sqSb_q9}--\eqref{eq:sqLG_q9},
\begin{equation}
\max_{T}\,\Psi_{T}\left(\tuBT\right)\sim \left(B^{3}\right)^{\gamma}, \qquad 
\max_{T}\,\Phi_{T}^q\left(\tuBT\right)\sim \left(B^{p}\right)^{\gamma} \quad \text{for} \ q > 3. 
\label{eq:gamma}
\end{equation} 
The dependence of the  exponent $\gamma$ on $q$ is shown in Figure 
\ref{fig:gamma} where we plot $(\gamma - 1)$ versus $q$ to highlight 
the part of the growth of $||\u(t)||_{L^q}$ that can be attributed to 
the nonlinear amplification rather than merely to the increase of the 
constraint parameter $B$. We see that, interestingly, in solutions of 
Problems \ref{pb:PhiLq} and \ref{pb:PsiLq}, $(\gamma-1)$ is an 
increasing function of $q$. In solutions of Problems \ref{pb:PhiHs} and 
\ref{pb:PsiHs} the general trend is similar, but the increase with $q$ 
is not monotonic. The change of the trend at $q = 5$ could possibly be 
attributed to the solutions of Problem \ref{pb:PhiHs} at this value of 
$q$ being only suboptimal local maximizers.

We now turn to the behavior of the enstrophy in the extreme flows. 
Since in the flows considered here the enstrophy tends to initially 
decrease, cf.~Figures \ref{fig:L3t} and \ref{fig:L9t} (right columns), 
we define the quantity 
\begin{equation}
\Delta \E:=\max_{t\in\left[\argmin_{s\in[0,T]}\E(s),T\right]}\E(t)-{\E_{\text{min}}}, \qquad \text{where} \qquad \E_{\text{min}}:=\min_{t\in[0,T]}\E(t).
\label{eq:deltaE}
\end{equation}
measuring the total increase of the enstrophy relative to its minimum 
$\E_{\text{min}}$ which may be attained at an intermediate time (rather 
than at the initial time $t = 0$). It is plotted as a function of 
$\E_{\text{min}}$ for all extreme flows found by solving Problems 
\ref{pb:PhiHs}, \ref{pb:PhiLq}, \ref{pb:PsiHs} and \ref{pb:PsiLq} with 
$q = 3,4,5,9$ and three distinct values of $B$ in each case in Figure 
\ref{fig:deltaE}. The scatter evident in this plot is due to the fact 
that enstrophy is not directly  controlled in these problems, hence may in 
principle take arbitrary values. Despite this, we see that the envelope 
of all the data points exhibits a well-defined power-law relation 
$\Delta \E = \OO\left( \E_{\text{min}} ^{3/2} \right)$ consistent with 
what was already observed in \cite{KangYunProtas2020} when solving 
Problem \ref{pb:maxET}, cf.~\eqref{eq:maxEt}, and in 
\cite{KangProtas2021} when solving Problem  \ref{pb:PhiHs} for $q = 4$.


\section{Discussion and Conclusions}
\label{sec:final}

Here we summarize  the main findings of this study and draw some 
conclusions. Building on a body of earlier work 
\cite{ap11a,ap13a,ap13b,ap16,Yun2018,KangYunProtas2020} surveyed in 
\cite{p21a}, this investigation considers what is arguably the most 
fundamental family of conditional regularity results, the 
Ladyzhenskaya-Prodi-Serrin conditions \eqref{eq:LPS}--\eqref{eq:LPS3}. 
Singularity formation in a classical solution of the Navier-Stokes 
system \eqref{eq:NS} is signalled by the norm  $||\u(t)||_{L^{3}}$ and 
the integral in \eqref{eq:LPSblowup} becoming unbounded as the 
singularity time is approached. We systematically searched for initial 
data $\u_{0}$  that might lead to such a behavior by solving different 
PDE-constrained optimization problems for a broad range of parameters. 
While no evidence of unbounded growth of the quantities of interest, 
and hence also for singularity formation, was detected, we were able to 
quantify the transient behavior arising in such worst-case scenarios. 
This observation does not rule out the possibility of blow-up, since 
Problems \ref{pb:PhiHs}--\ref{pb:PsiLq} are nonconvex and therefore 
solutions found in Section \ref{sec:results} may only be (possibly 
nonunique) suboptimal local maximizers. In other words, it is possible 
there may be other, local or global, maximizers producing larger, 
possibly unbounded, growth of the quantities in question. Another 
possibility is that the values of the constraint parameter $B$ we used 
were too low such that the resulting flows would have been in the 
``small-data" regime where global existence of classical solutions can 
be asserted a priori. In this context we note that an analogous 
approach based on solution of variational optimization problems was 
recently applied to search for singularities in 3D inviscid Euler flows 
producing flows consistent with a finite-time blow-up \cite{zp23}.

A key methodological innovation of the present study is an approach to 
solve PDE optimization problems defined on general Banach spaces 
without Hilbert structure. Introduced in Section \ref{sec:gradLq}, it 
made it possible to solve Problems \ref{pb:PhiLq} and \ref{pb:PsiLq} 
where the optimal initial data $\tuBT$ was sought in the Lebesgue spaces 
$L^{q}(\Omega)$, $q \ge 3$. Solutions of optimization problems defined 
on these spaces and on Hilbert-Sobolev spaces $H^{s}(\Omega)$ embedded 
in them, cf.~\eqref{eq:SobEmbb}, are distinct which can be attributed 
to the different topologies these spaces are endowed with. In this 
regard, we found that local maximizers of Problems \ref{pb:PhiHs} and 
\ref{pb:PsiHs} are also local maximizers of Problems \ref{pb:PhiLq} 
and \ref{pb:PsiLq}, respectively. However, the converse is not true. 
This observation can be rigorously justified as the following theorem
\begin{theorem}
\label{thm:localmax}
Consider the following optimization problems defined on the Banach spaces $X$
and $Y$
\begin{align}
& \max_{\z\in X} \boldsymbol{\varphi}(\z) \qquad \text{and}  \tag{P1} \label{p1} \\
& \max_{\z\in Y} \boldsymbol{\varphi}(\z), \tag{P2}  \label{p2}
\end{align}
\noindent where $X$ is densely embedded in $Y$ and 
$\boldsymbol{\varphi}(\z)$ is a Fr\'echet differentiable objective functional 
defined on both $X$ and $Y$.
If $\z_{0}$ is a local solution of \eqref{p1},
then $\z_{0}$ is also a local solution of \eqref{p2}.
\end{theorem}
\noindent Although Problems \ref{pb:PhiHs}--\ref{pb:PsiLq}
are constrained, they can be converted to an unconstrained form
using Lagrange multipliers such that one can invoke Theorem \ref{thm:localmax}.
The proof of this theorem is given in Appendix \ref{a2}.

While  for smaller values of $q$ solutions of Problems \ref{pb:PhiHs} 
and \ref{pb:PsiHs} defined on the Sobolev spaces generally yield larger 
values of the objective functional than solution of Problems 
\ref{pb:PhiLq} and \ref{pb:PsiLq} defined on Lebesgue spaces, these 
values are comparable for $q = 9$. The first observation is 
counterintuitive  since Lebesgue spaces can be viewed as ``larger" than 
the Sobolev spaces embedded in them.

In regard to the transient behavior observed in the extreme flows, 
Figures \ref{fig:RoCL3} and \ref{fig:RoCL9} together with the estimates 
on the maximum rate of growth of the norms $||\u(t)||_{L^q}$ and of the 
enstrophy $\E(t)$ from Sections \ref{sec:dEdt} and \ref{sec:dLPSdt} 
indicate that the extreme flows do enter a regime where these 
quantities are amplified at a rate consistent with singularity 
formation in finite time, but this growth is not sustained long enough 
for singularities to form due to depletion of nonlinearity. In the 
absence of singularities, quantities \eqref{eq:Psi}--\eqref{eq:Phi} 
appear bounded for all values of $B$ and $T$, and were found to scale 
in proportion to $(B^{p})^{\gamma}$, where $\gamma$ is an increasing 
function of $q$, cf.~Figure \ref{fig:gamma}. On the other hand, the 
maximum enstrophy in the extreme flows  scales as $\OO\left( 
\E_{\text{min}} ^{3/2} \right)$, cf.~Figure \ref{fig:deltaE}, which is 
the same behavior as was obtained in \cite{KangYunProtas2020} by 
solving Problem \ref{pb:maxET}. The fact that the extreme flows 
obtained by solving Problem \ref{pb:maxET} and Problems 
\ref{pb:PhiHs}--\ref{pb:PsiLq} exhibit the same scaling of the maximum 
enstrophy is intriguing since this quantity is not directly controlled 
in the latter case and this growth is realized through entirely 
different physical mechanisms. We add that the same bound is now 
rigorously established for the growth of enstrophy in 1D viscous 
Burgers flows \cite{ap11a,AlbrittonDeNitti2023}. Perhaps, if one day 
the 3D Navier-Stokes system \eqref{eq:NS} is shown to be globally well 
posed, the a priori bound on the maximum growth of enstrophy might have 
this form.

A comparison of the results from Sections \ref{sec:L3}--\ref{sec:L9} 
and Figure \ref{fig:gamma} indicates that solutions of Problems 
\ref{pb:PhiHs} and \ref{pb:PhiLq} obtained for increasing values of $q$ 
lead to a more ``extreme" behavior. Therefore, an interesting open 
question that we intend to consider in the future is what happens for 
larger values of $q > 9$.

As regards the numerical approach, cf.~Appendix \ref{sec:numer}, while 
pseudo-spectral methods provide excellent accuracy in approximating 
PDEs with smooth solutions, their main shortcoming is that they do not 
allow for local mesh refinement necessary to resolve very fine flow 
structures underlying extreme flow behaviors. Overcoming this 
limitation is necessary in order to be able to solve Problems 
\ref{pb:maxET}--\ref{pb:PsiLq} for larger values of the constraint 
parameter $B$ or, equivalently, larger values of the Reynolds number 
\eqref{eq:Re} where other amplification mechanisms may become active. 
To address this problem, in the future we will employ discretization 
techniques, such as finite-difference or finite-element methods, which 
do allow for local mesh refinement, especially adaptive mesh refinement.

\section*{Acknowledgements}

The authors thank Fabian Bleitner, Alexey Cheskidov, Mimi Dai, Roman 
Shvydkoy and Xinyu Zhao for many enlightening and enjoyable 
discussions. Partial funding for this research was provided through an 
NSERC (Canada) Discovery Grant RGPIN-2020-05710. Computational 
resources were provided by the Digital Research Alliance of Canada 
under its Resource Allocation Competition .

\begin{appendices}

\section{Evaluation of the Gradient in $L^2$}
\label{sec:L2grad}

When viewed as a function of its second argument, the G\^{a}teaux 
differential $(\Phi^q_T)'(\u_0;\cdot)$ is a bounded linear functional 
on $L^{2}(\Omega)$ which allows us to use the Riesz identity 
\eqref{eq:rieszL2} to extract the $L^{2}$ gradient from it, 
cf.~\eqref{eq:gradL2}. We can use its definition \eqref{eq:dPhi} to 
evaluate the G\^{a}teaux differential of the objective functional 
\eqref{eq:Phi} as
\begin{equation}
(\Phi^q_T)'(\u_0;\u_0') = \frac{2q}{(q-3)T} \int_0^T \left( \|\u(t)\|_{L^q}^{\frac{q(5-q)}{q-3}} 
\int_{\Omega} |\u(t,\x)|^{q-2} \u(t,\x) \cdot \u'(t,\x) d\x \right) dt,
\label{eq:dPhiT}
\end{equation}
where the perturbation field $\u' = \u'(t,\x)$ is a solution of the
Navier-Stokes system linearized around the trajectory $\u$ corresponding to
the initial data $\u_0$ \cite{g03}, i.e.,
\begin{subequations}
\label{eq:lNSE3D}
\begin{align}
 \mathcal{L}\begin{bmatrix} \u' \\ p' \end{bmatrix} := 
& \begin{bmatrix}
\partial_{t}\u'+\u'\cdot\bnabla_{\x}\,\u+\u\cdot\bnabla_{\x}\,\u'+\bnabla_{\x}\, p'-\nu\laplacian\u' \\
\bnabla_{\x}\cdot\u'
\end{bmatrix} = \begin{bmatrix} \mathbf{0} \\ 0\end{bmatrix}, \label{eq:lNSE3Da} \\
 \u'(0)= &\u_0' \label{eq:lNSE3Db}
\end{align}
\end{subequations}
which is subject to the periodic boundary conditions and where $p'$ is 
the perturbation to the pressure. We note that expression 
\eqref{eq:dPhiT} for the G\^{a}teaux differential is not yet in the Riesz 
form \eqref{eq:rieszL2}, because the perturbation $\u_0'$ of the 
initial data does not enter in it explicitly as a factor, but instead 
appears in the initial condition \eqref{eq:lNSE3Db} of the perturbation 
system. In order to transform the G\^{a}teaux differential to the 
required Riesz form \eqref{eq:rieszL2}, we introduce the {\em adjoint 
states} $\u^* \, : \, [0,T]\times\Omega  \rightarrow \RR^3$ and $p^* \, 
: \, [0,T]\times\Omega  \rightarrow \RR$, and the following 
duality-pairing relation
\begin{equation}
\begin{aligned}
\left\langle \mathcal{L}\begin{bmatrix} \u' \\ p' \end{bmatrix}, \begin{bmatrix} \u^* \\ p^* \end{bmatrix} \right\rangle
:= & \int_0^T \int_{\Omega} \mathcal{L}\begin{bmatrix} \u' \\ p' \end{bmatrix} \cdot \begin{bmatrix} \u^* \\ p^* \end{bmatrix} \, d\x \, dt 
= \overbrace{\left\langle \begin{bmatrix} \u' \\ p' \end{bmatrix}, \mathcal{L}^*\begin{bmatrix} \u^* \\ p^* \end{bmatrix}\right\rangle}^{(\Phi^q_T)'(\u_0;\u_0')} + \\
\phantom{=} & \int_\Omega \u'(t,\x)\cdot\u^*(T,\x)  \,d\x - 
\int_\Omega \u'(0,\x)\cdot\u^*(0,\x)  \,d\x = 0,
\end{aligned}
\label{eq:dual}
\end{equation}
where ``$\cdot$'' in the first integrand expression denotes the 
Euclidean dot product evaluated at $(t,\x)$. Using \eqref{eq:lNSE3D}, 
performing integration by parts with respect to both space and time 
(where all boundary terms resulting from integration by parts with 
respect to space vanish due to periodicity) and defining the adjoint 
system as in \eqref{eq:aNSE3D} with a judicious choice of the source 
term $\f$ and the terminal condition $\u^{*}(T)$ leaves us with the 
identity
\begin{equation}
(\Phi^q_T)'(\u_0;\u_0') = \int_\Omega \u'_0(\x)\cdot\u^*(0,\x)  \,d\x=\Big\langle \u^*(0), \u'_{0} \Big\rangle_{L^2}
\label{eq:dET2}
\end{equation}
which is a Riesz representation of the G\^{a}teaux differential. Noting 
that the perturbation $\u_{0}'$ is arbitrary, expression 
\eqref{eq:gradL2} for the $L^{2}$ gradient is finally obtained.

\section{Numerical Implementation}
\label{sec:numer}

Here we provide some information about the numerical approaches used to 
implement different elements of Algorithm \ref{alg:optimAlg}. First, we 
describe the techniques employed to discretize the the governing and 
adjoint system \eqref{eq:NS} and \eqref{eq:aNSE3D}, followed by a 
discussion of the computation of the Lebesgue gradients and solution of 
the arc-search problem \eqref{eq:tau_nHS}. We refer the reader to 
\cite{r2025} for additional details and validation tests.

\subsection{Discretization of PDE Problems \eqref{eq:NS} and \eqref{eq:aNSE3D}}

Discretization of systems \eqref{eq:NS} and \eqref{eq:aNSE3D} in space
is based on the standard pseudospectral Fourier-Galerkin method with 
dealiasing \cite{canuto:SpecMthd,boyd2001chebyshev}. Suppose 
$\u:[0,T]\times\Omega\rightarrow\RR^3$ is a vector field and consider 
the set $\mathcal{W}_{N}=\left\{\k\in\ZZ^3:\,k\leq N\right\}$ where 
$N\in\mathbb{N}^{+} $ is the spatial resolution. The Galerkin 
approximation of $\u(t,\x)$ is defined as
\begin{equation}
\u_{N}(t,\x):=\sum_{\k\in\mathcal{W}_{N}} \widehat{\u}_{\k}(t)e^{2\pi i \k\cdot\x},
\end{equation}
where $\widehat{\u}_{\k}(t)\in\CC^3$ are the Fourier coefficients. For 
$\u_{N}$ to be an approximate solution of \eqref{eq:NS}, its Fourier 
coefficients must  satisfy the following finite-dimensional system of 
ordinary differential equations (ODEs) obtained by performing a 
Galerkin projection of system \eqref{eq:NS} onto the subspace spanned 
by the Fourier modes with wavenumber{s} lying in $\mathcal{W}_{N}$
\begin{subequations}\label{eq:NSfourier}
\begin{alignat}{2}
\frac{d\widehat{\u}_{\k}(t)}{dt}+\boldsymbol{A}\,\widehat{\u}_{\k}(t)+\boldsymbol{r}(\widehat{\u}_{\k})(t)&=\0, &\\ 
\widehat{\u}_{\k}(t)\cdot\k&=0, \qquad &\mbox{for all}\;\k\in\mathcal{W}_{N},\,t\in(0,T),\label{eq:NSfourierb} \\
\widehat{\u}_{\k}(0)&=\left[\widehat{\u_{0}}\right]_{\k}, &
\end{alignat}
\end{subequations}
where the linear and the nonlinear operators $\boldsymbol{A}$ and $\boldsymbol{r}$ are defined as
\begin{align}
\boldsymbol{A}\,\widehat{\u}_{\k}(t)&:=(2\pi)^2\,k^2\,\widehat{\u}_{\k}(t), \\
\boldsymbol{r}(\widehat{\u}_{\k})(t)&:=\left[\widehat{\left(\u\cdot\bnabla_{\x}\right)\u}\right]_{\k}-2\pi\,i\,\k\,\widehat{p}_{\k}, \label{eq:op_r}
\end{align}
respectively. The Fourier coefficients of the nonlinear term 
$\left[\widehat{\left(\u\cdot{\bnabla_{\x}}\right)\u}\right]_{\k}$
are evaluated by {first computing} the product 
$\left(\u\cdot{\bnabla_{\x}}\right)\u$ in the physical space  and then 
calculating its Fourier transform with dealiasing. The Fourier 
coefficients of the pressure term, $\widehat{p}_{\k}$, are computed by 
solving the Poisson equation $\Delta p=-{\bnabla_{\x}}\,\cdot\left(\left(\u\cdot{\bnabla_{\x}}\right)\u\right)$ 
 subject to periodic boundary conditions, which is also performed in the Fourier space as
 \begin{equation}
 \widehat{p}_{\k}={\frac{\k}{k^2}\cdot\left[\widehat{\left(\u\cdot\bnabla_{\x}\right)\u}\right]_{\k}.}
 \end{equation}

For the discretization of the adjoint system \eqref{eq:aNSE3D}, the 
Fourier coefficients of the solution satisfy the following system of 
ODEs 
\begin{subequations}\label{eq:aNSfourier}
\begin{align}
\frac{d\widehat{\u}^{*}_{\k}(t)}{dt}+{\boldsymbol{S}\,\widehat{\u}^{*}_{\k}(t)}+{\boldsymbol{R}\,\widehat{\u}^{*}_{\k}(t)}&=\0, &\\
\widehat{\u}^{*}_{\k}(t)\cdot\k&=0, \qquad &\mbox{for all}\;\k\in\mathcal{W}_{N},\,t\in(0,T),\label {eq:aNSfourierb} \\
\widehat{\u}^{*}_{\k}(T)&=\0 &
\end{align}
\end{subequations}
where the linear operators {$\boldsymbol{S}$ and $\boldsymbol{R}$} are defined as
\begin{align}
\boldsymbol{S}\,\widehat{\u}^{*}_{\k}(t)&:=\,(2\pi)^2\,k\,\widehat{\u}^{*}_{\k}(t),\\
\boldsymbol{R}\,\widehat{\u}_{\k}(t)&:=\,\left[\widehat{({\bnabla_{\x}}\,\u^*)\,\u}\right]_{\k}(t)
+\left[\widehat{({\bnabla_{\x}}\,\u^*)^T\,\u}\right]_{\k}(t)-\widehat{\f}_{\k}(t)
\end{align}
and $\widehat{\f}_{\k}(t)$ are the Fourier coefficients of the source term 
$\f(t,\x)$ in \eqref{eq:aNSE3Df}. The terms 
$\widehat{({\bnabla_{\x}}\,\u^*)\,\u}$ and 
$\widehat{({\bnabla_{\x}}\,\u^*)^T\,\u}$ are first evaluated in 
physical space as in \eqref{eq:op_r}, and then transformed (with 
dealiasing) to the Fourier space. Unlike system \eqref{eq:NSfourier}, 
system \eqref{eq:aNSfourier} is a terminal-value problem and needs 
to be integrated backwards in time. The coefficients and the source 
term in the adjoint system \eqref{eq:aNSE3D} are given in terms of the 
solution $\u = \u(t,\x)$ of the Navier-Stokes system \eqref{eq:NS} 
around which linearization is performed at each iteration is Algorithm 
\ref{alg:optimAlg}. This solution therefore needs to be saved at all 
time steps discretizing the time interval $[0,T]$.

The ODE systems \eqref{eq:NSfourier} and \eqref{eq:aNSfourier} are 
integrated in time using a hybrid approach combining a second-order 
implicit Crank–Nicolson method applied to the linear terms and a 
third-order explicit Runge-Kutta method applied to the nonlinear terms 
\cite{NumRenaissance}, which offers an optimal trade-off between 
accuracy, stability and the computational cost. The incompressibility 
conditions \eqref{eq:NSfourierb} and \eqref{eq:aNSfourierb} are applied 
at each substep in the process. The implementation is massively 
parallel and relies in on the Message-Passing Interface (MPI) and the 
library {\tt FFTW} to perform parallel Fast Fourier Transforms 
\cite{fftw}.

\subsection{Computation of the Lebesgue Gradients}

To find a solution of  system \eqref{eq:LGsystem}, we employ an 
iterative splitting method proposed in \cite{pbh04}, where at each 
iteration we first  apply Newton’s method with globalization to 
equation solve \eqref{eq:LGsystema} with  $\eta$ fixed and then we 
apply a standard Poisson solver to equation \eqref{eq:LGsystemb} to 
update $\eta$.{}

\subsection{Solution of Arc-Search Problems}

The arc-search problems \eqref{eq:tau_nHS}--\eqref{eq:tau_nLq} are 
solved using a variant of derivative-free Brent's algorithm combining 
the golden-section search with inverse quadratic interpolation 
\cite{nw00,press2007numerical}.

\section{Proof of Theorem \ref{thm:localmax}}
\label{a2}

In this appendix we prove Theorem \ref{thm:localmax}. Suppose that $\z_{0}\in X$ is a local solution of \eqref{p1}. Therefore, it satisfies the optimality condition
\begin{equation}
\label{eq:opticond}
\boldsymbol{\varphi}'(\z_{0},\boldsymbol{w}')=\lim_{\epsilon\rightarrow0}\frac{\boldsymbol{\varphi}(\z_{0}+\epsilon\boldsymbol{w}')-\boldsymbol{\varphi}(\z_{0})}{\epsilon}=0,\quad\text{for all }\,\boldsymbol{w}'\in X.
\end{equation} 

\noindent We wish to verify that the following optimality condition  holds as well
\begin{equation}
\boldsymbol{\varphi}'(\z_{0},\z')=\lim_{\epsilon\rightarrow0}\frac{\boldsymbol{\varphi}(\z_{0}+\epsilon\z')-\boldsymbol{\varphi}(\z_{0})}\epsilon=0, \quad \text{for all } \z'\in Y.
\end{equation}

\noindent Given that $X$ is dense in $Y$, every $\z'{\in Y}$ can be approximated with a sequence of 
elements in $X$, i.e., there exists a sequence $\{\z'_{n}\}_{n\in\mathbb{N}} \in X$ such that $$\lim_{n\rightarrow\infty}\z'_{n}=\z'.$$ 
In particular, we have that for any $n\in\mathbb{N}$ expression \eqref{eq:opticond} implies
\begin{align*}
0=&\lim_{n\rightarrow\infty}\lim_{\epsilon\rightarrow0}\frac{\boldsymbol{\varphi}(\z_{0}+\epsilon\z_n')-\boldsymbol{\varphi}(\z_{0})}{\epsilon}\\
=&\lim_{\epsilon\rightarrow0}\lim_{n\rightarrow\infty}\frac{\boldsymbol{\varphi}(\z_{0}+\epsilon\z_n')-\boldsymbol{\varphi}(\z_{0})}{\epsilon}\\
=&\lim_{\epsilon\rightarrow0}\frac{\boldsymbol{\varphi}(\z_{0}+\epsilon\lim_{n\rightarrow\infty}\z'_{n})-\boldsymbol{\varphi}(\z_{0})}{\epsilon}\quad\text{(By the continuity of $\boldsymbol{\varphi}$)}\\
=&\lim_{\epsilon\rightarrow0}\frac{\boldsymbol{\varphi}(\z_{0}+\epsilon\z')-\boldsymbol{\varphi}(\z_{0})}{\epsilon}\\
=&\boldsymbol{\varphi}'(\z_{0},\z').
\end{align*}

\noindent Hence, $\z_{0}$ is also a local solution of Problem \eqref{p2}.\\ \qed

%
%
%
%
%
%
\end{appendices}


\bigskip
\noindent
{\bf Declaration of Interests.} The authors report no conflict of interest.

\bigskip
\noindent
{\bf Data Availability.} Data from the paper will be made available upon request.





\end{document}